\documentclass[11pt]{article}

\usepackage[a4paper,margin=1in]{geometry}
\usepackage{amsmath,amssymb,amsthm}
\usepackage{fancyvrb}
\usepackage{bbm}
\usepackage{stmaryrd}
\usepackage{changepage}
\usepackage{verbatim}
\usepackage[scr=boondoxo]{mathalfa}
\usepackage{nicematrix}
\usepackage{setspace}
\usepackage{multirow}
\usepackage{float}
\usepackage{comment}
\usepackage{enumitem}
\usepackage[all]{xy}
\usepackage{shuffle}
\usepackage{pdfrender}
\usepackage{xcolor}
\usepackage{csquotes}
\usepackage{scalerel}
\usepackage{accents}
\usepackage{subcaption}
\usepackage[hidelinks]{hyperref}

\newcommand{\boldchar}[1]{%
  \textpdfrender{TextRenderingMode=FillStroke,LineWidth=.6pt}
  {\ensuremath{#1}}%
}

\title{Nonlinear Stochastic Filtering with Volterra Gaussian Noises}

\author{
Thomas Cass\thanks{
Imperial College London, United Kingdom.
\texttt{thomas.cass@imperial.ac.uk}}
\and
Dan Crisan\thanks{
Imperial College London, United Kingdom.
\texttt{d.crisan@imperial.ac.uk}}
\and
Andrea Iannucci\thanks{
Imperial College London, United Kingdom.
\texttt{andrea.iannucci22@imperial.ac.uk}.
Corresponding author.}
}

\date{}

\newcommand{\floor}[1]{\left\lfloor #1\right\rfloor}

\usepackage{amsthm}

\newtheorem{theorem}{Theorem}[section]
\newtheorem{proposition}[theorem]{Proposition}
\newtheorem{lemma}[theorem]{Lemma}
\newtheorem{corollary}[theorem]{Corollary}

\theoremstyle{definition}
\newtheorem{definition}[theorem]{Definition}

\newtheorem{example}[theorem]{Example}

\theoremstyle{remark}
\newtheorem{remark}[theorem]{Remark}

\theoremstyle{definition}
\newtheorem{condition}[theorem]{Condition}

\theoremstyle{remark}

\allowdisplaybreaks

\begin{document}

\maketitle

\begin{abstract}
We develop a nonlinear filtering theory for signal--observation systems driven
by Volterra Gaussian processes, covering both the Young and genuinely rough
regimes. The dynamics are formulated as a rough differential equation in
which the observation has a signal-dependent Volterra drift, a structure
naturally induced by an equivalent change of measure. We establish global
well-posedness of the coupled system and derive a
Kallianpur--Striebel formula. We then obtain a robust pathwise representation
of the filter. In the one-dimensional setting, we characterise the
unnormalised conditional density through a rough Zakai equation and establish
its well-posedness using an extension of the rough-viscosity framework.
Finally, under a partial H\"ormander-type condition, we prove that the
conditional distribution of the signal admits a smooth density.
\end{abstract}



\section{Introduction}
Nonlinear filtering concerns the estimation of the conditional distribution of an unobserved
signal given the information generated by a noisy observation process. More
precisely, let \((X,Y)\) be a signal--observation pair defined on a probability
space \((\Omega,\mathcal F,\mathbb P)\). For each \(t\geq 0\), the information
generated by the observation process up to time \(t\) is represented by the
completed filtration
\[
    \mathcal Y_t
    :=
    \sigma\{Y_s:0\leq s\leq t\}
    \vee
    \mathscr Z,
\]
where \(\mathscr Z\) denotes the collection of all subsets of
\(\mathbb P\)-null sets in \(\mathcal F\). The filter at time \(t\) is the
conditional distribution of \(X_t\) given \(\mathcal Y_t\), characterised
through its action on bounded measurable functions by
\[
    \pi_t(\phi)
    :=
    \mathbb E^{\mathbb P}
    \left[
        \phi(X_t)
        \,\middle|\,
        \mathcal Y_t
    \right],
\]
where \(\phi:\mathbb R^{d_X}\to\mathbb R\) is bounded and measurable.
In the classical nonlinear filtering framework, the signal and observation
are typically specified as the solution of a system of It\^o stochastic
differential equations driven by independent Brownian motions. In particular,
\((X,Y)\) is a continuous semimartingale, so that the filtering problem can be
studied using It\^o calculus, Girsanov's theorem and martingale representation
methods.

Brownian motion is, however, unsuitable for modelling noise with temporal
dependence. A natural alternative is provided by Volterra Gaussian processes,
namely Gaussian processes admitting a representation of the form
\[
    I_t
    =
    \int_0^t K(t,s)\,dV_s,
\]
where \(V\) is a Brownian motion and \(K\) is a deterministic Volterra kernel.
This class includes fractional Brownian motion, Riemann--Liouville fractional
Brownian motion and a broad range of Gaussian processes with non-stationary
increments and memory effects.

Volterra Gaussian processes are generally neither
Markov processes nor semimartingales. Consequently, the It\^o calculus,
martingale representation arguments and semimartingale change-of-measure
techniques underlying the classical filtering theory are no longer directly
available. When the sample paths are sufficiently regular, the dynamics may still be interpreted using Young integration.
Much of the existing filtering literature for fractional Brownian noise is
concentrated in this regular regime.

In the linear setting, Kleptsyna, Kloeden and Anh studied filtering models
involving fractional Brownian motion, first under the restriction
\(H>3/4\) in \cite{kleptsyna1998linear}, and subsequently for \(H>1/2\) in
\cite{kleptsyna1999linear}. Kleptsyna, Le Breton and Roubaud developed a more
general approach to filtering systems with fractional Brownian observation
noise and applied it to linear models in \cite{ml2000general}.

The ideas in \cite{ml2000general} also form the basis of work by Mandrekar and Gawarecki \cite{Gawarecki2001}, who use them to study nonlinear filtering when the signal solves a martingale problem and the observation noise is a fractional Brownian motion of arbitary path regularity.

Further developments in the regime \(H>1/2\) include the work of Amirdjanova
and Chivoret \cite{amirdjanova2006new}. For a nonlinear filtering problem with
additive fractional Brownian observation noise, they represent the optimal
filter through multiple stochastic fractional integral expansions. Their
results provide a fractional counterpart of the classical multiple stochastic
integral expansions obtained by Ocone and Pardoux in
\cite{ocone1983multiple}.

A model particularly relevant to the present work was introduced by Coutin and
Decreusefond in \cite{coutin1999abstract}. Their
signal-observation system takes the schematic form
\begin{align}
    X_t
    &=
    x_0
    +
    \int_0^t K(t,s)\gamma(X_s)\,ds
    +
    \int_0^t K(t,s)\alpha(X_s)\,d\beta_s,
    \label{equation_coutin_decreusefond_signal}
    \\
    Y_t
    &=
    \int_0^t K(t,s)h(X_s)\,ds
    +
    \int_0^t K(t,s)\tau(Y_s)\,dW_s,
    \label{equation_coutin_decreusefond_observation}
\end{align}
where $K$ is the kernel associated to the fractional Brownian motion with Hurst parameter $H \geq 1/2$ and \(\beta\) and \(W\) are Brownian motions and the stochastic integrals are
understood in the It\^o sense. Thus, the system is formulated as a
stochastic Volterra equation driven by the Brownian motions underlying the
fractional noises.

The approach of Coutin and Decreusefond, together with the methods employed in
the works discussed above, does not extend directly to driving paths of
Brownian or lower regularity. Indeed, in this regime the processes are not semimartingales and Young integrability
condition generally fails.
Rough path theory provides the corresponding framework that allows to formulate a 
notion of solution of differential equation driven by signals of such low regularity.

Motivated by the Volterra structure of
\eqref{equation_coutin_decreusefond_observation}, we consider a
signal-observation pair satisfying
\begin{equation}
\label{equation_introduction_system}
    \begin{pmatrix}
        X_t\\
        Y_t
    \end{pmatrix}
    =
    \begin{pmatrix}
        x_0\\
        0
    \end{pmatrix}
    +
    \begin{pmatrix}
        0\\
        Kb(X, Y)(t)
    \end{pmatrix}
    +
    \int_0^t
    \begin{pmatrix}
        \sigma(X, Y)_s & 0\\
        0 & \mathbbm 1
    \end{pmatrix}
    d\mathbf I_s,
\end{equation}
with
\(
    Kb(X, Y)(t)
    :=
    \int_0^t K(t,s)b(X_s, Y_s)\,ds
\) and $\mathbf I$ the rough path constructed from $I = (B, U)$ where $B$ and $U$ are respectively the signal and the observation noise.

Our observation equation retains precisely the structural compatibility used
by Coutin and Decreusefond. Indeed, when \(I\) is a fractional Brownian
motion and \(b=h\), \eqref{equation_introduction_system} coincides with
the additive-noise specialisation of
\eqref{equation_coutin_decreusefond_observation}. The signal equations are,
however, of a different nature: whereas
\eqref{equation_coutin_decreusefond_signal} is a stochastic Volterra equation
defined through It\^o integration against an underlying Brownian motion,
\eqref{equation_introduction_system} is a rough differential equation driven
directly by the Volterra Gaussian process. Our model should therefore be viewed as a rough-path extension of the
additive-observation framework of Coutin and Decreusefond: it preserves its
central Volterra structure while allowing for general Volterra Gaussian
drivers, covering both the Young regime and lower sample-path regularities
beyond it.

More generally, the observation equation in
\eqref{equation_introduction_system} is naturally suited to situations in
which the observed process becomes a centred Volterra Gaussian process under
an equivalent reference measure. Indeed, the corresponding change of measure
acts on the Brownian motion underlying the Volterra representation and
therefore induces, under the original measure, a drift transported through
the kernel \(K\). The resulting observation dynamics consequently take the
form \(U+Kb(X)\), with the latent state \(X\) affecting the observation through the Volterra drift term \(Kb\).

A concrete example of this arises from rough-volatility modelling.
One may regard \(Y\) as an observable volatility factor extracted from market
data, such as an implied-volatility or forward-variance factor, and \(X\) as an
unobserved state governing the evolution of volatility or volatility of
volatility. Under a suitable reference measure, the observable factor may be
described by a fractional Brownian motion, or more generally by a Volterra
Gaussian process (\cite{Livieri02092018} and \cite{floch2022roughnessimpliedvolatility}). Returning to the original measure introduces the drift
\(Kb(X)\), through which the latent state affects the observed volatility
factor. At the same time, the dependence of the signal coefficient on \(Y\)
allows the hidden state to evolve in response to the observed rough-volatility
environment. The framework therefore accommodates a two-way interaction
between an observable rough factor and an unobservable state variable.

The remainder of the paper is organised as follows. Section~2 introduces the
analytic framework used throughout the paper. We recall the elements of rough
path theory needed to formulate differential equations driven by signals of
low regularity, together with the Gaussian analytic machinery used to construct
rough path of the signal and observation noises. We then introduce the class of Volterra Gaussian processes
considered in this work and describe the covariance and kernel assumptions
under which it is possible to construct a rough path above these processes.

Section~3 formalises the coupled signal--observation dynamics as a rough
differential equation with a path-dependent Volterra drift. We introduce a
notion of solution for this class of equations based on the theory presented in the earlier section and establish global
existence and uniqueness of the resulting system. The latter is achieved by an extension of
the classical Yamada-Watanabe argument to the present rough path setting.

Section~4 develops the corresponding nonlinear filtering theory. We construct
a reference probability measure under which the observation process is reduced
to a driftless Volterra Gaussian process and derive a
Kallianpur--Striebel formula for the conditional distribution of the signal.
The resulting representation expresses the filter through an unnormalised
conditional expectation involving a likelihood process defined with respect
to the Brownian motion underlying the Volterra observation.

Section~5 turns this probabilistic representation into a pathwise one. We
construct a version of the filter which depends continuously on the
observation path and prove its stability in the rough path topology. To this
end, we develop the robust filtering method of Clark and Crisan
\cite{clark2005robust} in the Volterra Gaussian rough-path setting and combine
it with the integrability estimates of Cass, Litterer, and Lyons
\cite{cass2013integrability}. The resulting representation remains well
defined under smooth approximation, discretisation, and perturbation of the
observed path, and agrees almost surely with the original conditional
expectation when evaluated along the model observation.

In Section~6, under the assumption that \(X\), \(Y\), and \(B\) are
one-dimensional, we study the density of the unnormalised conditional measure.
We show that this density is characterised by a rough analogue of the Zakai
equation, which we derive by passing to the limit in equations driven by smooth
approximations of the observation. We establish the well-posedness of the
resulting equation by extending the rough viscosity solution framework
introduced by Friz and Oberhauser in \cite{friz2010rough} to the class of rough
partial differential equations arising in the present Volterra Gaussian
filtering problem. This yields a pathwise evolution equation for the
conditional density.

Finally, Section~7 studies the regularity of the conditional law of the signal.
Under a partial H\"ormander-type condition, we establish sufficient conditions
for the conditional distribution to admit a smooth density with respect to
Lebesgue measure. The main difficulty is to prove the required non-degeneracy
of the Malliavin covariance after conditioning on the observation. To overcome
this, we combine the conditional density criterion of
Nualart and Zakai \cite{nualart1989partial} with the inverse-moment method of
Cass, Hairer, Litterer, and Tindel \cite{cass2015smoothness}.

\section{Preliminaries}\label{section_preliminaries}
The filtering problem studied in this paper involves signal--observation
dynamics driven by Gaussian processes which need not be semimartingales and
may lie beyond the scope of classical pathwise integration. Before introducing
the equations themselves, we therefore recall the rough path framework used
to interpret such dynamics.

The issue is already visible in differential equations of the form
\[
    Z_t=z_0+\int_0^t g(Z_s)\,dI_s.
\]
When \(I\) has bounded variation, the integral is defined in the
Riemann--Stieltjes sense, while sufficiently regular paths can be treated by
Young integration. Rough path theory provides a unified framework which
recovers both regimes and extends them to rougher signals by enhancing \(I\)
to a rough path \(\mathbf I\), thereby yielding a continuous solution map
with respect to the enhanced driving signal.
In the present work, the driving signals are Volterra Gaussian processes.
These form a flexible class of Gaussian processes described through a Volterra
kernel \(K\), and include several examples of interest beyond the semimartingale
setting. Since the filtering equations we consider are driven by such
processes, we need to recall both the rough path machinery used to solve
equations driven by irregular signals and the Gaussian analytic framework used
to construct the corresponding rough path.

The goal of this preliminary section is therefore threefold. We first recall
the basic notions from rough path theory needed to formulate and solve rough
differential equations. 
We then introduce the class of Volterra Gaussian
processes used throughout the paper.
Finally, we review the Gaussian rough path construction on
which our approach to the filtering problem relies.

This material will be used to give a pathwise meaning to  the signal--observation system underlying the filtering
problem,
\[
\begin{pmatrix}
X_t \\ Y_t
\end{pmatrix}
=
\begin{pmatrix}
x_0 \\ Kf(X,Y)(t)
\end{pmatrix}
+
\int_0^t
\begin{pmatrix}
\sigma(X_s,Y_s) & 0 \\
0 & \mathbbm{1}
\end{pmatrix}
d
\mathbf I_s ,
\]
where the process $\mathbf I_s $ denotes the Gaussian rough path associated to the Volterra
Gaussian driving signal $I$. 
\subsection{Review of rough paths theory}\label{Intro_RP}

Let $\mathcal U$, $\mathcal V$ and $\mathcal W$ be real finite dimensional Banach spaces. 
For every $0 \leq a \leq b \leq T$ let $\Delta_{[a,b]}$ denote the simplex 
\[
\Delta_{[a,b]} := \{ (s,t) \in [a,b]^2 : s \le t \}.
\]

\begin{definition}[$p$-variation]
    Let \(p\ge 1\). For a two-parameter map
    \(h:\Delta_{[a,b]}\to \mathcal V \), we set
    \[
        \|h\|_{p-var;[a,b]}
        :=
            \sup_{\mathcal P }\left(\sum_{i}\|h_{s_i,s_{i+1}}\|^p_{\mathcal V}
            \right)^{1/p},
    \]
    where the supremum is taken over all finite partitions $\mathcal P = \{s_i\}$ of $[a,b]$.
 We then
    write
    \[
        \mathcal{C}^p(\Delta_{[a,b]},\mathcal V)
        :=
        \left\{
            h:\Delta_{[a,b]}\to \mathcal V:
            \|h\|_{p\text{-var};[a,b]}<\infty
        \right\}.
    \]
    
Finally, for a path \(g:[a,b]\to \mathcal V\), we say that
\(g\in \mathcal{C}^p([a,b],\mathcal V)\) whenever its increment map
\[
    g_{s,t}:=g_t-g_s
\]
belongs to \(\mathcal{C}^p(\Delta_{[a,b]},\mathcal V)\). In this case we set
\(
    \|g\|_{p\text{-var};[a,b]}
    :=
    \|g_{\cdot,\cdot}\|_{p\text{-var};[a,b]} .
\)
\end{definition}
Whenever there is no risk of ambiguity we will often write 
$\|h\|_{p}$ or $\|h\|_{p; [a,b]}$ in place of $ \|h\|_{p\text{-var};[a,b]}$.

For every $0\le a \le b \le T$, $\mathcal{C}^{0,p}([a,b],\mathcal V)$ is defined as the closure of $\mathcal{C}^1([a,b],\mathcal V)$ in $\mathcal{C}^p([a,b],\mathcal V)$.

For \(k\geq 0\), we write \(\mathcal V^{\otimes k}\) for the \(k\)-fold tensor product
of \(\mathcal V\), with the convention \(\mathcal V^{\otimes 0}=\mathbb{R}\). These tensor
products are always endowed with the projective tensor norm.

We denote by
\[
    T(\mathcal V):=\bigoplus_{k\geq 0} \mathcal V^{\otimes k}
\]
the tensor algebra over \(\mathcal V\), whose unit element is written \(\mathbf 1\).
We also write
\[
    T((\mathcal V)):=\prod_{k\geq 0} \mathcal V^{\otimes k}
\]
for the algebra of formal tensor series over \(\mathcal V\). We identify \(T(\mathcal V)\) with the subalgebra of \(T((\mathcal V))\) consisting of
series with only finitely many non-zero homogeneous components. For \(N\in\mathbb{N}\), we set
\[
    T^{\leq N}(\mathcal V):=\bigoplus_{k=0}^{N} \mathcal V^{\otimes k},
\]
and refer to it as the truncated tensor algebra up to level \(N\). Finally, \(G^N(\mathcal V)\) stands for the
step-\(N\) nilpotent group over \(\mathcal V\), equipped with its
Carnot--Carathéodory metric.
\begin{definition}[Geometric rough paths]
Let \(p\geq 1\). A geometric \(p\)-rough path over \(\mathcal V\) is a map
\(
    \mathbf I \in \mathcal C^{0,p}
    \bigl(\Delta_{[0,T]},G^{\lfloor p\rfloor}(\mathcal V)\bigr)
\)
satisfying
\[
    \mathbf I _{s,s}= \mathbf 1,
    \qquad
    \mathbf I_{s,u}
    =
    \mathbf I_{s,t}\otimes \mathbf I_{t,u},
    \qquad
    0\leq s\leq t\leq u\leq T .
\]
The space of such paths is denoted by
\(
    \Omega G_p([0,T],\mathcal V).
\)
\end{definition}

Geometric rough paths encode precisely the analytic bounds and algebraic
relations required to construct a pathwise integration theory for driving
signals whose sample paths are too irregular to have bounded variation.

Let $\mathrm{Hom}(\mathcal V, \mathcal U)$ denote the space of bounded linear maps endowed with the operator topology. 
For $\eta\in T((\mathcal V))$ and $k\in\mathbb{N}$, we write
\[
    \eta^{(k)}:=\operatorname{proj}_k(\eta), 
    \qquad
    \eta^{\leq k}:=\operatorname{proj}_{\leq k}(\eta),
    \qquad
    \eta^{<k}:=\operatorname{proj}_{<k}(\eta),
\]
where $\operatorname{proj}_k$ denotes the canonical projection onto
$\mathcal V^{\otimes k}$, while $\operatorname{proj}_{\leq k}$ and
$\operatorname{proj}_{<k}$ denote the canonical projections onto the levels
of degree at most $k$ and strictly less than $k$, respectively. If $g \in \mathrm{Hom}(T(\mathcal V), \mathcal U)$, for every non-negative integer $k$ we write $g^{(k)} = g \circ \operatorname{proj}_k$. 
Fixing a basis $(e_1,\ldots,e_d)$ of $\mathcal V$, for every word
$I=(i_1,\ldots,i_k)$ of length $k$ we write
\(e_I:=e_{i_1}\otimes\cdots\otimes e_{i_k}.
\)
The corresponding coordinate of $\eta^{(k)}$ is denoted by
$\eta^{(k), I}$. When $k=0$ we often write $\eta^{I}$ in place of $\eta^{(k), I}$.
Moreover, we define the left sided tensor multiplication operator \[
\mathbb L_{\eta}: T((\mathcal V)) \to T((\mathcal V)), \quad \mathbb L_{\eta} \xi \to \eta \otimes \xi.  
\] 

We are now ready to introduce controlled rough paths. This class provides the
natural functional framework for solutions of rough differential equations
driven by geometric rough paths. In particular, it will be the space in which
the signal and observation processes are shown to live.
\begin{definition}[$\mathbf{I}$-controlled paths]
Let \(p\in [1,\infty)\), set \(\kappa=\lfloor p\rfloor\), and let
\(\mathbf{I}\in \Omega G_{p}([0,T],
\mathcal V)\). Let
\(
    Z:[0,T]\to \mathrm{Hom}(T^{<\kappa}(\mathcal V),\mathcal U).
\)
We say that \(Z\) is an \(\mathbf{I}\)-controlled path if there exists a map
\(
    R:\Delta_{[0,T]}\to \mathrm{Hom}(T^{<\kappa}(\mathcal V), \mathcal U)
\)
such that, for all \(0\leq s\leq t\leq T\),
\begin{equation}\label{1D_controlled}
    Z_t
    =
    Z_s\circ \mathbb L_{\mathbf{I}_{s,t}^{<\kappa}}
    +
    R_{s,t}.
\end{equation}
Moreover, we require that  for
every \(j=0,\ldots,\kappa-1\),
\begin{equation}\label{regularity_constraint_remainder_1d}
    R^{(j)} \in
    \mathcal C^{p / (\kappa-j)}\left( \Delta_{[0, T]}, \mathrm{Hom}(\mathcal V ^{\otimes j}, \mathcal{U}) \right).
\end{equation}
\end{definition}
The space of controlled paths can be made into a Banach space (denoted by $\mathscr{D}_{\boldsymbol{I}, p}(\mathcal U)$) by introducing the norm 
\[
\|Z\|_{\mathscr{D}_{\mathbf{I}, p}; [s, t]} = \|Z_s\| + \sum_{j=0}^{\kappa-1}\|R^{(j)}\|_{\frac{p}{\kappa -j}; [s,t]}, \quad 0 \le s \le t \le T.
\]
Let \(\tilde{\mathbf I}\in \Omega G_{p}([0,T],\mathcal V)\) and let
\(\tilde Z\in \mathscr D_{\tilde{\mathbf I},p}(\mathcal U)\). Although \(Z\) and
\(\tilde Z\) are controlled by different rough paths, and therefore do not
belong to the same controlled path space, it is useful to introduce the
following comparison quantity. For \(0\leq s\leq t\leq T\), set
\[
    d_{[s,t]}(Z,\tilde Z)
    :=
    \|Z_s-\tilde Z_s\|
    +
    \sum_{j=0}^{\kappa-1}
    \bigl\|
        R^{Z,(j)}-R^{\tilde Z,(j)}
    \bigr\|_{\frac{p}{\kappa-j};[s,t]} .
\]
We refer to \(d_{[s,t]}(Z,\tilde Z)\) as a distance only by abuse of
terminology: it is not a metric on a single space, since \(Z\) and
\(\tilde Z\) are elements of the different spaces
\(\mathscr D_{\mathbf I,p}(\mathcal U)\) and
\(\mathscr D_{\tilde{\mathbf I},p}(\mathcal U)\), respectively.

We refer to $Z^{(0)}$ as the trace of $Z$ and to $R$ as the remainder of $Z$. Whenever ambiguity can arise from the context we will use the notation $R^Z$ in place of $R$.

The importance of the structure in the definition of a controlled rough path is that it gives the integrand a local
expansion in terms of the iterated increments of the driving rough path. This
allows one to define integration against \(\mathbf I\) through compensated
Riemann sums, leading to the following statement.
\begin{proposition}[Rough integral]\label{prop:rough_integral}
Let \(p\geq 1\), set \(\kappa=\lfloor p\rfloor\), and let
\(\mathbf I\in \Omega G_{p}([0,T],\mathcal V)\). If
\(
    Z \in \mathscr{D}_{\mathbf I,p}
    \bigl(\operatorname{Hom}(\mathcal V,\mathcal U)\bigr),
\)
then, for every \(0\leq s\leq t\leq T\), the limit
\[
    \lim_{|\mathcal P|\to 0}
    \sum_k\sum_{j=0}^{\kappa-1}
    Z_{u_k}^{(j)}
    \mathbf I_{u_k,u_{k+1}}^{(j+1)}
\]
exists and is independent of the sequence of partitions
\(\mathcal P=\{u_k\}\) of \([s,t]\), where
\(
    |\mathcal P|
    :=
    \max_{k}(u_{k+1}-u_k)
\)
denotes the mesh of \(\mathcal P\). We define the rough integral of \(Z\)
against \(\mathbf I\) by
\[
    \int_s^t Z_r\,d\mathbf I_r
    :=
    \lim_{|\mathcal P|\to 0}
    \sum_k\sum_{j=0}^{\kappa-1}
    Z_{u_k}^{(j)}
    \mathbf I_{u_k,u_{k+1}}^{(j+1)}.
\]
\end{proposition}

For \(n\in\mathbb N\) and \(\alpha\in(0,1]\), we denote by
\(C^{n+\alpha}(\mathcal V,\mathcal U)\) the space of functions
\(g\in C^n(\mathcal V,\mathcal U)\) whose \(n\)-th Fréchet derivative is
globally \(\alpha\)-Hölder continuous. We write
\(C_b^{n+\alpha}(\mathcal V,\mathcal U)\) for the subspace consisting of
those functions for which \(D^k g\) is bounded for every
\(0\leq k\leq n\), with the convention \(D^0g=g\), and
\(C_c^{n+\alpha}(\mathcal V,\mathcal U)\) for the subspace of compactly
supported functions. We also use the standard notation
\[
    C^\infty(\mathcal V,\mathcal U)
    :=
    \bigcap_{n\geq 1} C^n(\mathcal V,\mathcal U),
\]
and define \(C_b^\infty(\mathcal V,\mathcal U)\) and
\(C_c^\infty(\mathcal V,\mathcal U)\) by requiring, respectively, that
all derivatives, including the function itself, are bounded, or that the
function has compact support. 

Let \({\Delta}_{\shuffle}\) denote the shuffle coproduct on
\(T^{<\kappa}(\mathcal V)\), namely the unique algebra morphism
\[
    {\Delta}_{\shuffle}
    :
    T^{<\kappa}(\mathcal V)
    \to
    T^{<\kappa}(\mathcal V)\boxtimes T^{<\kappa}(\mathcal V)
\]
such that
\[
    {\Delta}_{\shuffle}(\mathbf 1)
    =
    \mathbf 1\boxtimes \mathbf 1,
    \qquad
    {\Delta}_{\shuffle}(\eta)
    =
    \eta\boxtimes \mathbf 1
    +
    \mathbf 1\boxtimes \eta,
    \qquad \eta\in \mathcal V.
\]
We write \({\Delta}^{m}_{\shuffle}\) for its \(m\)-fold iteration and $\tilde \Delta_\shuffle$ its reduced coproduct i.e. $\tilde \Delta_{\shuffle} \eta = \Delta_{\shuffle} \eta - \mathbf{1}\otimes \eta - \eta \otimes \mathbf 1$;
see Section~2 of \cite{cass2022combinatorial} for further details.

We are ready to introduce a fundamental property of controlled rough paths:  stability
under composition. This allows one to apply sufficiently regular nonlinear
maps to controlled paths, which is essential in the formulation of rough
differential equations.

Let \(g\in C^{\kappa}(\mathcal V,\mathcal W)\) and let
\(Z\in \mathscr D_{\mathbf I,p}(\mathcal U)\). Then \(g(Z)\) is again controlled by
\(\mathbf I\). Its controlled rough path structure is given by
\[
    g(Z)^{(0)}_t
    :=
    g(Z_t),
\]
and, for \(1\le k\le \kappa-1\),
\[
    g(Z)^{(k)}_t
    :=
    \sum_{m=1}^{k}
    \frac{1}{m!}
    D^m g(Z_t)
    \circ
    \bigl(Z_t\bigr)^{\otimes m}
    \circ
    \tilde {\Delta}^{m}_{\shuffle\mid \mathcal V^{\otimes k}} .
\]

For \(k\geq 1\), throughout this work, we use the canonical identification
\[
    \operatorname{Hom}
    \bigl(
        V^{\otimes(k-1)},
        \operatorname{Hom}(\mathcal V,\mathcal U)
    \bigr)
    \simeq
    \operatorname{Hom}(\mathcal V^{\otimes k},\mathcal U), \quad 
    g(v_1\otimes\cdots\otimes v_k)
    :=
    g(v_1\otimes\cdots\otimes v_{k-1})(v_k).
\]

\begin{definition}[Solution of a rough differential equation]
    Let $p\geq 1$, $\kappa = \floor p$, $\mathbf I \in \Omega G_{p}([0, T], \mathcal V)$ and $g \in C^\kappa(\mathcal U, \mathrm{Hom}(\mathcal V, \mathcal U))$. We say that $Z  \in \mathscr{D}_{\mathbf I, p}(\mathcal U) $ solves the rough differential equation
    \[
    Z_t = Z_0 + \int_0^t g(Z)_s d\mathbf{I}_s 
    \]
     if for every $t \in [0, T]$
     \begin{align*}
         Z^{(0)}_t &= Z^{(0)}_0 + \int_0^t g(Z)_s d\mathbf{I}_s,\\
         Z^{(k)}_t &= g(Z)^{(k-1)}_t,\qquad k \geq 1.
     \end{align*}
\end{definition}
A sufficient condition for the global well-posedness of a Rough Differential Equation (RDEs) is that $g \in C^{\kappa+1}_b(\mathcal U, \mathrm{Hom}(\mathcal V, \mathcal U))$; see Theorem 8.8 in \cite{GUBINELLI2010693}.

\subsection{Gaussian spaces and isonormal processes }\label{intro_malliavin}
In this section we recall the Gaussian spaces framework used
throughout the paper. The driving signal considered below is a Volterra
Gaussian process, and therefore its rough path lift is constructed from the
underlying Gaussian structure. We introduce the associated canonical Gaussian
space, its Cameron--Martin space, and the isonormal Gaussian process used to
define Malliavin derivatives later in this work. These objects will also be used later in the
proofs, in particular when deriving estimates for the lifted driving signal.

For a detailed introduction to Malliavin calculus we refer the reader to
\cite{nualart2009malliavin} and \cite{nourdin2012normal}.

We work on the canonical path space
\(
    \Omega=C_0([0,T],\mathbb R^d),
\)
equipped with the supremum norm and its Borel \(\sigma\)-algebra
\(\mathscr B(\Omega)\). For \(t\in[0,T]\), define the canonical coordinate
process \(I=\{I_t:t\in[0,T]\}\) by
\[
    I_t(\omega):=\omega(t),
    \qquad \omega\in\Omega.
\]
Let \(\mu\) be a centred Gaussian probability measure on
\((\Omega,\mathscr B(\Omega))\) under which the components of \(I\) are
independent and identically distributed. Since each coordinate map
\(\omega\mapsto I_t^k(\omega)\) is a continuous linear functional on
\(\Omega\), the process \(I\) is, under \(\mu\), a centred continuous
Gaussian process.

We denote its covariance function by
\[
    C_I(s,t)
    :=
    \mathbb E^\mu[I_s\otimes I_t]
    =
    \sum_{k,j=1}^d
    C_I^{k,j}(s,t)\,e_k\otimes e_j.
\]
Here \((e_1,\ldots,e_d)\) denotes the canonical basis of $\mathbb R^d$.

In particular, by independence and identical distribution of the components,
one has
\(
    C_I^{k,j}(s,t)=\delta_{k,j}C(s,t)
\)
for a scalar covariance function \(C\).

The Cameron--Martin space associated with \(\mu\) will be denoted by
\(\mathfrak H^d\). It is the reproducing kernel Hilbert space generated by
the covariance of \(I\). More precisely, for \(s\in[0,T]\) and
\(k\in[1:d]\), set
\[
    C_s^k(\cdot)
    :=
    \mathbb E^\mu[I_{\cdot} I_s^k]
    =
    \sum_{j=1}^d C_I^{j,k}(\cdot,s)e_j .
\]
Then \(\mathfrak H^d\) is the completion of the linear span of the functions
\begin{align}\label{span_e}
    \mathcal E
    :=
    \left\{ C_s^k:\, s\in[0,T],\ k\in[1:d] \right\}
\end{align}
with respect to the inner product determined by
\[
    \langle C_s^k, C_t^j\rangle_{\mathfrak H^d}
    =
    C_I^{k,j}(s,t).
\]
Since the components of \(I\) are independent, \(\mathfrak H^d\) is
isometrically isomorphic to the Hilbert direct sum
\(
    \mathfrak H^{(1)}\oplus\cdots\oplus\mathfrak H^{(d)},
\)
where \(\mathfrak H^{(k)}\) denotes the Cameron--Martin space associated
with the \(k\)-th component of \(I\). Whenever $d=1$ we will omit the superscript and write $\mathfrak{H}$ in place of $\mathfrak{H}^1$.

We shall also use the canonical Hilbert space \(\mathfrak H_1^d\) associated
with \(I\). This is defined as the completion of the linear span of
\[
    \mathcal E_1
    :=
    \left\{
        \mathbbm 1^k_{[0,s]}:\,
        s\in[0,T],\ k\in [1:d]
    \right\},
\]
where \(\mathbbm 1^k_{[0,s]}\) denotes the \(\mathbb R^d\)-valued indicator
function \(r\to \mathbbm 1_{[0,s]}(r)e_k\). Its inner product is defined
by
\[
    \left\langle
        \mathbbm 1^k_{[0,s]},
        \mathbbm 1^j_{[0,t]}
    \right\rangle_{\mathfrak H_1^d}
    =
    C_I^{k,j}(s,t).
\]
The map
\[
    \overline{\psi}:\mathfrak H_1^d\to \mathfrak H^d,
    \qquad
    \overline{\psi}\bigl(\mathbbm 1^k_{[0,s]}\bigr):=C_s^k,
\]
extends by linearity and completion to an isometric isomorphism.

Finally, the map
\[
    \vartheta:
    \mathbbm 1^k_{[0,t]}
    \to
    I_t^k
\]
extends to an isometry from \(\mathfrak H_1^d\) into
\(L^2(\Omega,\mathscr B(\Omega),\mu)\). Its image is the first Wiener chaos
associated with \(I\). For \(h\in\mathfrak H_1^d\), we write
\(
    I(h):=\vartheta(h).
\)
The family
\(
    \{I(h): h_1\in\mathfrak H_1^d\}
\)
is the isonormal Gaussian process over \(\mathfrak H_1^d\), that is,
\[
    \mathbb E^\mu[I(h)I(g)]
    =
    \sum_{k,j=1}^d
    \mathbb E^\mu\!\left[I^k(h^k)I^j(g^j)\right]
    =
    \sum_{k=1}^d
    \langle h^k,g^k\rangle_{\mathfrak H_1}
    =
    \langle h,g\rangle_{\mathfrak H_1^d}.
\]
The triple
\(
    (\Omega,\mathfrak H^d,\mu)
\)
is the abstract Wiener space associated with the canonical Gaussian process
\(I\).  
\subsection{Volterra Gaussian processes}\label{sec:VGP}

We now specialise the preceding Gaussian framework to the class of Volterra
Gaussian processes used as driving signals in this work. 

Let \(I=\{I_t:t\in[0,T]\}\) be a $\mathbb R^d$-valued centred continuous Gaussian process and
covariance
\(
    C_I(s,t).
\)
We say that \(I\) is a Volterra Gaussian process if there exists a kernel
\(K\in L^2([0,T]^2)\), with \(K(t,s)=0\) for \(s>t\), such that for every $k,j \in [1:d]$
\[
    C^{k, j}_I(s,t)
    =
    \delta_{k,j}\int_0^T K(t,r)K(s,r)\,dr
    =
    \delta_{k,j} \int_0^{s\wedge t} K(t,r)K(s,r)\,dr .
\]
Equivalently, on a suitable probability space, \(I\) admits the Wiener
integral representation
\[
    I_t
    =
    \int_0^t K(t,s)\,dV_s,
    \qquad t\in[0,T],
\]
for a $\mathbb R^d$-valued Brownian motion \(V\).

Associated with \(K\) is the Hilbert--Schmidt operator, still denoted by
\(K\),
\begin{equation}\label{kernel_operator}
    (Kh)(t)
    :=
    \int_0^T K(t,s)h_s\,ds
    =
    \int_0^t K(t,s)h_s\,ds,
    \qquad h\in L^2([0,T], \mathbb R^d).
\end{equation}
We also define the adjoint-type operator \(K^*\) initially on step functions
by
\begin{align}\label{K(t)_star}
    K^*\mathbbm 1_{[0,t]} = K(t,\cdot),
    \qquad t\in[0,T],
\end{align}
and extend it by linearity and completion. Thus, for the canonical Hilbert
space \(\mathfrak H_1^d\) associated with \(I\),
\begin{align}\label{isometry_h_h_1}
    \left\langle
        \mathbbm 1^k_{[0,t]},
        \mathbbm 1^j_{[0,s]}
    \right\rangle_{\mathfrak H_1^d}
    =
    C^{k, j}_I(t,s)
    =
    \left\langle
        K^*\mathbbm 1^k_{[0,t]},
        K^*\mathbbm 1^j_{[0,s]}
    \right\rangle_{L^2}.
\end{align}
Consequently, \(K^*\) extends to an isometry from \(\mathfrak H^d\) onto the
closed subspace $\mathfrak H_2^d$ of \(L^2([0,T], \mathbb R^d)\) that is generated by
\[
    \{K(t,\cdot)e_k:t\in[0,T], k \in [1:d]\}.
\]
When the kernel is sufficiently regular, this operator admits the usual
representation
\[
    (K^*h)(t)
    =
    h_tK(T,t)
    +
    \int_t^T \bigl(h_s-h_t \bigr)\,K(ds,t),
\]
for every \(h\) in the domain for which the right-hand side is well defined, see Section 2 in \cite{alos2001stochastic}.

If for every $k \in [1:d]$ and $t\in[0,T]$ 
\(
    \mathbbm 1_{[0,t]}\in \operatorname{Range}(K^*)  
\),
then one finds that it is possible to recover the standard Brownian motion $V$ via the isonormal process $I$, more precisely
\[
    V^k_t
    =
    I\bigl((K^*)^{-1}\mathbbm 1^k_{[0,t]}\bigr),
    \qquad t\in[0,T].
\]
In particular, \(V\) and \(I\) generate the same
filtration.

We now present the two main classes of processes that will serve as examples throughout this work: the Mandelbrot--Van Ness fractional Brownian motion, also known as Type I fractional Brownian motion, and the Riemann--Liouville fractional Brownian motion, also known as Type II fractional Brownian motion. 
For completeness, we recall the fractional operators used below. For
\(\alpha>0\), the left- and right-sided Riemann--Liouville fractional
integrals of \(f\in L^1([0,T])\) are defined, for almost every
\(t\in[0,T]\), by
\[
    (I_{0+}^{\alpha}f)(t)
    :=
    \frac{1}{\Gamma(\alpha)}
    \int_0^t (t-s)^{\alpha-1}f(s)\,ds,
    \qquad
    (I_{T-}^{\alpha}f)(t)
    :=
    \frac{1}{\Gamma(\alpha)}
    \int_t^T (s-t)^{\alpha-1}f(s)\,ds.
\]
Here $\Gamma$ denotes the Gamma function.
We adopt the convention that 
\(I_{0+}^0=I_{T-}^0=\operatorname{Id_{L^1}}\), where $\operatorname{Id_{L^1}}$ is the identity operator on $ L^1([0, T])$.
Setting \(n=\lceil\alpha\rceil\), the corresponding fractional
derivatives are defined by
\[
    D_{0+}^{\alpha}f
    :=
    \frac{d^n}{dt^n}I_{0+}^{\,n-\alpha}f,
    \qquad
    D_{T-}^{\alpha}f
    :=
    (-1)^n\frac{d^n}{dt^n}I_{T-}^{\,n-\alpha}f,
\]
for \(f\in I_{0+}^{\alpha}(L^1([0,T]))\) and
\(f\in I_{T-}^{\alpha}(L^1([0,T]))\), respectively.

\begin{example}[Fractional Brownian motion]\label{ex:fbm}
Let \(H\in(0,1)\). Fractional Brownian motion with Hurst parameter
\(H\) is the centred Volterra Gaussian process \(I\)
whose covariance function is
\[
    C_I(t,s)
    =
    \frac12
    \left(
        t^{2H}+s^{2H}-|t-s|^{2H}
    \right),
    \qquad s,t\in[0,T].
\]

The associated Volterra operator is an isomorphism
\[
    K:
    L^2([0,T])
    \longrightarrow
    I_{0+}^{H+\frac12}\bigl(L^2([0,T])\bigr).
\]
Suppressing the multiplicative normalising constants, it is given by
\[
    K h
    =
    \begin{cases}
        I_{0+}^{2H}
        \left(
            t^{\frac12-H}
            I_{0+}^{\frac12-H}
            \left(t^{H-\frac12}h\right)
        \right),
        & H<\frac12,\\[0.5em]
        I_{0+}^{1}h,
        & H=\frac12,\\[0.5em]
        I_{0+}^{1}
        \left(
            t^{H-\frac12}
            I_{0+}^{H-\frac12}
            \left(t^{\frac12-H}h\right)
        \right),
        & H>\frac12.
    \end{cases}
\]
Here and below, multiplication by powers of \(t\) is understood
pointwise.

The \(L^2([0,T])\)-adjoint of \(K\) is given by
\[
    K^* h
    =
    \begin{cases}
        t^{H-\frac12}
        I_{T-}^{\frac12-H}
        \left(
            t^{\frac12-H}
            I_{T-}^{2H}h
        \right),
        & H<\frac12,\\[0.5em]
        I_{T-}^{1}h,
        & H=\frac12,\\[0.5em]
        t^{\frac12-H}
        I_{T-}^{H-\frac12}
        \left(
            t^{H-\frac12}
            I_{T-}^{1}h
        \right),
        & H>\frac12.
    \end{cases}
\]
On \(\operatorname{Range}(K_H^*)\), its inverse is given, up to the
corresponding multiplicative constants, by
\[
    (K_H^*)^{-1}h
    =
    \begin{cases}
        D_{T-}^{2H}
        \left(
            t^{H-\frac12}
            D_{T-}^{\frac12-H}
            \left(t^{\frac12-H}h\right)
        \right),
        & H<\frac12,\\[0.5em]
        D_{T-}^{1}h,
        & H=\frac12,\\[0.5em]
        D_{T-}^{1}
        \left(
            t^{\frac12-H}
            D_{T-}^{H-\frac12}
            \left(t^{H-\frac12}h\right)
        \right),
        & H>\frac12.
    \end{cases}
\]
\end{example}

\begin{example}[Riemann--Liouville fractional Brownian motion]
\label{ex:rl-fbm}
Let \(H\in(0,1)\). The Riemann--Liouville fractional Brownian motion
with Hurst parameter \(H\) is the centred Volterra Gaussian process
\(I\) whose covariance function is
\[
    C_I(t,s)
    =
    \frac{1}{\Gamma\left(H+\frac12\right)^2}
    \int_0^{t\wedge s}
        (t-r)^{H-\frac12}
        (s-r)^{H-\frac12}
    \,dr,
    \qquad s,t\in[0,T].
\]

The associated Volterra operator is the isomorphism
\[
    K:
    L^2([0,T])
    \longrightarrow
    I_{0+}^{H+\frac12}\bigl(L^2([0,T])\bigr)
\]
given by
\[
    Kh
    =
    I_{0+}^{H+\frac12}h.
\]
Its \(L^2([0,T])\)-adjoint is
\[
    K^*h
    =
    I_{T-}^{H+\frac12}h,
\]
and hence
\[
    \operatorname{Range}(K^*)
    =
    I_{T-}^{H+\frac12}\bigl(L^2([0,T])\bigr).
\]
On \(\operatorname{Range}(K^*)\), the inverse of \(K^*\) is given by
\[
    (K^*)^{-1}h
    =
    D_{T-}^{H+\frac12}h.
\]
\end{example}

It is immediate to see that both classes include standard Brownian motion as the special case \(H= 1 /2\).

\subsection{Gaussian rough paths}
We now recall the Gaussian rough path construction used to enhance the driving
process. Throughout, let \(I=(I^1,\ldots,I^d)\) be a centred continuous Gaussian
process with independent and identically distributed components, and let \(C_I\)
denote its covariance function. The construction of the
enhancement is governed by the two-dimensional variation of \(C_I\), measured in
terms of its rectangular increments.

More precisely, for \(0\leq a<b\leq T\), \(0\leq c<d\leq T\), and a grid partition
\(
    \mathcal Q=\{(s_i,u_j):0\leq i\leq m,\ 0\leq j\leq n\}
\)
of \([a,b]\times[c,d]\), with
\(
    a=s_0<\cdots<s_m=b,\,
    c=u_0<\cdots<u_n=d,
\)
we write
\[
    C_I
    \begin{pmatrix}
        s_i & s_{i+1} \\
        u_j & u_{j+1}
    \end{pmatrix}
    :=
    C_I(s_i,u_j)+C_I(s_{i+1},u_{j+1})
    -C_I(s_i,u_{j+1})-C_I(s_{i+1},u_j).
\]
For \(p,q\geq 1\), we say that \(C_I\) has finite mixed
\((p,q)\)-variation on \([a,b]\times[c,d]\) if
\[
    \|C_I\|_{(p,q)\operatorname{-var};[a,b]\times[c,d]}
    :=
    \sup_{\mathcal Q}
    \left(
        \sum_i
        \left(
            \sum_j
            \left|
            C_I
            \begin{pmatrix}
                s_i & s_{i+1} \\
                u_j & u_{j+1}
            \end{pmatrix}
            \right|^p
        \right)^{q/p}
    \right)^{1/q}
    <\infty,
\]
where the supremum is taken over all grid partitions \(\mathcal Q\) of
\([a,b]\times[c,d]\).

When \(p=q=\rho\), we write
\[
    \|C_I\|_{\rho\operatorname{-var};[a,b]\times[c,d]}
    :=
    \|C_I\|_{(\rho,\rho)\operatorname{-var};[a,b]\times[c,d]},
\]
and refer to this quantity as the two-dimensional \(\rho\)-variation of
\(C_I\).

The construction of the rough path separates naturally into two regimes. If the driving signal is
regular enough to be treated in the Young regime, then the path itself already
defines the required level-one rough path and no additional Gaussian
enhancement is needed. The covariance assumptions below are only used in the
genuinely rough regime, where one constructs the higher-order iterated
integrals from the Gaussian structure of \(I\).

\begin{condition}\label{condition_covariance_rho_var}
Let \(I=(I^1,\ldots,I^d)\) be a centred continuous Gaussian
process with independent and identically distributed components, and let \(C_I\)
be its covariance function.
We assume that
one of the following two alternatives holds.

\begin{enumerate}
    \item[(A)] There exists \(p\in[1,2)\) such that, almost surely,
    \[
        (s,t)\mapsto I_{s,t}
        \in
        \mathcal C^{0,p}\bigl(\Delta_{[0,T]},\mathbb R^d\bigr)
        =
        \mathcal C^{0,p}\bigl(\Delta_{[0,T]},G^1(\mathbb R^d)\bigr).
    \]
    In this case, \(I\) itself defines a level-one geometric \(p\)-rough path.
    \item[(B)] There exists \(\rho\in[1,2)\) such that the covariance function
    \(C_I\) of one component of \(I\) has finite two-dimensional
    \(\rho\)-variation on \([0,T]^2\), that is,
    \[
        \|C_I\|_{\rho\operatorname{-var};[0,T]^2}<\infty .
    \]
    In this case, the canonical Gaussian rough path lift of \(I\) is
    constructed from its covariance structure.
\end{enumerate}
\end{condition}
In particular, condition \emph{(A)} holds if, for some \(p\in(1,2)\), the
sample paths of \(I\) belong almost surely to
\(\mathcal C^{p'}\bigl([0,T],\mathbb R^d\bigr)\) for every \(1 \leq p'< p\).

Under Condition \emph{(B)}, for every
\(p\in(2\rho,4)\), the process \(I\)
admits a canonical geometric \(p\)-rough path enhancement such that, for
almost every \(\omega\in\Omega\),
\[
    \mathbf I(\omega)
    \in
    \Omega G_p ([0, T], \mathbb{R}^d).
\]
The enhancement can be realised as a measurable map, defined
\(\mu \)-almost surely on the path space of \(B\),
\[
    \mathscr{S} :
    C([0,T],\mathbb R^d)
    \to
    \Omega G_p ([0, T], \mathbb{R}^d).
\]
We shall refer to \(  \mathscr{S} (I)\) as the measurable enhancement of \(I\).
See Theorem~35 and Section~6.3 in \cite{friz2010differential}.
In particular, for almost
every \(\omega\in\Omega\), equations driven by \(I\) can be interpreted
pathwise as rough differential equations driven by the deterministic rough
path \(\mathbf I(\omega)\).\newline
It is possible to show that for every $H \in (0, 1)$ both the fractional Brownian motion and the Riemann-Liouville frational Brownian motion satisfy Condition \ref{condition_covariance_rho_var} (see Appendix \ref{app:examples}).

\begin{remark}
Let
\(
k_t(r):=K(t,r),
\)
a simple sufficient condition for \(C_I\) to have finite two-dimensional
\(\rho\)-variation, \(\rho\geq1\), can be deduced from the inequality

\[
|C_I|_{\rho\text{-var};[0,T]^2}
\leq
\int_0^T
   \bigl|t\mapsto k_t(r)\bigr|_{\rho\text{-var};[0,T]}^2\,dr,
\]
which is  a consequence of Minkowski inequality.
For more singular kernels, let \(\mathcal Q:=\{(t,t)\}\) and define, on
rectangles contained in \(\mathcal Q^c\),
\[
\mu_K((s,t]\times(u,v])
:=\langle k_t-k_s,k_v-k_u\rangle_{L^2}.
\]
Suppose that \(\mu_K\) extends to a signed Radon measure on \(\mathcal Q^c\), whose
positive part \(\mu_K^+\) is finite and has continuous distribution function
\[
(a,b)\to
\mu_K^+\bigl((0,a]\times(0,b]\cap \mathcal Q^c\bigr).
\]
If, for some \(C,h>0\),
\(
\|k_t-k_s\|_{L^2}^2\leq C|t-s|^{1/\rho}
\) and \(
\langle k_t-k_s,k_v-k_u\rangle_{L^2}\geq0
\)
whenever \([u,v]\subseteq[s,t]\) and \(t-s\leq h\), then \(C\) has finite two-dimensional \(\rho\)-variation;
see Theorem~2.2 in \cite{friz2016jain}.
\end{remark}

\section{Setup and well-posedness of the signal-observation system} \label{section_setup}

Let \((\Omega,\mathcal F,\mathbb P)\) be a probability space, and let \(
    B=(B^1,\ldots,B^{d_B})
    \) and \(
    U=(U^1,\ldots,U^{d_Y})
\)
be Volterra Gaussian processes defined on this space. We assume that the
collection
\(
    \{B^1,\ldots,B^{d_B}, \newline U^1,\ldots,U^{d_Y}\}
\)
consists of mutually independent and identically distributed Gaussian
processes. We also assume that these
processes are associated with the deterministic kernel \(K\) and  satisfy Conditions
\ref{condition_covariance_rho_var}.

Let \(x_0\) be an \(\mathbb R^{d_X}\)-valued random variable such that
\[
    x_0\in \bigcap_{r\geq 1} L^r(\Omega,\mathbb R^{d_X}),
\]
with law \(\mu_0\), and assume that \(x_0\) is independent of the driving
Gaussian processes. We endow the probability space with the usual augmentation
\((\mathcal F_t)_{t\in[0,T]}\) of the filtration generated by \(x_0\), \(B\),
and \(U\). More precisely, if
\[
    \mathcal F_t^0
    :=
    \sigma(x_0)\vee \sigma(B_s:s\le t)\vee \mathcal \sigma(U_s:s\le t),
\]
then
\[
    \mathcal F_t
    :=
    \bigcap_{r>t}\left(\mathcal F_r^0\vee\mathscr Z\right),
    \qquad t\in[0,T],
\]
where \(\mathscr Z\) denotes the collection of \(\mathbb P\)-null sets in $\mathcal F$. In
particular, \((\mathcal F_t)_{t\in[0,T]}\) is complete and right-continuous.

On the filtered probability space
\(
    (\Omega,\mathcal F,(\mathcal F_t)_{t\in[0,T]},\mathbb P)
\), for \(\mathbb P\)-almost every \(\omega\in\Omega\), we consider the rough differential equation
\begin{equation}\label{sensor_equation}
\begin{pmatrix}
X_t \\ Y_t
\end{pmatrix}
=
\begin{pmatrix}
x_0 \\ Kf(X^{(0)},Y^{(0)})(t)
\end{pmatrix}
+
\int_0^t
\begin{pmatrix}
\sigma(X,Y)_u & 0 \\
0 & \mathbbm{1}
\end{pmatrix}
d\boldchar{(}\mathbf B,\mathbf U\boldchar{)}_u .
\end{equation}
To simplify the notation, we suppress the dependence of all random variables on
\(\omega\) throughout the paper, unless it is required for clarity. Here and
throughout, \(\boldchar{(}\mathbf B,\mathbf U\boldchar{)}\) denotes a joint
enhancement of \((B,U)\) to a \(p\)-rough path, for some \(p\in[1,4)\) chosen
in accordance with Condition~\ref{condition_covariance_rho_var}. We refer to
\(\boldchar{(}\mathbf B,\mathbf U\boldchar{)}\) as the joint enhancement, or
equivalently the joint lift, of \(B\) and \(U\).
For notational convenience, we also write
\[
    Z_t := (X_t,Y_t),
    \qquad
    \mathbf I := \boldchar{(}\mathbf B,\mathbf U\boldchar{)} \in G \Omega_p([0, T], \mathbb R^{d_I}),
\]
and \(Kf(X, Y)\) instead of \(Kf\bigl(X^{(0)}, Y^{(0)}\bigr)\) whenever no confusion can arise.
Then \eqref{sensor_equation} can be written compactly as
\[
    Z_t
    =
    z_0 + Kb(Z)(t)
    +
    \int_0^t g(Z)_u\,d\mathbf I_u,
\]
where
\(
    z_0=(x_0,0) \in \mathbb R^{d_Z},
\)
and the coefficients \(b\) and \(g\) are induced by
\eqref{sensor_equation}.

The aim of this section is to formulate an appropriate notion of solution for \eqref{sensor_equation} within the framework of controlled rough paths, and to establish global well-posedness.

\begin{definition}[Solution for Volterra RDE]\label{operator_drift_solution} Let $I = (I^1, \cdots, I^{d_I})$ a Volterra Gaussian process with kernel $K$ that satisfies Condition \ref{condition_covariance_rho_var}. Let $b \in C^1(\mathbb  R^{d_Z}, \mathbb  R^{d_Z})$ $g \in C^\kappa(\mathbb  R^{d_Z}, \mathrm{Hom}(\mathbb  R^{d_I}, \mathbb  R^{d_Z}))$ and $\mathbf I$ be the geometric $p$-rough path constructed from $I$. We say that $Z  \in \mathscr{D}_{\mathbf I, p}(\mathbb  R^{d_Z}) $ solves the rough differential equation
    \[
        Z_t
    =
    z_0
    +
    Kb\bigl(Z\bigr)(t)
    +
    \int_0^t g(Z)_s\,d\mathbf I_s, 
    \]
     if, for every $t \in [0, T]$,
     \begin{align*}
         Z^{(0)}_t
    &=
    z_0
    +
    Kb\bigl(Z\bigr)(t)
    +
    \int_0^t g(Z)_s\,d\mathbf I_s,\\
         Z^{(k)}_t &= g(Z)^{(k-1)}_t.
     \end{align*}
\end{definition}

To verify that \(Z\) is a controlled path in the sense of the preceding definition, we first require sufficient regularity of the deterministic Volterra term \(Kb\). This is guaranteed by the following condition.
\begin{condition}\label{cameron_martin_condition}
Let \(\mathfrak H^d\) denote the Cameron--Martin space associated with the centered Gaussian Volterra process \(I=(I^1,\ldots,I^d)\) with i.i.d components. There exist
\(q \geq 1\) and \(C>0\) such that, for every \(h \in \mathfrak{H}^d\),
\[
    \|h\|_{q\text{-var};[0,T]}
    \leq C \|h\|_{\mathfrak{H}^d},
\]
where
\(
    \frac1p+\frac1q>1
    \).
\end{condition}
When \(p\in[1,2)\), Condition \ref{cameron_martin_condition} is automatic,
provided \(I\) is realised as a Gaussian random variable in
\(\mathcal C^{p}\). Indeed, the Cameron--Martin space
then embeds continuously into the same \(p\)-variation path space. Hence the
condition holds with \(q=p\), and \(p<2\) gives
\(
    1 /p+1/q=2 / p>1.
\) 

Both fractional Brownian motion and Riemann--Liouville fractional Brownian
motion satisfy Condition \ref{cameron_martin_condition}. The proof of this fact
is deferred to the Appendix \ref{app:examples}.
The role of the condition is twofold. First, through Proposition \ref{regularity_Kf}, this condition gives the pathwise regularity needed to treat \(Kb(Z)\) as
a drift of finite \(q\)-variation, with \(q\) complementary to the \(p\)-variation
regularity of the rough driver. Second, the same assumption is the standard
Cameron--Martin embedding used in the Gaussian rough path theory, and will later
provide the integrability estimates needed for the results of Section \ref{section_robustness} and \ref{section_zakai}.

Despite this assumption, there is a minor issue at the level of the controlled-path expansion. The current regularity of \(Kb(Z)\) does not guarantee that this term
can to be absorbed into the zero-level remainder. Under Condition
\ref{cameron_martin_condition}, however, \(Kb(Z)\) has finite \(q\)-variation
whenever \(b(Z)\) is square-integrable in time; in particular, this is the case
when \(b\) is bounded. Thus, by choosing the controlled-path topology with a
slightly larger roughness exponent \(p_*>p\), we can ensure that \(Kb(Z)\) has
the required regularity to be included in the zero-level remainder. The higher
levels of the controlled path can then be treated as in the classical RDE
setting as required by Definition \ref{operator_drift_solution}.

More precisely, fix $\kappa = \floor{p}$ an auxiliary exponent \(p_*>p\) such that
\[
    q \leq \frac{p_*}{\kappa}
    \qquad \text{and} \qquad
    \frac{p_*}{\kappa} < \frac{p}{p-1}.
\]
Such a choice is possible. Indeed, by the complementary Young condition
\(
    1 / p+ 1 / q>1,
\)
we have
\(
    q <  p / (p-1).
\)
Moreover, since \(\kappa=\lfloor p\rfloor\), one has \(\kappa>p-1\), and hence
\(
    {p} / {\kappa} < {p} / ({p-1}).
\)
Therefore the interval
\[
    \left( \max\{p,\kappa q\}, \frac{\kappa p}{p-1} \right)
\]
is non-empty, and we may choose \(p_*\) in this interval. Then \(p_*>p\),
\(q < p_*/\kappa\), and
\[
    \frac{p_*}{\kappa}<\frac{p}{p-1}.
\]

To establish global well-posedness for \eqref{sensor_equation}, we require the
filtration generated by the Volterra Gaussian process \(I\) to coincide with
that generated by its associated Brownian motion \(V\). As discussed in
Section~\ref{sec:VGP}, this is equivalent to the following condition.

\begin{condition}\label{condition_filtration}
For every \(t\in[0,T]\),
\(
    \mathbbm{1}_{[0,t]}\in\operatorname{Range}(K^*).
\)
\end{condition}

In Appendix~\ref{app:examples}, we verify that both Type~I and Type~II
Volterra Gaussian processes satisfy Condition~\ref{condition_filtration}.

Finally, throughout this part of the work, we impose the following regularity
assumptions on the coefficients \(f\) and \(\sigma\).

\begin{condition}\label{vf_condition_section_well_posedness}
The coefficients satisfy
\(
    f\in C_b^2
    \bigl(
        \mathbb R^{d_X}\times\mathbb R^{d_Y},
        \mathbb R^{d_Y}
    \bigr)
\)
and
\(
    \sigma\in C_b^{\kappa+1}
    \bigl(
        \mathbb R^{d_X}\times\mathbb R^{d_Y},
        \operatorname{Hom}(\mathbb R^{d_B},\mathbb R^{d_X})
    \bigr).
\)
\end{condition}

We are now in a position to state the well-posedness result for the
signal--observation dynamics considered in this work.

\begin{theorem}\label{weaK_solution_existence}
Suppose that Conditions~\ref{condition_covariance_rho_var},
\ref{cameron_martin_condition}, \ref{condition_filtration}, and
\ref{vf_condition_section_well_posedness} hold. Then there exist an
\((\mathcal F_t)\)-adapted process
\[
    Z\in
    \mathscr D_{\mathbf B,p_*}
    \bigl([0,T],\mathbb R^{d_Z}\bigr)
\]
and a set \(\Omega_2\subseteq\Omega\) of full measure such that, for every
\(\omega\in\Omega_2\), the path \(Z(\omega)\) solves
\eqref{sensor_equation} in the sense of
Definition~\ref{operator_drift_solution}.
\end{theorem}

The well-posedness theory underlying
Theorem~\ref{weaK_solution_existence} is developed in
Appendix~\ref{app:proofs}. We defer this analysis to the appendix in order to
keep the main text focused on the filtering problem. The argument is divided
into two parts. First, we establish local well-posedness and global uniqueness
for \eqref{sensor_equation} by means of a ball-preserving argument. We then
extend the classical Yamada--Watanabe argument to the rough-path setting,
thereby obtaining the full-measure set \(\Omega_2\) and the
\((\mathcal F_t)\)-adapted process \(Z\) appearing in
Theorem~\ref{weaK_solution_existence}.

\begin{remark}[Additional finite-variation drifts]
For clarity of exposition, we formulate the model without separate
finite-variation drift terms. This is not an essential restriction. One may,
for example, consider
\[
\begin{pmatrix}
\widetilde X_t \\ \widetilde Y_t
\end{pmatrix}
=
\begin{pmatrix}
x_0 \\ Kf(\widetilde X,\widetilde Y)(t)
\end{pmatrix}
+
\int_0^t
\begin{pmatrix}
\varphi\bigl(\widetilde X^{(0)}_u,\widetilde Y^{(0)}_u\bigr)\\
\varsigma\bigl(\widetilde Y^{(0)}_u\bigr)
\end{pmatrix}
\,du
+
\int_0^t
\begin{pmatrix}
\sigma(\widetilde X,\widetilde Y)_u & 0 \\
0 & \mathbbm{1}
\end{pmatrix}
d\boldchar{(}\mathbf B,\mathbf U\boldchar{)}_u ,
\]
provided that, for each of the results below, \(\varphi\) and \(\varsigma\)
satisfy the corresponding regularity, boundedness, and growth assumptions
imposed there on \(f\).

The inclusion of these terms creates no additional rough-integration
difficulty. Indeed, the system may be regarded as an RDE driven by the
augmented rough path
\[
    \boldchar{(}t,\mathbf B,\mathbf U\boldchar{)},
\]
where the time component has bounded variation and its mixed iterated
integrals with the rough components are canonically defined. In particular,
the signal drift \(\varphi\) does not affect the Girsanov transformation used
for the observation process; it only produces the usual additional
first-order term in the associated generators and filtering equations.

An additional observation drift can be treated by introducing
\[
    \overline Y_t
    :=
    \widetilde Y_t
    -
    \int_0^t
    \varsigma\bigl(\widetilde Y^{(0)}_u\bigr)\,du.
\]
Under the assumptions above, this is an invertible causal transformation:
given \(\overline Y\), the path \(\widetilde Y\) is uniquely recovered from
\[
    \widetilde Y_t
    =
    \overline Y_t
    +
    \int_0^t
    \varsigma\bigl(\widetilde Y^{(0)}_u\bigr)\,du.
\]
Consequently, \(\widetilde Y\) and \(\overline Y\) generate the same completed
filtration. The Girsanov transformation can therefore be applied to
\(\overline Y\), while the additional drift contributes only
finite-variation terms. We omit these terms from the main formulation in
order to keep the notation and the subsequent filtering equations concise.
\end{remark}

\section{The Kallianpur-Striebel formula} \label{section_KS_formula}
Having established well-posedness of the signal-observation pair, we now
turn to the associated filtering problem. The aim is to describe, at each
time \(t\), the conditional distribution of the signal \(X_t\) given the
observations accumulated up to time \(t\). This information is encoded by
the observation filtration
\[
    \mathcal{Y}_t
    :=
    \sigma\{Y_s \,:\, 0\leq s\leq t\}\vee\mathscr Z,
\]
where \(\mathscr Z\) denotes the collection of $\mathbb P$-null sets of  \(\mathcal F\).

More precisely, for every bounded measurable test function
\(\phi:\mathbb R^{d_X}\to\mathbb R\), the filter is defined by
\[
    \pi_t(\phi)
    :=
    \mathbb E^{\mathbb P}
    \bigl[\phi(X^{(0)}_t)\mid \mathcal Y_t\bigr].
\]
Equivalently, \(\pi_t\) is the conditional law of \(X^{(0)}_t\) given the
observation history, evaluated against the test function \(\phi\). To keep the notation light, whenever it is clear from the context we write $\phi(X_t)$ in place of $\phi(X^{(0)}_t)$.  

The first step in this analysis is to establish the Kallianpur--Striebel
formula. This formula rewrites the filter under a reference probability
measure \(\mathbb Q\), under which the observation process has a simpler
law. More precisely, let \(Y(\cdot)\) denote the isonormal Gaussian process over
\(\mathfrak H^{d_Y}\) canonically associated with the
\(\mathbb R^{d_Y}\)-valued Volterra Gaussian process \(Y\). Assume $Y$ satisfies Condition \ref{condition_filtration}. Setting
\[
    \mathcal Y
    :=
    \sigma\{Y_s:0\leq s\leq T\} \vee\mathscr Z,
    \qquad
    W_t^i
    :=
    Y\left((K^*)^{-1}\mathbbm 1_{[0,t]}^i\right),
    \qquad
    i\in[1:d_Y],
\]
where
\(
    \mathbbm 1_{[0,t]}^i\) is defined as in Section \ref{section_preliminaries},
one obtains
\[
    \pi_t(\phi)
    =
    \frac{\widetilde{\pi}_t(\phi)}
         {\widetilde{\pi}_t(1)}
    :=
    \frac{
        \mathbb E^{\mathbb Q}
        \bigl[\phi(X_t)\Lambda_t\mid\mathcal Y\bigr]
    }{
        \mathbb E^{\mathbb Q}
        \bigl[\Lambda_t\mid\mathcal Y\bigr]
    },
\]
where
\[
    \Lambda_t
    :=
    \frac{d\mathbb P}{d\mathbb Q}\bigg|_{\mathcal F_t}
    =
    \exp\Bigg(
        \sum_{j=1}^{d_Y}
        \int_0^t f_j(X_s,Y_s)\,dW_s^j
        -
        \frac12
        \sum_{j=1}^{d_Y}
        \int_0^t
        \bigl(f_j(X_s,Y_s)\bigr)^2\,ds
    \Bigg).
\]
Here $W$ is the standard Brownian motion associated to $Y$ (see Section \ref{section_preliminaries}) and the stochastic integral is understood in  the It\^o sense.
Despite the non-Markovian Volterra structure of the model, the resulting Kallianpur--Striebel formula is formally identical to its classical diffusion counterpart: under the reference measure, the conditional distribution is obtained by weighting with the usual exponential likelihood ratio. See, for instance, Proposition 3.16 in \cite{bain2009fundamentals}.

By the results of the previous section, under the measure \(\mathbb Q\), the
process \(Y\) is a Volterra Gaussian process with kernel \(K\), while the
law of \((B^j)_{j\in[1:d_B]}\) remains unchanged. Moreover, by an argument
analogous to that used in Theorem~\ref{weaK_solution_existence}, the
pair \(Z=(X,Y)\) satisfies
\begin{equation}\label{no_sensor_equation}
\begin{pmatrix}
X_t \\ Y_t
\end{pmatrix}
=
\begin{pmatrix}
x_0 \\ 0
\end{pmatrix}
+
\int_0^t
\begin{pmatrix}
\sigma(X,Y)_s & 0 \\
0 & \mathbbm{1}
\end{pmatrix}
d\boldchar{(}\mathbf B, \mathbf Y\boldchar{)}_s ,
\end{equation}
Equivalently, writing \(Z_t=(X_t,Y_t)\) and \(z_0=(x_0,0)\), we obtain the
compact form
\[
    Z_t
    =
    z_0
    +
    \int_0^t g(Z)_s\,d\mathbf J_s,
\]
where the function \(g\) is implicitly defined by the block matrix appearing in
\eqref{no_sensor_equation} and $\mathbf J = \boldchar{(}\mathbf B, \mathbf Y\boldchar{)}$.

Since \(\mathbb{Q}\) and \(\mathbb{P}\) are equivalent measures, we can  write by conditional Bayes formula 
\[
\pi_t(\phi) = \frac{\mathbb{E}^{\mathbb{Q}}\big[ \phi(X_t) \Lambda_t \, \big| \, \mathcal{Y}_t\big]}{\mathbb{E}^{\mathbb{Q}}\big[\Lambda_t \, \big| \, \mathcal{Y}_t\big]}.
\]  
The remaining task is to demonstrate that, in this expression, the filtration \(\mathcal{Y}_t\) can be replaced by \(\mathcal{Y}\). As the next Proposition makes precise, this is just a consequence of Condition \ref{condition_filtration}.

\begin{proposition}[Kallianpur--Striebel formula]
\label{Kallianpur_Striebel_formula}
Suppose that Conditions~\ref{condition_covariance_rho_var},
\ref{cameron_martin_condition}, \ref{condition_filtration}, and
\ref{vf_condition_section_well_posedness} hold. Then, for every bounded
measurable function \(\phi:\mathbb R^{d_X}\to\mathbb R\) and every
\(t\in[0,T]\),
\[
    \pi_t(\phi)
    =
    \frac{
        \mathbb E^{\mathbb Q}
        \bigl[\phi(X_t)\Lambda_t\mid\mathcal Y\bigr]
    }{
        \mathbb E^{\mathbb Q}
        \bigl[\Lambda_t\mid\mathcal Y\bigr]
    }
    \qquad \mathbb Q\text{-a.s.}
\]
\end{proposition}

\begin{proof}
Under Condition~\ref{condition_filtration}, the filtration generated by the
\(\mathbb Q\)-Volterra Gaussian process \(Y\) coincides with the filtration
\((\mathcal F_t^W)_{t\in[0,T]}\) generated by the associated Brownian motion
\(W^i_t := Y((K^*)^{-1} \mathbbm 1^i_{[0, t]})\) for $i : [1:d_Y]$. Hence conditioning on the observation history is equivalent to
conditioning on the Brownian filtration.

By Proposition~3.15 in \cite{bain2009fundamentals}, for every integrable
\(\mathcal F^W_t\)-measurable random variable \(H\),
\[
    \mathbb E^{\mathbb Q}
    \bigl[H\mid\mathcal F_t^W\bigr]
    =
    \mathbb E^{\mathbb Q}
    \bigl[H\mid\mathcal F^W\bigr],
\]
where \(
    \mathcal F^W:=\sigma\{W_s:0\leq s\leq T\}.
\)
Applying this identity with
\(
    H=\phi(X_t)\Lambda_t
    \) and \(
    H=\Lambda_t
\),
and using again that \(W\) and \(Y\) generate the same filtration, yields
\[
    \mathbb E^{\mathbb Q}
    \bigl[\phi(X_t)\Lambda_t\mid\mathcal Y_t\bigr]
    =
    \mathbb E^{\mathbb Q}
    \bigl[\phi(X_t)\Lambda_t\mid\mathcal Y\bigr],
\]
and similarly for \(\Lambda_t\). Therefore the usual Kallianpur--Striebel
formula with conditioning on \(\mathcal Y_t\) may equivalently be written
with conditioning on \(\mathcal Y\), giving the desired identity.
\end{proof}

The preceding discussion also highlights that in the particular case where the signal coefficient does not depend on the
observation process, the framework reduces to a more familiar filtering
setting. Under \(\mathbb Q\), the law of the signal driver \(B\) is
unchanged, while the observation process is equivalent to the Brownian
observation, i.e.
\(
    \tilde  Y_t := W_t .
\)
Consequently, under \(\mathbb Q\), the filtering problem is equivalent to
the signal-observation system
\[
    X_t
    =
    x_0+\int_0^t \sigma(X)_s\,d \boldchar{(}\mathbf B, \mathbf W\boldchar{)}_s,
    \qquad
    \tilde Y_t = W_t .
\]
Thus, in this uncoupled case, the signal equation is driven only by the
rough path \(\mathbf B\), while the observation equation is driven by the
Brownian motion \(W\). In particular, one does not need to construct a
rough path enhancement involving both the signal driver and the observation
path.

This observation places the present framework in relation to existing
filtering models. 
It may be viewed as the Stratonovich counterpart of the model introduced in \cite{coutin1999abstract}, extended here to possibly rougher driving signals and it contains, as a special case, filtering problems in
which the signal is driven by fractional Brownian motion and the observation
is driven by Brownian motion. The general formulation considered here goes
beyond this setting in two directions: the signal noise may be any Volterra
Gaussian process admitting a rough path enhancement, and the signal dynamics
may depend on the observation path itself. The joint rough path construction
is precisely what allows this latter, fully coupled, situation to be treated
pathwise.

\section{Robust representation of the filter}\label{section_robustness}

The next step in the analysis of the filtering problem is to construct  a robust version of the filter.
This is important for both theoretical and practical reasons. For each fixed
\(t\), the conditional expectation
\[
    \pi_t(\phi)
    =
    \mathbb E^\mathbb{P}[\phi(X_t)\mid \mathcal Y_t]
\]
is defined only up to null sets. Equivalently, one may write
\[
    \pi_t(\phi)=\mathfrak F_t^\phi(Y)
\]
for some Borel measurable functional \(\mathfrak F_t^\phi\) of the observation path, but
such a functional is not unique: it may be modified on any set which is null
under the law of \(Y\). This lack of uniqueness is harmless if one works only
under the exact reference model, but it becomes problematic when the filter is
evaluated on observed data.

Indeed, in applications the observation path is never available as an ideal
continuous-time sample from the model. It is observed through discrete data,
numerical approximations, measurement errors, or under a model which is only an
approximation of the true data-generating mechanism. In such situations, a
purely measurable version of the conditional expectation gives no stability
with respect to perturbations of the observation path. A robust representation
remedies this by selecting a version of the filter which is obtained by
evaluating a continuous functional on path space. The resulting object can then
be applied not only to paths sampled exactly from the model, but also to nearby
or approximating observation paths. In this sense, robustness provides a
pathwise form of stability under discretisation, approximation, and moderate
model misspecification.

In the present rough-path setting, the relevant input is not merely the
observation path \(Y\), but the enhanced path associated with the pair
\(\hat Y := (Y,W)\), where \(W\) is the \(\mathbb Q\)-Brownian motion appearing in the
innovation representation of the observation. Although \(Y\) and \(W\) are in
general distinct processes, this does not introduce additional unobservable
information. By Condition \ref{condition_filtration}, the filtrations
generated by \(Y\) and \(W\) coincide, up to completion. Hence \(W\) is
observable from \(Y\), and the joint enhancement
\(\hat{\mathbf{Y}} = (\mathbf Y,\mathbf W)\), constructed as in Appendix \ref{app:enlarged_VGP}, represents the same information as the observation
filtration. We therefore seek, for each test
function \(\phi\), a continuous functional in the \(p\)-variation rough-path
topology,
\(
    {\mathfrak p}^{\phi}_t:
     \Omega G_{p}([0, T], \mathbb{R}^{2 d_Y})
    \to \mathbb R,
\)
such that, for every \(t\in[0,T]\),
\[
    \mathfrak p_t^\phi(\hat{\mathbf{Y}})
    =
    \pi_t(\phi),
    \qquad
    \mathbb P\text{-a.s}.
\]
The continuity of
\(\mathfrak p_t^\phi\) is the robust property: whenever a sequence of enhanced
observation paths converges to \(\hat{\mathbf{Y}}\) in the rough-path
topology, the corresponding values of the filter converge to
\(\pi_t(\phi)\).

The argument follows the classical robustness approach to semimartingale filtering developed in \cite{clark1978design} and \cite{clark2005robust}; see also Chapter~5 of \cite{bain2009fundamentals}. It is combined with the integrability estimates for Gaussian rough differential equations established in \cite{cass2013integrability}. The main idea is to construct a version of the filter that depends on the observation only through its rough-path lift. To this end, we introduce a decoupling construction by enlarging the probability space with an independent copy of the Volterra Gaussian noise driving the signal, while retaining the original observation process and initial condition. On the resulting product space, the signal dynamics can be solved conditionally on a fixed realization of the observation, with the remaining randomness carried by the independent signal noise. This separates the probabilistic roles of the signal and observation and allows the conditional expectations defining the filter to be represented as pathwise functionals of the observed rough path.

A central role in this section is played by the \enquote{rough version} of the process $\Lambda$ introduced in the previous section. We denote this process by $\Theta$. As shown below, $\Theta$ is constructed using the rough path $\hat{\mathbf{J}}$, as defined in Proposition \ref{existence_geometric_lift}, which is the geometric Gaussian rough path constructed above the processes 
\[
\hat J  :=\big(B^1, \dots, B^{d_B}, Y^1, \dots, Y^{d_Y}, W^1, \dots, W^{d_Y}\big).
\]

In this framework, the process $\Theta$ is constructed by replacing the stochastic integral in the exponential martingale $\Lambda$ with the corresponding integral with respect to $\hat{\mathbf{J}}$.

Explicitly, we define  $\Theta = \exp(\Xi)$ where 
\begin{align*}
\Xi(\omega)_t & := \int_0^t F(X, Y)_u(\omega) \, d\hat{\mathbf{J }}_u(\omega)  
 - \frac{1}{2} \sum_{j=1}^{d_Y} \int_0^t (f_j(X_s, Y_s)(\omega))^2 \, ds 
\end{align*}
and  we define \[F: \mathbb R^{d_X} \times \mathbb R^{d_Y} \to \mathrm{Hom}(\mathbb R^{d_X + 2d_Y}, \mathbb R), \qquad
F(x,y)[\eta] := \sum_{j=1}^{d_Y} f_j(x,y)\eta^{j + d_X + d_Y}.
\] 
\begin{proposition}\label{ito_conversion}
Suppose that Conditions~\ref{condition_covariance_rho_var},
\ref{cameron_martin_condition}, \ref{condition_filtration}, and
\ref{vf_condition_section_well_posedness} hold.
Then for any \(0 \leq s \leq t \leq T\), the identity
\[
\sum_{j=1}^{d_Y}\int_s^t f_j(X_u,Y_u)\,dW_u^j
=
\int_s^t F(X,Y)_u\,d\hat{\mathbf J}_u
\]
holds \(\mathbb Q\)-almost surely. In particular $\Theta$ is a version of $\Lambda$.
\end{proposition}
\begin{proof}
By the definitions of the It\^o integral and of the rough integral
(Proposition~\ref{prop:rough_integral}), it is enough to show that, for any
sequence \(\mathcal P^n\) of finite partitions of \([s,t]\) with mesh size tending to
zero,
\begin{equation}\label{brownian_integration_rough}
\sum_{t_k \in \mathcal P^n}
F^{(1)}_{t_k}\hat{\mathbf J}^{(2)}_{t_k,t_{k+1}}
\to 0,
\qquad
\sum_{t_k \in \mathcal P^n}
F^{(2)}_{t_k}\hat{\mathbf J}^{(3)}_{t_k,t_{k+1}}
\to 0,
\end{equation}
in \(L^2(\mathbb Q)\).

We first consider the second-level correction. By construction, the only
non-zero coordinates of \(F^{(1)}\) which contribute to the above sum are those
for which the last coordinate of the iterated integral corresponds to a
Brownian component of \(\hat J\), while the preceding coordinate corresponds to
a Volterra Gaussian component. Since
\[
F^{(1)}(X,Y)_u
=
DF(X_u,Y_u)
\begin{pmatrix}
\sigma(X_u,Y_u)\\
0
\end{pmatrix}
\]
is bounded, it is enough to estimate, for such coordinates \(i,j\),
\[
\sum_k F^{(1),ij}_{t_k}
\hat{\mathbf J}^{(2),ij}_{t_k,t_{k+1}} .
\]
Using the definition of $\hat J$, we have the It\^o integral
\[
\hat{\mathbf J}^{(2),ij}_{t_k,t_{k+1}}
=
\int_{t_k}^{t_{k+1}}
\hat J^i_{t_k,u}\,d\hat J_u^j ,
\]
and therefore, by the It\^o isometry,
\begin{align*}
\mathbb E\left[
\left|
\sum_k F^{(1),ij}_{t_k}
\hat{\mathbf J}^{(2),ij}_{t_k,t_{k+1}}
\right|^2
\right]
&\leq
C_F
\sum_k
\mathbb E\left[
\left|
\int_{t_k}^{t_{k+1}}
\hat J^i_{t_k,u}\,d\hat J_u^j
\right|^2
\right] \\
&=
C_F
\sum_k
\int_{t_k}^{t_{k+1}}
\mathbb E\left[
\left|\hat J^i_{t_k,u}\right|^2
\right]\,du \\
&\leq
C_F |t-s|
\max_k
\sup_{t_k \leq r \leq v \leq t_{k+1}}
\mathbb E\left[
\left|\hat J^i_{r,v}\right|^2
\right].
\end{align*}
Since the Volterra Gaussian components are mean-square continuous, the final
quantity tends to \(0\) as \(|\mathcal P^n| \to 0\). Hence the first convergence
in \eqref{brownian_integration_rough} follows.

The third-level correction is treated in the same way. The only contributing
coordinates are those whose last coordinate is Brownian, so that the relevant
third-level increments can be written as It\^o integrals of second-level
increments. Using the boundedness of \(F^{(2)}\) and the It\^o isometry gives
\[
\mathbb E\left[
\left|
\sum_k
F^{(2)}_{t_k}
\hat{\mathbf J}^{(3)}_{t_k,t_{k+1}}
\right|^2
\right]
\leq
C_F |t-s|
\max_k
\sup_{t_k \leq r \leq v \leq t_{k+1}}
\mathbb E\left[
\left|\hat{\mathbf J}^{(2)}_{r,v}\right|^2
\right],
\]
which tends to \(0\) by the \(L^2\)-continuity of the second level of
\(\hat{\mathbf J}\). This proves \eqref{brownian_integration_rough}, and
therefore the desired identity.
\end{proof}

As stated above, for the remainder of this section, we employ a decoupling argument by introducing the probability space \((\check{\Omega}, \check{\mathcal{F}}, \check{\mathbb{Q}})\), defined as the product  
\[
(\check{\Omega}, \check{\mathcal{F}}, \check{\mathbb{Q}}) := (\Omega \times \overline{\Omega}, \mathcal{F} \otimes \overline{\mathcal{F}}, \mathbb{Q} \otimes \overline{\mathbb{Q}}),
\]
where \((\overline{\Omega}, \overline{\mathcal{F}}, \overline{\mathbb{Q}})\) is an independent copy of \((\Omega, \mathcal{F}, \mathbb{Q})\). The space $(\overline{\Omega}, \overline{\mathcal{F}}, \overline{\mathbb{Q}})$ supports the Volterra Gaussian process \((B^j)_{j \in [1:  d_B]}\). We define analogously the sets $\overline \Omega_1$ and $\overline \Omega_2$
as the copy of the sets $\Omega_1$ and $ \Omega_2$ introduced previously in this work.

The processes \(Y\), \(B\)  and \(x_0\) are naturally extended to this new probability space by defining  
\[
\overline{B}(\omega, \overline{\omega}) := (B^1(\overline{\omega}), \dots, B^{d_B}(\overline{\omega})),\quad  \check{Y}(\omega, \overline{\omega}) := Y(\omega), \quad \check{W}(\omega, \overline{\omega}) := W(\omega), \quad \check{x}_0(\omega, \overline{\omega}) := x_0(\omega).
\]

Applying the decoupling argument to  \(\hat{\mathbf{J}}\), we obtain its decoupled counterpart, which we denote by \(\check{\mathbf{J}}\).

Moreover, by Condition~\ref{condition_filtration}, the Brownian motion \(W\),
regarded as a
\(C_0([0,T],\mathbb R^{d_Y})\)-valued random variable, is measurable with
respect to \(Y\). Since the corresponding path spaces are Polish (see Proposition 5.38 in \cite{friz2010multidimensional}), there
exists a Borel measurable map
\[
\mathscr W:
C_0\bigl([0,T],\mathbb R^{d_Y}\bigr)
\to
C_0\bigl([0,T],\mathbb R^{d_Y}\bigr)
\]
such that
\(
    W=\mathscr W(Y)\),
    $\mathbb Q$-almost surely.
We fix one such version of \(\mathscr W\).

Recall that the Gaussian rough-path lift is realised by a Borel measurable
map
\[
\hat{\mathscr S}:
C_0\bigl([0,T],\mathbb R^{d_B+2d_Y}\bigr)
\to
\Omega G_p\bigl([0,T],\mathbb R^{d_B+2d_Y}\bigr),
\]
such that
\(
    \hat{\mathbf J}
    =
    \hat{\mathscr S}(B,Y,W),
    \) \(\mathbb Q\)-almost surely.

We may therefore define
\[
\mathscr{\hat S}:
C_0\bigl([0,T],\mathbb R^{d_Y}\bigr)
\to
\bigcap_{r\geq1}
L^r\left(
    \overline\Omega,
    \Omega G_p\bigl([0,T],\mathbb R^{d_B+2d_Y}\bigr)
\right)
\]
by
\[
\mathscr{\hat S}(y)
:=
\left[
\overline\omega
\to
\hat{\mathscr S}\bigl(
    B(\overline\omega),
    y,
    \mathscr W(y)
\bigr)
\right],
\]
where the three components are arranged in the same order as in
\(\hat J=(B,Y,W)\). Thus, the third argument in the definition above is the
Brownian path associated with \(y\).

In particular,
\[
\hat{\mathbf Y}(\omega)
\xrightarrow{\ \operatorname{proj}_Y\ }
Y(\omega)
\xrightarrow{\ \mathscr{\hat S}\ }
\check{\mathbf J}(\omega,\cdot),
\]
where
\[
\begin{aligned}
\check{\mathbf J}(\omega,\cdot)
&:=
\mathscr{\hat S}\bigl(Y(\omega)\bigr) \\
&=
\left[
\overline\omega
\to
\hat{\mathscr S}\bigl(
    B(\overline\omega),
    Y(\omega),
    \mathscr W(Y(\omega))
\bigr)
\right] \\
&=
\left[
\overline\omega
\to
\hat{\mathscr S}\bigl(
    B(\overline\omega),
    Y(\omega),
    W(\omega)
\bigr)
\right],
\end{aligned}
\]
where the last equality holds \(\mathbb Q\)-almost surely.

Here, \(\operatorname{proj}_Y\) denotes the canonical projection onto the
coordinates \(1,\ldots,d_Y\) of the first level of the tensor algebra. This
projection is continuous with respect to the \(p\)-variation topology. Since
\(\mathscr W\) and \(\hat{\mathscr S}\) are Borel measurable, the map
\[
(y,\overline\omega)
\to
\hat{\mathscr S}\bigl(
    B(\overline\omega),
    y,
    \mathscr W(y)
\bigr)
\]
is jointly measurable. Together with the corresponding \(L^r\)-integrability
estimate, this implies that
\(
\omega
\to
\check{\mathbf J}(\omega,\cdot)
\)
defines an
\(
L^r\left(
    \overline\Omega,
    \Omega G_p\bigl([0,T],\mathbb R^{d_B+2d_Y}\bigr)
\right)
\)
-valued random variable for every \(r\geq1\).
We extend the notational convention adopted for $\mathbf J$ to all the other realisations of joint rough paths appearing in this section.

We define the decoupled process \(\check{\Xi}\) to be the measurable map \(\check{\Xi}: \overline{\Omega} \to \mathbb{R}\)  
\begin{align*}
& \check{\Xi}(\omega)_t(\cdot) := \int_0^t F(\check{X}, \check{Y})_s(\omega, \cdot) \, d\check{\mathbf{J }}_s(\omega, \cdot) - \frac{1}{2}  \sum_{j=1}^{d_Y} \int_0^t ( f(\check{X}_s, \check{Y}_s)(\omega, \cdot))^2 \, ds ,
\end{align*}  
where $\check{X}$ is the solution to the RDE \eqref{no_sensor_equation} driven by the decoupled version of ${\mathbf{J}}$ and with $\check{Y}$ in place of $Y$.

The 
function \({\mathfrak p}^{\phi}\), providing the robust version of the filter, is defined as follows:  
\begin{equation}\label{robust_filter_definition}
 \mathfrak p^{\phi}_t(\hat{\mathbf{Y}} (\omega)) := \frac{\tilde {\mathfrak p}^{\phi}_t(\hat{\mathbf{Y}} (\omega))}{\tilde {\mathfrak p}^{1}_t(\hat{\mathbf{Y}} (\omega))} := \frac{\mathbb{E}^{\overline{\mathbb Q}}\left[ \phi(\check
 {X}_t)(\omega, \cdot) \exp \check{\Xi}(\omega)_t \right]}{\mathbb{E}^{\overline{\mathbb Q}}\left[  \exp \check{\Xi}(\omega)_t \right]}.
\end{equation}
Before proceeding with the next proof, we recall the notion of a greedy
partition, introduced in
\cite{cass2013integrability}. This construction decomposes the time interval
into successive subintervals on which the \(p\)-variation of the driving
rough path is controlled by a prescribed threshold. 
\begin{definition}
    Let $\boldsymbol \zeta \in  \Omega G_{p}([0, T], \mathcal V)$, partition $[0, T]$ as follows:
    \begin{align*}
    & s^\zeta(\alpha)_0 := 0 \quad \quad s^\zeta(\alpha)_i := \inf\Big\{ z>s^\zeta(\alpha)_{i-1} : \|\boldsymbol{\zeta}\|^p_{p; [ s^\zeta(\alpha)_{i-1}, z]} \geq \alpha \Big\} \wedge T,
\end{align*}
where we assume $\inf \varnothing  = \infty$. We define $N^{\zeta}_{\alpha}$ as 
\[
N^{\zeta}_{\alpha} := \inf\{ i \in \mathbb{N}_0 \mid s^\zeta(\alpha)_i < T \}.
\]
\end{definition}
An immediate consequence of the definition of $N^{\zeta}_\alpha$ is that we have the bound
\begin{equation}\label{upper_bound_N}
\alpha N^{\zeta}_\alpha \leq \|\boldsymbol{\zeta}\|^p_{p; [0, T]}.
\end{equation}
\begin{lemma} \label{integrability_N_alpha}
Suppose that Conditions~\ref{condition_covariance_rho_var},
\ref{cameron_martin_condition}, \ref{condition_filtration}, and
\ref{vf_condition_section_well_posedness} hold.
Then there exists a set \(\Omega_3\subset\Omega\) with
\(\mathbb P(\Omega_3)=1\) such that, for every
\(\omega\in\Omega_3\) and $\alpha > 0$
\[
    \exp\bigl(N^{(\check W,\check B)}_\alpha(\omega,\cdot)\bigr),
    \qquad
    \exp\bigl(N^{(\check Y,\check B)}_\alpha(\omega,\cdot)\bigr)
\]
belong to \(L^r(\overline\Omega,\overline{\mathbb P})\) for every
\(r\geq 1\).
\end{lemma}

\begin{proof}
By Condition~\ref{cameron_martin_condition}, the Cameron--Martin spaces
associated with \((\check W,\check B)\) and \((\check Y,\check B)\) are
continuously embedded into the space of paths of finite \(q\)-variation, for
every \(q\geq 1\) satisfying
\(
    1 / p+1/q>1 .
\)
Moreover, the corresponding enhanced Gaussian rough paths are obtained as
limits of their piecewise linear approximations (see Appendix \ref{app:enlarged_VGP}). Hence the assumptions of
Theorem~6.3 in \cite{cass2013integrability} are satisfied. Consequently,
\[
    \exp\bigl(N^{(\check W,\check B)}_\alpha\bigr),
    \exp\bigl(N^{(\check Y,\check B)}_\alpha\bigr)
    \in
    \bigcap_{r\geq 1} L^r(\check\Omega,\check{\mathbb P}).
\]

In particular, for every \(n\in\mathbb N\),
\[
    \int_{\Omega}
    \int_{\overline\Omega}
    \exp\bigl(nN^{(\check W,\check B)}_\alpha(\omega,\bar\omega)\bigr)
    \,d\overline{\mathbb P}(\bar\omega)\,
    d\mathbb P(\omega)
    <\infty,
\]
and the same estimate holds with \((\check W,\check B)\) replaced by
\((\check Y,\check B)\). By Fubini's theorem, for each \(n\in\mathbb N\)
there exists a set \(\Omega_{(n)}\subset\Omega\) with
\(\mathbb P(\Omega_n)=1\) such that, for every \(\omega\in\Omega_{(n)}\),
both
$\exp\bigl(N^{(\check W,\check B)}_\alpha(\omega,\cdot)\bigr),
    $ and $
    \exp\bigl(N^{(\check Y,\check B)}_\alpha(\omega,\cdot)\bigr)
$
belong to \(L^n(\overline\Omega,\overline{\mathbb P})\).

Set
\[
    \Omega_3:=\bigcap_{n=1}^\infty \Omega_{(n)} .
\]
Then \(\mathbb P(\Omega_3)=1\). For every \(\omega\in\Omega_3\), the two
random variables above belong to \(L^n(\overline\Omega)\) for every
\(n\in\mathbb N\).
\end{proof}
\begin{lemma}\label{rough_exp_martingale_bounded}
Suppose that Conditions~\ref{condition_covariance_rho_var},
\ref{cameron_martin_condition}, \ref{condition_filtration}, and
\ref{vf_condition_section_well_posedness} hold.
Then there exists a set \(\Omega_4\) of full measure such that, for every
\(\omega\in\Omega_4\) and every \(r\geq 1\), there exists a set
\(O_\omega\), bounded in the \(p\)-variation topology, such that
\[
    \sup_{t\in[0,T]}
    \sup_{\tilde \omega \in  \hat{\mathbf Y}^{-1}(O_\omega)}
    \mathbb E^{\overline{\mathbb Q}}
    \left[
        \exp\left(r\check{\Xi}(\tilde\omega)_t\right)
    \right]
    <\infty.
\]
\end{lemma}

\begin{proof}
It is enough to consider the case $r=1$, since the remaining cases follow
by the same argument.

Define
\[
\check{\Xi}^1(\omega)_t(\cdot)
:=
\int_0^t
F(\check{X}_s,\check{Y}_s)(\omega,\cdot)
\,d\check{\mathbf{J }}_s(\omega,\cdot),
\]
and set $\check{\Theta}^1=\exp(\check{\Xi}^1)$.

The argument relies on the fact that, since $\check X$ and $\check Y$ are
solutions of an RDE and $\exp$ is smooth, the process $\check \Theta^1$ can
be viewed as a controlled rough path. More precisely, the pair
$(\check Z,\check \Theta^1)$ solves the RDE
\[
\begin{pmatrix}
    \check{Z}_t \\
    \check{\Theta}^1_t
\end{pmatrix}
= \begin{pmatrix}
    \check{z}_0 \\
    1
\end{pmatrix} + 
\int_0^t \begin{pmatrix}
    g(\check{Z})_s \\
    \check{\Theta}^1 \hat{b}(\check{Z})_s
\end{pmatrix}
d\check{\mathbf{J }}_s, \quad  .
\]
See Lemma 2.26 and Theorem 2.27 in \cite{cass2022combinatorial}.

Moreover, by the results of \cite{cass2013integrability}, there exists a
positive constant $C$, depending on $f,\sigma,T,p$, such that the trace
of $\check \Theta$ can be estimated in terms of
$\|\check{\mathbf J }\|_{p}$ and $N^{\check B}$. Suppressing the superscript
identifying level $0$, we obtain for any $\overline \omega \in \overline \Omega_2$
\begin{equation}\label{linear_RDE}
|\check{\Theta}^1(\cdot)_t(\overline{\omega})|
\leq
C
\|\check{\mathbf{J }}(\cdot,\overline{\omega})\|^p_{p;[0,t]}
\exp\left(
    C\left(N^{\check{B}(\cdot,\overline{\omega})}_C+1\right)
\right).
\end{equation}

Let $q$ be as in Condition \ref{cameron_martin_condition}.
For any \(M>0\), define
\[
    O_\omega
    :=
    \left\{
        \hat{\mathbf Y}(\tilde\omega)
        :
        \tilde\omega\in \Omega,\ 
        \left\|
            \hat{\mathbf Y}(\tilde\omega)
            -
            \hat{\mathbf Y}(\omega)
        \right\|_{p \wedge q;[0,T]}
        \leq M^{1/p}
    \right\}.
\]
This set is non-empty. Indeed, by
Condition~\ref{cameron_martin_condition}, it contains the shifts of the
Gaussian rough path by elements of the Cameron-Martin space \(\mathfrak H^{d_Y}\) of sufficiently
small norm. This can be shown by via calculations analogous to the ones in Theorem \ref{weaK_solution_existence}.

By the definition of $O_\omega$, we have via Young-Loeve estimate that $\check{\mathbb Q}$-a.s.
\[
\left\|
\check{\mathbf{J }}(\tilde{\omega},\overline{\omega})
\right\|_{p;[s,t]}
\leq
(1+1\vee M)
\left\|
\check{\mathbf{J }}(\omega,\overline{\omega})
\right\|_{p;[s,t]},
\qquad [s,t]\subseteq[0,T].
\]
Consequently, for every $\alpha>0$ and every
$i\leq N^{\check{B}(\omega,\overline{\omega})}_{\alpha}$,
\[
s^{\check{B}(\omega,\overline{\omega})}
\left(
    \frac{\alpha}{(1+1\vee M)^p}
\right)_i
\leq
s^{\check{B}(\tilde{\omega},\overline{\omega})}(\alpha)_i .
\]
Therefore, by the definition of $N^{\check B}_{\cdot}$, setting
\(
M_1:=\alpha / {(1+1\vee M)^p},
\)
we obtain
\(
N^{\check{B}(\tilde \omega,\overline{\omega})}_{\alpha}
\leq
N^{\check{B}(\omega,\overline{\omega})}_{M_1}.
\)
Combining this estimate with \eqref{linear_RDE}, and setting
$M_2=C/M_1$, gives
\begin{equation}\label{bound_theta_2}
\begin{aligned}
|\check{\Theta}^1(\tilde{\omega})_t(\overline{\omega})|
&\leq
C
(1+1\vee M)^{2p}
\|\check{\mathbf{J }}(\omega,\overline{\omega})\|^p_{p;[0,t]}  \\
&\qquad \times
\exp\left(
    C\left(
        N^{\check{B}(\omega,\overline{\omega})}_{M_2}+1
    \right)
\right).
\end{aligned}
\end{equation}

Fix now $\alpha>0$ and consider the partitions
\[
\{s^{(\check{W},\overline{B})}(\alpha)_i\},
\qquad
\{s^{(\check{Y},\overline{B})}(\alpha)_i\},
\qquad
\{s^{(\check{Y},\check{W})}(\alpha)_i\}.
\]
Their union, denoted by $\mathcal P$, is such that, for any
$t_i\in\mathcal P$,
\[
\|\check{\mathbf{J }}(\omega,\overline{\omega})\|^p_{p;[t_i,t_{i+1}]}
\leq 3\alpha .
\]
It follows that
\begin{align*}
N^{\check{B}(\omega,\overline{\omega})}_{\alpha}
&\leq
N^{(\check{W},\overline{B})(\omega,\overline{\omega})}_{\alpha/3}
+
N^{(\check{Y},\overline{B})(\omega,\overline{\omega})}_{\alpha/3}
+
N^{(\check{Y},\check{W})(\omega,\overline{\omega})}_{\alpha/3}
\\
&\leq
N^{(\check{W},\overline{B})(\omega,\overline{\omega})}_{\alpha/3}
+
N^{(\check{Y},\overline{B})(\omega,\overline{\omega})}_{\alpha/3}
+
\frac{\alpha}{3}
\|\boldchar{(}\check{\mathbf{Y}},\check{\mathbf{W}}\boldchar{)}
(\omega)\|^p_{p;[0,T]},
\end{align*}
where in the last step we used \eqref{upper_bound_N}.

Together with Lemma \ref{integrability_N_alpha}, this estimate implies that,
for $\omega$ a.s.,
\[
\exp\left(
    N^{\check{B}(\omega,\cdot)}_{M_2}
\right)
\in
\bigcap_{r\geq 1} L^r(\overline{\Omega}).
\]
Using Fernique's estimate from \ref{existence_geometric_lift} and the estimate
\eqref{bound_theta_2}, we conclude that
\begin{align*}
\mathbb{E}^{\overline{\mathbb Q}}
\left[
    |\check{\Theta}^1(\tilde{\omega})_t|
\right]
&\leq
\mathbb{E}^{\overline{\mathbb Q}}
\left[
    C(1+1\vee M)^{2p}
    \|\check{\mathbf{J }}(\omega,\cdot)\|^p_{p;[0,t]}
    \exp\left(
        C\left(
            N^{\check{B}(\omega,\cdot)}_{M_2}+1
        \right)
    \right)
\right]
\\
&\leq
\mathbb{E}^{\overline{\mathbb Q}}
\left[
    C^2(1+1\vee M)^{2p}
    \|\check{\mathbf{J }}(\omega,\cdot)\|^{2p}_{p;[0,t]}
\right]^{1/2}
\\
&\qquad \times
\mathbb{E}^{\overline{\mathbb Q}}
\left[
    \exp\left(
        2C\left(
            N^{\check{B}(\omega,\cdot)}_{M_2}+2
        \right)
    \right)
\right]^{1/2}.
\end{align*}
The right-hand side is finite. This proves the desired bound for
$\check{\Theta}^1$.

The remaining terms in
$\exp(\check{\Xi}-\check{\Xi}^1)$ are handled by the same argument, and we
therefore omit the details.
\end{proof}
From the previous Lemma and the estimate $|\exp(a) - \exp(b)| \leq (\exp(a) + \exp(b))|a-b| $ we can deduce the following result:
\begin{corollary}
In the setting of the preceding lemma, for every \(\omega\in\Omega_4\) and every \(\tilde\omega\) such that
\(\hat{\mathbf Y}(\tilde\omega)\in O_\omega\), we have
\[
    \sup_{t\in[0,T]}
    \mathbb{E}^{\overline{\mathbb Q}}
    \left[
        \left|
            \check{\Xi}(\omega)_t
            -
            \check{\Xi}(\tilde\omega)_t
        \right|^2
    \right]^{1/2}
    \leq
    C_{V,p,T,\hat Y,Z}
    \left\|
        \hat{\mathbf Y}(\omega)
        -
        \hat{\mathbf Y}(\tilde\omega)
    \right\|_{p;[0,T]}.
\]
The constant \(C_{V,p,T,\hat Y,Z}\) may be chosen independently of
\(\tilde\omega\) subject to
\(\hat{\mathbf Y}(\tilde\omega)\in O_\omega\).
\end{corollary}

\begin{proposition}
Suppose that Conditions~\ref{condition_covariance_rho_var},
\ref{cameron_martin_condition}, \ref{condition_filtration}, and
\ref{vf_condition_section_well_posedness} hold.
Then for any $t \in [0, T]$ and  realization $\omega \in \Omega_4$ the function $\tilde{\mathfrak{p}}_t^\phi$ is locally Lipschitz  from the space $(\Omega G_p([0, T], \mathbb{R}^{2d_B}), \|\cdot \|_{p})$ into $(\mathbb{R}, |\cdot|)$.   
\end{proposition}
\begin{proof}
We follow the method of proof detailed in Lemma 5.6 of \cite{bain2009fundamentals} and in Theorem 6 in \cite{crisan2013robust}.  
Fix $O_\omega$ as in Lemma \ref{rough_exp_martingale_bounded} and let $\hat{\mathbf{Y}}(\tilde{\omega})$ be an element of this set.

We have
\begin{align*}
& \left| \tilde{\mathfrak p}_t^{\phi}(\hat{\mathbf{Y}}(\omega)) - \tilde{\mathfrak p}_t^{\phi}(\hat{\mathbf{Y}}(\tilde{\omega})) \right| \\
& \leq \mathbb{E}^{\overline{\mathbb Q}}\left[ \left| \phi(\check{X}(\omega, \cdot)_t) \exp(\check{\Xi}(\omega)_t) -  \phi(\check{X}(\tilde{\omega}, \cdot)_t) \exp(\check{\Xi}(\tilde{\omega})_t) \right| \right] \\
& \leq \mathbb{E}^{\overline{\mathbb Q}}\left[ \left| \phi(\check{X}(\omega, \cdot)_t) - \phi(\check{X}(\tilde{\omega}, \cdot)_t)\right| \exp(\check{\Xi}(\omega)_t) \right] \\
&\qquad +   \mathbb{E}^{\overline{\mathbb Q}}\left[ \left| \phi(\check{X}(\tilde{\omega}, \cdot)_t) \right| \left| \exp(\check{\Xi}(\omega)_t) - \exp(\check{\Xi}(\tilde{\omega})_t)  \right|\right] \\ 
& \leq \mathbb{E}^{\overline{\mathbb Q}}\left[ \left| \phi(\check{X}(\omega, \cdot)_t) - \phi(\check{X}(\tilde{\omega}, \cdot)_t)\right|^2\right]^{\frac{1}{2}}\mathbb{E}^{\overline{\mathbb Q}}\left[\exp(2\check{\Xi}(\omega)_t) \right]^{\frac{1}{2}} \\
& \qquad +  \left\| \phi\right\|_{\infty} \mathbb{E}^{\overline{\mathbb Q}}\left[  \left| \exp(\check{\Xi}(\omega)_t) - \exp(\check{\Xi}(\tilde{\omega})_t)  \right|\right] \\ 
& \leq \mathbb{E}^{\overline{\mathbb Q}}\left[ \left| \phi(\check{X}(\omega, \cdot)_t) - \phi(\check{X}(\tilde{\omega}, \cdot)_t)\right|^2\right]^{\frac{1}{2}}\mathbb{E}^{\overline{\mathbb Q}}\left[\exp(2\check{\Xi}(\omega)_t) \right]^{\frac{1}{2}} \\
& \qquad +  \left\| \phi\right\|_{\infty} \mathbb{E}^{\overline{\mathbb Q}}\left[  \left| \check{\Xi}(\omega)_t -\check{\Xi}(\tilde{\omega})_t  \right|^2\right]^{\frac{1}{2}}\\
&\qquad \times \mathbb{E}^{\overline{\mathbb Q}}\left[  \left| \exp(\check{\Xi}(\omega)_t) + \exp(\check{\Xi}(\tilde{\omega})_t)  \right|^2\right]^{\frac{1}{2}}. 
\end{align*}
Using an argument similar to Lemma \ref{rough_exp_martingale_bounded} for the first term in the last expression and the results in Lemma \ref{rough_exp_martingale_bounded} and its corollary we deduce
\[
\left| \tilde{\mathfrak p}_t^{\phi}(\hat{\mathbf{Y}}(\omega)) - \tilde{\mathfrak p}_t^{\phi}(\hat{\mathbf{Y}}(\tilde{\omega})) \right| \leq C_{Y, W, M, p, T, b, \sigma, \phi} \|\hat{\mathbf{Y}}(\omega)  - \hat{\mathbf{Y}}(\tilde{\omega})\|_{p; [0, t]},
\]
which gives the desired result.
\end{proof}

\begin{proposition}
Suppose that Conditions~\ref{condition_covariance_rho_var},
\ref{cameron_martin_condition}, \ref{condition_filtration}, and
\ref{vf_condition_section_well_posedness} hold. Then for every $t \in [0, T] $  the identity $\mathfrak p_t^{\phi}(\hat{\mathbf{Y}}) = \tilde{\pi}_t(\phi)$ holds $\mathbb{P} $ and $\mathbb{Q}\,$ a.s.
\end{proposition}
\begin{proof}
We will prove that for any $\varsigma: \mathbb{C}_b(\mathbb{C}([0, t], \mathbb{R}^{d_Y}), \mathbb{R})$, $t \in [0,T]$, we have
\begin{equation}\label{robustness_identity}
\mathbb{E}^{\mathbb Q}\left[ \tilde{\pi}_t(\phi) \varsigma(Y_.)\right] = \mathbb{E}^{\mathbb Q}\left[ \tilde{\mathfrak p}_t^{\phi}(\hat{\mathbf{Y}}) \varsigma(Y_.)\right] ,
\end{equation}
with $Y. := \{Y_s : s\leq t\}$. From the construction of the decoupled probability space we see that
\[
\begin{pmatrix}
x_0 \\
B \\
Y 
\end{pmatrix}_{\mathbb{Q}} \stackrel{d}{\sim} \begin{pmatrix}
x_0 \\
\overline{B} \\
Y 
\end{pmatrix}_{\check{\mathbb{Q}}},
\]
which implies immediately from the continuity of the solution of RDE with respect to the driving rough path (see Theorem 4.20 in \cite{friz2018differential})
\[
\begin{pmatrix}
X \\
Y 
\end{pmatrix}_{\mathbb{Q}} \stackrel{d}{\sim} \begin{pmatrix}
\check{X} \\
Y 
\end{pmatrix}_{\check{\mathbb{Q}}},
\]
Now for the left side of equation \eqref{robustness_identity}, Proposition \ref{ito_conversion} yields
\begin{align} \label{lhs_robustness}
\begin{aligned}
    &\mathbb{E}^{\mathbb Q}\left[  \phi(X_t) \Lambda_t \varsigma(Y_.) \right] \\
    &= \mathbb{E}^{\check {\mathbb Q}}\left[  \phi(\check{X}_t) \check{\Lambda}_t \varsigma(Y_.) \right] \\
    &= \mathbb{E}^{\check {\mathbb Q}}\left[  \phi(\check{X}_t) \exp(\check{\Xi}_t) \varsigma(Y_.) \right],
\end{aligned}
\end{align}
where $\check{\Lambda}$ denotes the decoupled version of $\Lambda$.

For the right side of equation \eqref{robustness_identity}, by using the definition of $\tilde{ \mathfrak p}$ first and the definition of decoupled space one can see that
\begin{align}\label{rhs_robustness}
\begin{aligned}
& \mathbb{E}^{\mathbb Q}\left[ \mathfrak p_t^{\phi}(\hat{\mathbf{Y}})\varsigma(Y_.) \right] \\
&= \mathbb{E}^{\mathbb Q}\left[ \mathbb{E}^{\overline{\mathbb Q}}\left[ \phi(\check{X}_t) \exp(\check{\Xi}_t)\right] \varsigma(Y_.) \right] \\
&= \mathbb{E}^{\check {\mathbb Q}}\left[  \phi(\check{X}_t) \exp(\check{\Xi}_t) \varsigma(Y_.) \right].
\end{aligned}
\end{align}
Comparing the final expression obtained in \eqref{lhs_robustness} and \eqref{rhs_robustness} we obtain the desired result by equivalence of the probability measures $\mathbb{Q}$ and $\mathbb{P}$.

\end{proof}
We may now conclude the section by identifying the nonlinear filter with its robust pathwise representation.
\begin{corollary}[Existence of the robust filter]
In the setting of the preceding proposition, for every $t \in [0, T] $  the identity $\mathfrak p^\phi_t(\hat{\mathbf{Y}}) = \pi_t(\phi)$ holds $\mathbb{P} $ and $\mathbb{Q}$ a.s.
\end{corollary}

\section{The Zakai equation}
\label{section_zakai}

The purpose of this section is to pass from the robust representation of the
unnormalised filter to an equation for its density. More precisely, under the assumption that $X, Y$ and $B$ are one-dimensional, we study the
unnormalised conditional measure
\[
    \tilde{\pi}_t(\phi)
    =
    \mathbb E^{\mathbb Q}\left[
        \phi(X_t)\Lambda_t
        \,\middle|\,
        \mathcal Y
    \right],
\]
and we show that, whenever this measure admits a density, this density is
characterised by a rough analogue of the Zakai equation. We recall that in the case where $I$ is a Brownian motion, then the Zakai equation reads

\begin{align}\label{classical_mortensen_zakai_stratonovich}
\begin{aligned}
\tilde{\varrho}(t,x)
&=
\tilde{\varrho}_0(x)
+
\int_0^t
A_{Y_s}^*\tilde{\varrho}(s,x)\,ds
+
\int_0^t
f(x,Y_s)\tilde{\varrho}(s,x)\circ dY_s
\\
&\quad
-
\frac12
\int_0^t
\left(
    \partial_y f(x,Y_s)
    +
    f(x,Y_s)^2
\right)
\tilde{\varrho}(s,x)\,ds .
\end{aligned}
\end{align}
where the stochastic integral is understood in the Stratonovich sense and
\(A^*_y\) denotes the formal adjoint of the generator associated with the
Stratonovich signal
\[
    A_y^*\varphi(x)
    =
    \frac12
    \partial_{x,x}
    \left(
        \sigma(\cdot,y)^2\varphi
    \right)(x)
    -
    \frac12
    \partial_x
    \left(
        \sigma(\cdot,y)
        \partial_x\sigma(\cdot,y)
        \varphi
    \right)(x).
\]

There are two distinct features which have to be kept separate. First, the
signal is driven by a Volterra Gaussian process rather than by a semimartingale.
Consequently, the second-order part of the density equation is not generated by
the Lebesgue clock \(dt\), but by the diagonal covariance clock
\(c_t=C_B(t,t)\). Second, the observation enters through the likelihood
process \(\Lambda\), and, after conditioning, the resulting equation is driven
pathwise by the enhanced observation rough path \(\hat{\mathbf Y}\). By
construction, the \(Y\)-\(W\) cross-iterated integral in \(\hat{\mathbf Y}\)
is defined in the It\^o sense. The dependence of \(f\) on \(Y\) is therefore
already encoded in the chosen enhancement, and no additional
\(-\frac12\partial_y f\) correction appears. The compensator
\(-\frac12 f^2\), originating from the It\^o likelihood, remains explicit.
Thus, the density equation contains both a parabolic component, governed by
the covariance of the signal, and a rough multiplicative component, governed
by the observation.

We impose the following additional regularity assumption on the covariance
clock.

\begin{condition}\label{condition_7}
Let
\(
    c_t:=C_B(t,t),
    t\in[0,T].
\)
Assume that
\(
    c\in C([0,T],\mathbb R)\cap C^1((0,T],\mathbb R).
\)
Moreover, for every \(\delta>0\), there exists \(\epsilon_\delta>0\) such
that, for all \(a,b\in[0,\epsilon_\delta]\), the map
\[
    t\to C_B(t-a,t-b)
\]
belongs to the Sobolev space \(W^{1,1}([\delta,T])\), and
\(
    \left\|
        \frac{d}{dt}C_B(\,\cdot-a,\cdot-b\,)-\dot c
    \right\|_{L^1([\delta,T])}
    \to 0
    \text{ as }(a,b)\to(0,0).
\)
\end{condition}
We defer to Appendix \ref{app:examples} the proof that Type I and Type II fractional Brownian motion satisfy Condition \ref{condition_7}.

This assumption ensures that the diagonal covariance has sufficient regularity
to define the deterministic second-order part of the density equation. On the
set \(\Omega_4\) introduced in Lemma
\ref{rough_exp_martingale_bounded}, the formal candidate for the density
equation is the following Rough Partial Differential Equation (RPDE):
\begin{align}\label{rough_mortensen_zakai}
\begin{aligned}
\tilde{\varrho}(t,x)
&=
\tilde{\varrho}_0(x)
+
\int_0^t A_{Y_s}^*\tilde{\varrho}(s,x)\,dc_s
+
\int_0^t
\mathscr f(x,Y_s)\tilde{\varrho}(s,x)\,d\hat{\mathbf Y}_s
\\
&\quad
-
\frac12
\int_0^t
\bigl(f(x,Y_s)\bigr)^2
\tilde{\varrho}(s,x)\,ds ,
\end{aligned}
\end{align}
where
\(
    \mathscr f(x,y)
    :=
    \begin{pmatrix}
        0\\
        f(x,y)
    \end{pmatrix}.
\)

The equation above should not be viewed merely as a formal SPDE. Our
construction proceeds by approximation. Let \(\hat J^\epsilon\) be a smooth
approximation of the Volterra Gaussian driver, and let
\((X^\epsilon,Y^\epsilon)\) be the corresponding smooth signal-observation
system. In this regularised setting, the unnormalised filter is classical, and as we will show,
its density \(\tilde{\varrho}^\epsilon\) solves
\begin{align}\label{smooth_mortensen_zakai}
\begin{aligned}
\tilde{\varrho}^\epsilon(t,x)
&=
\tilde{\varrho}^\epsilon_0(x)
+
\int_0^t A_{Y^\epsilon_s}^*\tilde{\varrho}^\epsilon(s,x)\,dc^\epsilon(s,s)
+
\int_0^t
f(x,Y^\epsilon_s)\tilde{\varrho}^\epsilon(s,x)\,dW^\epsilon_s
\\
&\quad
-
\frac12
\int_0^t
\bigl(f(x,Y^\epsilon_s)\bigr)^2
\tilde{\varrho}^\epsilon(s,x)\,ds .
\end{aligned}
\end{align}

For fixed \(\epsilon>0\) and \(\omega\in\Omega_4\), the smoothed equation is a deterministic linear second-order parabolic equation with smooth time-dependent coefficients. Indeed, since \(c^\epsilon\), \(Y^\epsilon\), and \(W^\epsilon\) are smooth, it can be rewritten in differential form as
\[
    \partial_t \tilde\varrho^\epsilon(t,x)
    =
    \dot \partial_t c^\epsilon_t
    A^*_{Y_t^\epsilon}\tilde\varrho^\epsilon(t,x)
    +
    \left(
        f(x,Y_t^\epsilon)\dot W_t^\epsilon
        -
        \frac12 f(x,Y_t^\epsilon)^2
    \right)
    \tilde\varrho^\epsilon(t,x),
    \qquad
    \tilde\varrho^\epsilon(0,\cdot)=\tilde\varrho_0 .
\]
Since the equation may be degenerate, the natural notion of solution is that
of viscosity solution. More precisely, a bounded continuous function
\(\tilde\varrho^\epsilon\) is a viscosity subsolution if, whenever
\(\varphi\in C^{1,2}((0,T)\times\mathbb R)\) and
\(\tilde\varrho^\epsilon-\varphi\) has a local maximum at
\((t_0,x_0)\in(0,T)\times\mathbb R\), one has
\[
    \partial_t\varphi(t_0,x_0)
    \leq
    \dot c^\epsilon_{t_0}
    A^*_{Y_{t_0}^\epsilon}\varphi(t_0,x_0)
    +
    \left(
        f(x_0,Y_{t_0}^\epsilon)\dot W_{t_0}^\epsilon
        -
        \frac12 f(x_0,Y_{t_0}^\epsilon)^2
    \right)
    \varphi(t_0,x_0).
\]
It is a viscosity supersolution if the reverse inequality holds whenever
\(\tilde\varrho^\epsilon-\varphi\) has a local minimum at
\((t_0,x_0)\), and it is a viscosity solution if it is both a subsolution and
a supersolution, with the initial condition attained continuously at \(t=0\). Under boundedness and smoothness assumptions on the coefficients,
and assuming \(\tilde\varrho_0\in BUC(\mathbb R) \cap L^1(\mathbb R)\), standard comparison
and stability results for degenerate parabolic equations yield existence and
uniqueness of a bounded continuous viscosity solution; see, for instance,
Crandall, Ishii and Lions \cite{crandall1992users}. Thus, for every fixed
\(\epsilon>0\), the smoothed equation admits a unique bounded continuous
viscosity solution.

Following \cite{friz2010rough}, the rough equation \eqref{rough_mortensen_zakai} is then obtained as the 
limit of the smooth equations \eqref{smooth_mortensen_zakai}. This approximation
procedure serves two purposes. It provides a well-posedness theory for the
limiting RPDE, and it identifies its solution with the density of the
unnormalised conditional measure \(\tilde{\pi}_t\). In this sense, the density
equation is not postulated independently of the filtering problem; it is the
pathwise limit of the corresponding smooth Zakai equations.

More precisely, the next definition introduces the class of approximating paths
which will be used to define the limiting RPDE.

\begin{definition}\label{approximation_class_rpde}
    We say that
    \(
        \hat J^\epsilon \subset
        C^1([0,T],\mathbb R^3)
    \)
    is a smooth approximating sequence for \(\hat J\) if the following two
    conditions hold.

    First, the canonical enhancement
    \[
        \mathscr S_{\kappa}(\hat  J ^\epsilon)
        =
        \left(
            1,
            \hat J^\epsilon_{s,t},
            \int_{s<u_1<u_2<t}
                d\hat J^\epsilon_{u_1}\otimes d\hat J^\epsilon_{u_2},
            \ldots,
            \int_{s<u_1<\cdots<u_\kappa<t}
                d\hat J^\epsilon_{u_1}\otimes\cdots\otimes d\hat J^\epsilon_{u_\kappa}
        \right),
    \]
     satisfies
    \(
        \mathscr S_{\kappa}(\hat J ^\epsilon)
        \to
        \hat{\mathbf J}
    \)
    $\mathbb P$-a.s. in the \(p\)-rough path topology.

    Second, $c^\epsilon \in \mathcal C^1([0, T], \mathbb{R})$, 
    \(
        c^\epsilon \to c     
        \text{ uniformly on }[0,T],
    \)
    and, for every \(\delta>0\),
    \(
        \dot c^\epsilon \to \dot c
        \text{ in } L^1([\delta,T]).
    \)
\end{definition}

We defer the example showing how to construct such a class under Condition \ref{condition_7} to Appendix \ref{app:examples}.

The reason for introducing this approximation class is that the approximation
problem here is slightly different from the standard one appearing in rough
viscosity theory (see \cite{friz2010rough}). In the usual rough viscosity framework, the second-order
operator is fixed, and only the rough driving signal is approximated. Thus the
main stability requirement is convergence of the smooth lifts to the limiting
rough path.

In the present problem, however, smoothing the Volterra Gaussian driver changes
two objects at the same time. It changes the rough driver, through
\(
    \mathscr S_{\kappa}(\hat J^\epsilon),
\)
but it also changes the deterministic parabolic clock appearing in the
second-order part of the density equation. Indeed, the regularised equation
contains the term
\(
    \int_0^t A^*\tilde{\varrho}^\epsilon(s,\cdot)\,dc^\epsilon_s
    =
    \int_0^t
        \dot c^\epsilon_s
        A^*\tilde{\varrho}^\epsilon(s,\cdot)\,ds,
\)
whereas the limiting equation contains
\(
    \int_0^t A^*\tilde{\varrho}(s,\cdot)\,dc_s .
\)
Therefore, convergence of
\(\mathscr S_{\kappa}(\hat J^\epsilon)\) to \(\hat{\mathbf J}\) controls the
rough part of the equation, but it does not by itself identify the limiting
second-order operator.

For this reason, we require both rough path convergence of the enhanced smooth
paths and convergence of the covariance clocks. The uniform convergence
\(c^\epsilon\to c\) identifies the integrated parabolic clock, while the local
$L^1$ convergence \(\dot c^\epsilon\to \dot c\) on \([\delta,T]\), for every
\(\delta>0\), allows one to apply classical stability estimates for the
regularised parabolic equations away from the origin. This local formulation is
natural for Volterra Gaussian processes, since the derivative of the diagonal
covariance may be singular at \(0\), even though the clock \(c\) itself is
continuous on \([0,T]\).
We are now ready to provide the definition for rough viscosity solution for the Zakai equation associated to our problem

\begin{definition}\label{rough_viscosity_solution}
    We say that the continuous functional
    \(\tilde{\varrho}\) solves
    \eqref{rough_mortensen_zakai} if, for every sequence of smooth approximating paths
    \((\hat J^\epsilon)_{\epsilon>0}\) that satisfies Definition \ref{approximation_class_rpde}, the unique solutions
    \(\tilde{\varrho}^\epsilon\) to the associated smooth
    problems \eqref{smooth_mortensen_zakai} converge to
    \(\tilde{\varrho}\) locally uniformly on
    \([0,T]\times\mathbb R\).
    Uniqueness is understood in the class of continuous functionals obtained
    as such locally uniform limits.
\end{definition}

The next theorem  adapts Theorem 1 of \cite{friz2010rough} to the case
where the approximation affects not only the rough driver, but also the
deterministic covariance clock appearing in the second-order part of the
equation. This allows us to pass from the family of smoothed viscosity solutions
\((\tilde\varrho^\epsilon)_{\epsilon>0}\) to a well-defined limiting object,
which we then identify as the solution of the rough Zakai equation.

Before stating and proving the theorem we need to impose an additional condition of $f, \sigma$ and $\varrho_0$ that in this section will replace \ref{vf_condition_section_well_posedness}.

\begin{condition}\label{cond:vector_field_reg_zakai}
    The coefficients satisfy
    \(
        f,\sigma\in C_b^{\kappa+3}(\mathbb R^2).
    \)
    Moreover,
    \(
        \tilde{\varrho}_0\in \mathrm{BUC}(\mathbb R)\cap L^1(\mathbb R),
    \)
    where \(\mathrm{BUC}(\mathbb R)\) denotes the space of bounded uniformly
    continuous functions on \(\mathbb R\).
\end{condition}

\begin{theorem}\label{thm:lifted_rough_mz_wellposedness}
Assume that conditions \ref{condition_covariance_rho_var}, \ref{condition_filtration}, \ref{condition_7} and \ref{cond:vector_field_reg_zakai} hold.
Then there exists a
 unique rough viscosity solution $\tilde \varrho$ of the Zakai
equation in the sense of definition \ref{rough_viscosity_solution}. 
Moreover,
\[
    \sup_{t\in[0,T]}
    \|\tilde{\varrho}(t,\cdot)\|_\infty
    \leq
    C\|\tilde{\varrho}_0\|_\infty,
\]
for $C$ dependent on $f, \sigma, p, T$. 
\end{theorem}
\begin{proof}
Let $\hat J^\epsilon = (B^\epsilon, Y^\epsilon, W^\epsilon)$ denote a path in the class defined in Definition \ref{approximation_class_rpde}.
Consider the equation
\[
\begin{aligned}
    \partial_t P^\epsilon(t,x,y)
    &=
    \dot c_t^\epsilon A_y^*P^\epsilon(t,x,y)
    -
    \frac12 f(x,y)^2P^\epsilon(t,x,y)        \\
    &\quad
    -
    \partial_yP^\epsilon(t,x,y)\dot Y_t^\epsilon
    +
    f(x,y)P^\epsilon(t,x,y)\dot W_t^\epsilon,
    \qquad
    P^\epsilon(0,x,y)= \tilde{\varrho}_0(x),
\end{aligned}
\]
where we recall that \(c^\epsilon(0)=0\) and
\(\dot c^\epsilon\geq0\). We refer to this equation as the smooth lifted equation in the $y$-variable. 

We first rewrite the deterministic part in viscosity form. Expanding \(A_y^*\),
we have
\[
    A_y^*\phi
    =
    \partial_{xx}(\sigma^2\phi)
    -
    \partial_x(\sigma\sigma_x\phi)
    =
    \sigma^2\phi_{xx}
    +
    3\sigma\sigma_x\phi_x
    +
    \chi \phi,
\]
where
\[
    \chi (x,y)
    :=
    \sigma_x(x,y)^2+\sigma(x,y)\sigma_{xx}(x,y).
\]
Writing \(z=(x,y)\), \(\mathscr p=(\mathscr p_x,\mathscr p_y)\), and
\(\Upsilon=(\Upsilon_{ij})_{i,j=1}^2\), the smooth lifted equation becomes
\[
    \partial_t P^\epsilon
    +
    A^\epsilon(t,z,P^\epsilon,DP^\epsilon,D^2P^\epsilon)
    =
    \Lambda_1(z,P^\epsilon,DP^\epsilon)\dot Y_t^\epsilon
    +
    \Lambda_2(z,P^\epsilon,DP^\epsilon)\dot W_t^\epsilon,
\]
with
\[
\begin{aligned}
    A^\epsilon(t,z,\xi,\mathscr p,\Upsilon)
    &=
    -\dot c_t^\epsilon\sigma(x,y)^2\Upsilon_{xx}
    -
    3\dot c_t^\epsilon\sigma(x,y)\sigma_x(x,y)\mathscr p_x        \\
    &\quad
    +
    \left(
        \frac12 f(x,y)^2
        -
        \dot c_t^\epsilon\chi(x,y)
    \right)\xi,
\end{aligned}
\]
and
\[
    \Lambda_1(z,\xi,p)=-\mathscr p_y,
    \qquad
    \Lambda_2(z,\xi,p)=f(z)\xi .
\]
The operator \(A^\epsilon\) is degenerate elliptic because
\(
    \dot c_t^\epsilon\sigma(x,y)^2\geq0.
\)

We now apply the transform used in Theorem 13 in \cite{friz2010rough} to remove the terms $\dot Y^\epsilon$ and $\dot W^\epsilon$.
For each \(\epsilon\), there exist a smooth flow
\(\Phi_t^\epsilon\) and a smooth multiplier \(m_t^\epsilon\) such that
\begin{alignat*}{2}
d\Phi^\epsilon_t(z)
&= -\dot{Y}^\epsilon_tdt,
\qquad
&\Phi^\epsilon_0(z)
&= z,\\
dm^\epsilon_t
&= m^\epsilon_t f\bigl(z,\Phi^\epsilon_t\bigr)dW^\epsilon_t,
\qquad
&m^\epsilon_0
&= 1.
\end{alignat*}
and 
\[
    P^\epsilon(t,z)
    =
    \Theta_t^\epsilon v^\epsilon(t,z)
    :=
    m_t^\epsilon(z)\,
    v^\epsilon(t,\Phi_t^\epsilon(z))\quad v^\epsilon(0,z)=\tilde{\varrho}_0(x).
\]
The transformed unknown \(v^\epsilon\) solves a deterministic viscosity
equation of the form
\[
    \partial_t v^\epsilon
    +
    G_0^\epsilon(t,z,v^\epsilon,Dv^\epsilon,D^2v^\epsilon)
    +
    \dot c_t^\epsilon
    G_1^\epsilon(t,z,v^\epsilon,Dv^\epsilon,D^2v^\epsilon)
    =
    0.
\]
The transformed operators satisfy the same structural assumptions as in the
Theorem 1 in \cite{friz2010rough}. In particular, for \(\xi\geq\zeta\),
\[
\begin{aligned}
    &G_i^\epsilon(t,z,\xi,\mathscr p,\Upsilon)
    -
    G_i^\epsilon(t,z,\zeta,\mathscr p,\Upsilon)        \\
    &\qquad\geq
    -\lambda_i(\xi-\zeta),
    \qquad
    i=0,1,
\end{aligned}
\]
with constants \(\lambda_0,\lambda_1\) independent of \(\epsilon\).
We next prove the initial-layer estimate. Let
\[
    \psi_0(z):=\tilde{\varrho}_0(x),
    \qquad z=(x,y).
\]
Fix \(\delta>0\). Since \(\tilde{\varrho}_0\in BUC(\mathbb R) \cap L^1(\mathbb R)\), there exists
\(\phi_0\in C_b^\infty(\mathbb R)\) such that
\[
    \|\phi_0-\tilde{\varrho}_0\|_\infty\leq\delta.
\]
Set
\(
    \phi(x,y):=\phi_0(x).
\)
Then \(\phi\in C_b^\infty(\mathbb R^2)\) and
\[
    \|\phi-\psi_0\|_\infty\leq\delta.
\]
Since \(\phi\) has bounded derivatives and the transformed coefficients are
bounded, there exists \(C_\phi>0\), independent of \(\epsilon\), such that
\[
    \left|
        G_0^\epsilon(t,z,\phi,D\phi,D^2\phi)
    \right|
    +
    \left|
        G_1^\epsilon(t,z,\phi,D\phi,D^2\phi)
    \right|
    \leq C_\phi
\]
for all \(t\in[0,T]\), \(z\in\mathbb R^2\), and all sufficiently small
\(\epsilon\).
Define
\[
    M^\epsilon(t)
    :=
    \int_0^t\lambda^\epsilon(s)\,ds
    =
    \lambda_0t+\lambda_1c^\epsilon(t),
\]
and
\[
    N^\epsilon(t)
    :=
    e^{M^\epsilon(t)}\delta
    +
    C_\phi
    \int_0^t
    e^{M^\epsilon(t)-M^\epsilon(s)}
    \bigl(ds+dc_s^\epsilon\bigr).
\]
Then
\[
    \overline v^\epsilon(t,z):=\phi(z)+N^\epsilon(t)
\]
is a viscosity supersolution of the transformed deterministic equation, while
\[
    \underline v^\epsilon(t,z):=\phi(z)-N^\epsilon(t)
\]
is a viscosity subsolution. Indeed, \(N^\epsilon\) satisfies
\[
    dN^\epsilon_t
    =
    \lambda^\epsilon(t)N^\epsilon(t)\,dt
    +
    C_\phi\bigl(dt+dc_t^\epsilon\bigr),
    \qquad
    N^\epsilon(0)=\delta,
\]
which absorbs both the residual obtained by inserting \(\phi\) into the
equation and the possible negative monotonicity in the zeroth-order variable.
At \(t=0\), since \(N^\epsilon(0)=\delta\) and
\(\|\phi-\psi_0\|_\infty\leq\delta\), we have
\[
    \underline v^\epsilon(0,z)
    \leq
    v^\epsilon(0,z)
    \leq
    \overline v^\epsilon(0,z),
    \qquad
    z \in\mathbb R^2.
\]
By comparison principle (see Appendix in \cite{friz2010rough}),
\[
    \underline v^\epsilon(t,z)
    \leq
    v^\epsilon(t,z)
    \leq
    \overline v^\epsilon(t,z),
    \qquad
    t\in[0,T],\ z\in\mathbb R^2.
\]
Hence
\[
    |v^\epsilon(t,z)-\psi_0(z)|
    \leq
    \delta+N^\epsilon(t),
    \qquad
    t\in[0,T],\ z\in\mathbb R^2.
\]
For \(t\leq\eta\), we have
\[
    N^\epsilon(t)
    \leq
    e^{\lambda_0\eta+\lambda_1c^\epsilon_\eta}\delta
    +
    C_\phi
    e^{\lambda_0\eta+\lambda_1c^\epsilon_\eta}
    \bigl(\eta+c^\epsilon_\eta\bigr).
\]
Since \(c^\epsilon\to c\) uniformly and \(c_0=0\),
\(
    \lim_{\eta\downarrow0}
    \limsup_{\epsilon\downarrow0}
    c^\epsilon_\eta=0.
\)
Therefore
\[
    \lim_{\eta\downarrow0}
    \limsup_{\epsilon\downarrow0}
    \sup_{0\leq t\leq\eta}
    \sup_{z\in\mathbb R^2}
    |v^\epsilon(t,z)-\psi_0(z)|
    \leq 2\delta.
\]
Letting \(\delta\downarrow0\), we obtain
\[
    \lim_{\eta\downarrow0}
    \limsup_{\epsilon\downarrow0}
    \sup_{0\leq t\leq\eta}
    \sup_{z\in\mathbb R^2}
    |v^\epsilon(t,z)-\psi_0(z)|
    =
    0.
\]
We now transfer this estimate back to \(P^\epsilon\). Let
\(\mathcal K\subset\mathbb R^2\) be compact. By continuity of the rough flow and by
the convergence of the smooth lifts to \(\hat{\mathbf J}\), there exists a compact
set \(\mathcal K_\eta\subset\mathbb R^2\) such that, for all sufficiently small \(\eta>0\)
and all sufficiently small \(\epsilon\),
\[
    \Phi_t^\epsilon(\mathcal K)\subset \mathcal K_\eta,
    \qquad
    0\leq t\leq\eta.
\]
Moreover,
\[
    a^\epsilon(\eta)
    :=
    \sup_{0\leq t\leq\eta}\sup_{z\in \mathcal K}
    \left(
        |\Phi_t^\epsilon(z)-z|
        +
        |m_t^\epsilon(z)-1|
    \right)
\]
satisfies
\(
    \lim_{\eta\downarrow0}
    \limsup_{\epsilon\downarrow0}
    a^\epsilon(\eta)=0.
\)
The multipliers are also uniformly bounded on \(\mathcal K\) for small times: there is
\(C_{\mathcal K}<\infty\) such that
\[
    \sup_{\epsilon}
    \sup_{0\leq t\leq\eta}
    \sup_{z\in \mathcal K}
    |m_t^\epsilon(z)|
    \leq C_{\mathcal K}
\]
for \(\eta\) sufficiently small.

For \(z\in \mathcal K\), we estimate
\[
\begin{aligned}
    |P^\epsilon(t,z)-\psi_0(z)|
    &=
    \left|
        m_t^\epsilon(z)
        v^\epsilon(t,\Phi_t^\epsilon(z))
        -
        \psi_0(z)
    \right|        \\
    &\leq
    |m_t^\epsilon(z)|
    \left|
        v^\epsilon(t,\Phi_t^\epsilon(z))
        -
        \psi_0(\Phi_t^\epsilon(z))
    \right|        \\
    &\quad
    +
    |m_t^\epsilon(z)|
    \left|
        \psi_0(\Phi_t^\epsilon(z))-\psi_0(z)
    \right|        \\
    &\quad
    +
    |m_t^\epsilon(z)-1|\,|\psi_0(z)|.
\end{aligned}
\]
Since \(\psi_0\) is uniformly continuous, this gives
\[
\begin{aligned}
    \sup_{0\leq t\leq\eta}\sup_{z\in \mathcal K}
    |P^\epsilon(t,z)-\psi_0(z)|
    &\leq
    C_{\mathcal K}
    \sup_{0\leq t\leq\eta}
    \sup_{z\in \mathcal K_R}
    |v^\epsilon(t,z)-\psi_0(z)|        \\
    &\quad
    +
    C_{\mathcal K}\varpi_{\psi_0}(a^\epsilon(\eta))
    +
    \|\psi_0\|_{\infty;K}a^\epsilon(\eta),
\end{aligned}
\]
where \(\varpi_{\psi_0}\) is a modulus of continuity of \(\psi_0\). Taking
\(\limsup_{\epsilon\downarrow0}\), then \(\eta\downarrow0\), we obtain
\[
    \lim_{\eta\downarrow0}
    \limsup_{\epsilon\downarrow0}
    \sup_{0\leq t\leq\eta}
    \sup_{z\in K}
    |P^\epsilon(t,z)-\tilde{\varrho}_0(x)|
    =
    0.
\]
This is the required initial-layer estimate for the lifted solutions.
We also record the uniform bound. Let
\(
    C_0:=\|\tilde{\varrho}_0\|_\infty.
\)
Since the transformed equation is homogeneous in $v^\epsilon$ and \(0\) is a solution, the
functions
\[
    \overline v^\epsilon(t,z)
    :=
    C_0e^{M^\epsilon(t)},
    \qquad
    \underline v^\epsilon(t,z)
    :=
    -C_0e^{M^\epsilon(t)}
\]
are respectively a supersolution and a subsolution. Additionally
\[
    \underline v^\epsilon(0,z)
    \leq
    v^\epsilon(0,z)
    \leq
    \overline v^\epsilon(0,z),
\]
so that comparison principle gives
\[
    |v^\epsilon(t,z)|
    \leq
    \|\tilde{\varrho}_0\|_\infty
    \exp(\lambda_0t+\lambda_1c^\epsilon_t).
\]
The multiplier \(m_t^\epsilon\) is uniformly bounded on
\([0,T]\times\mathbb R^2\), with a bound depending on the rough path control
of \(\hat{Y}^\epsilon\) and on the bounded coefficient \(f\). Since
\(\mathscr S_{\kappa }(\hat Y^\epsilon)\to  \hat{\mathbf Y}\), these bounds are uniform
for small \(\epsilon\). Thus there exists \(C_\Theta(T,\hat{\mathbf Y})<\infty\)
such that
\[
    \sup_{t\in[0,T]}\|P^\epsilon(t,\cdot)\|_\infty
    \leq
    C_\Theta(T,\hat{\mathbf Y})
    \exp(\lambda_0T+\lambda_1c^\epsilon_T)
    \|u_0\|_\infty .
\]
Since \(c^\epsilon\to c\) uniformly,
\(
    \sup_{\epsilon}c^\epsilon(T)<\infty.
\)
Consequently, for some finite constant \(C\), independent of \(\epsilon\),
\[
    \sup_{\epsilon}
    \sup_{t\in[0,T]}
    \|P^\epsilon(t,\cdot)\|_\infty
    \leq
    C\|\tilde{\varrho}_0\|_\infty.
\]
We now prove convergence. Let \(\mathcal K\subset\mathbb R^2\) be compact and fix
\(a>0\). From condition \ref{condition_7}, on every strip \([a,T]\), the deterministic clock densities converge
in the required sense
\(
    \dot c^\epsilon\to\dot c
    \text{ in }L^1([a,T]).
\)
Moreover, the canonical lifts of \(\hat{Y}^\epsilon\) converge to
\(\hat{\mathbf Y}\) in rough path topology. Therefore
Theorem 1 in \cite{friz2010rough}, together with the deterministic stability
of viscosity solutions under \(L^1\)-in-time perturbations of the coefficients (Theorem 1.1 in \cite{barles2006new}),
implies stability on \([a,T]\). More precisely, for \(0<\eta<a\),
\[
\begin{aligned}
    &\sup_{t\in[a,T]}\sup_{z\in \mathcal K}
    |P^\epsilon(t,z)-P^{\epsilon'}(t,z)|        \\
    &\qquad
    \leq
    C
    \sup_{q\in \mathcal K_R}
    |P^\epsilon(\eta,z)-P^{\epsilon'}(\eta,z)|
    +
    o_{\epsilon,\epsilon'}(1),
\end{aligned}
\]
where \(C\) is independent of \(\eta,\epsilon,\epsilon'\), and
\(o_{\epsilon,\epsilon'}(1)\to0\) as
\(\epsilon,\epsilon'\downarrow0\). The constant \(C\) depends on
\(T\), \(\mathcal K\), the rough path control, and the integrated properness bound
\[
    \sup_\epsilon
    \int_0^T
    \lambda^\epsilon(s)\,ds
    =
    \sup_\epsilon
    \left(
        \lambda_0T+\lambda_1c^\epsilon_T
    \right),
\]
but not on \(\|\dot c^\epsilon\|_\infty\). By the initial-layer estimate,
\[
    \lim_{\eta\downarrow0}
    \limsup_{\epsilon,\epsilon'\downarrow0}
    \sup_{q\in K_R}
    |P^\epsilon(\eta,q)-P^{\epsilon'}(\eta,q)|
    =
    0.
\]
Therefore \((P^\epsilon)_{\epsilon > 0 }\) is Cauchy locally uniformly on
\((0,T]\times\mathbb R^2\). Let \(P\) denote the locally uniform limit there.
The initial-layer estimate also implies that \(P\) extends continuously to
\(t=0\) by setting
\[
    P(0,x,y)=\tilde{\varrho}_0(x).
\]
Thus
\[
    P^\epsilon\to P
    \qquad
    \text{locally uniformly on }[0,T]\times\mathbb R^2.
\]
Passing to the limit in the uniform bound gives
\[
    \sup_{t\in[0,T]}
    \|P(t,\cdot)\|_\infty
    \leq
    C\|\tilde{\varrho}_0\|_\infty.
\]
On every strip \([a,T]\), the limit \(P\) is the rough
viscosity solution of the lifted equation with deterministic coefficient
\(\dot c\). Since \(a>0\) is arbitrary and \(P\) has the initial trace
\(\tilde{\varrho}_0\), \(P\) is the lifted rough viscosity solution of
\[
\begin{aligned}
    dP(t,x,y)
    &=
    A_y^*P(t,x,y)\,dr_t
    -
    \frac12 f(x,y)^2P(t,x,y)\,dt        \\
    &\quad
    -
    \partial_yP(t,x,y)\,d\mathbf Y_t
    +
    \mathscr f(x,y)P(t,x,y)\,d\hat{\mathbf {Y}}_t,
    \qquad
    P(0,x,y)=\tilde{\varrho}_0(x).
\end{aligned}
\]
Define
\[
    \tilde{\varrho}^\epsilon(t,x):=P^\epsilon(t,x,Y_t^\epsilon).
\]
Since \(\hat Y^\epsilon\) is smooth, the classical chain rule gives
\[
    \partial_t\tilde{\varrho}^\epsilon(t,x)
    =
    \partial_tP^\epsilon(t,x,Y_t^\epsilon)
    +
    \partial_yP^\epsilon(t,x,Y_t^\epsilon)\dot Y_t^\epsilon.
\]
Substituting the smooth lifted equation, the \(Y^\epsilon\)-transport terms
cancel. Hence
\[
\begin{aligned}
    \partial_t\tilde{\varrho}^\epsilon(t,x)
    &=
    \dot c_t^\epsilon
    A_{Y_t^\epsilon}^*\tilde{\varrho}^\epsilon(t,x)
    -
    \frac12 f(x,Y_t^\epsilon)^2\tilde{\varrho}^\epsilon(t,x)        \\
    &\quad
    +
    f(x,Y_t^\epsilon)\tilde{\varrho}^\epsilon(t,x)\dot W_t^\epsilon.
\end{aligned}
\]
Equivalently,
\[
\begin{aligned}
    \tilde{\varrho}^\epsilon(t,x)
    &=
    \tilde{\varrho}_0(x)
    +
    \int_0^t
    A_{Y_s^\epsilon}^*
    \tilde{\varrho}^\epsilon(s,x)\,dc_s^\epsilon        \\
    &\quad
    -
    \frac12
    \int_0^t
    f(x,Y_s^\epsilon)^2\tilde{\varrho}^\epsilon(s,x)\,ds
    +
    \int_0^t
f(x,Y_s^\epsilon)\tilde{\varrho}^\epsilon(s,x)\,dW_s^\epsilon .
\end{aligned}
\]

Define
\(
    \tilde{\varrho}(t,x):=P(t,x,Y_t).
\)
Rough path convergence implies, in particular,
\[
    Y^\epsilon\to Y
    \qquad
    \text{uniformly on }[0,T].
\]
Let \(\mathcal K\subset\mathbb R\) be compact. Since \(Y^\epsilon\to Y\) uniformly
and \(Y\) is continuous, there exists a compact interval
\(\mathcal K_Y\subset\mathbb R\) such that
\[
    Y_t\in \mathcal K_Y,
    \qquad
    Y_t^\epsilon\in \mathcal K_Y,
    \qquad
    t\in[0,T],
\]
for all sufficiently small \(\epsilon\). Therefore
\[
\begin{aligned}
    &\sup_{(t,x)\in[0,T]\times \mathcal K}
    |\tilde{\varrho}^\epsilon(t,x)-\tilde{\varrho}(t,x)|        \\
    &\quad\leq
    \sup_{(t,x)\in[0,T]\times \mathcal K}
    |P^\epsilon(t,x,Y_t^\epsilon)
        -P(t,x,Y_t^\epsilon)|        \\
    &\qquad
    +
    \sup_{(t,x)\in[0,T]\times \mathcal K}
    |P(t,x,Y_t^\epsilon)-P(t,x,Y_t)|.
\end{aligned}
\]
The first term converges to zero by locally uniform convergence of
\(P^\epsilon\) to \(P\) on \([0,T]\times \mathcal K\times \mathcal K_Y\). The second term
converges to zero because \(P\) is uniformly continuous on this compact set
and \(Y^\epsilon\to Y\) uniformly. Hence
\[
    \tilde{\varrho}^\epsilon\to \tilde{\varrho}
    \qquad
    \text{locally uniformly on }[0,T]\times\mathbb R.
\]
Moreover,
\[
    \sup_{t\in[0,T]}
    \|\tilde{\varrho}(t,\cdot)\|_\infty
    \leq
    \sup_{t\in[0,T]}
    \|P(t,\cdot)\|_\infty
    \leq
    C\|\tilde{\varrho}_0\|_\infty.
\]

Finally, uniqueness follows from the same stability estimate. If two limits
are obtained from two admissible smooth approximation schemes, their
approximations can be compared on every strip \([a,T]\). The stability
estimate gives equality on \([a,T]\) after passing to the limit, while the
initial-layer estimate identifies the same initial trace as \(a\downarrow0\).
Thus the lifted solution \(P\), and consequently the trace
\(
    \tilde{\varrho}(t,x)=P(t,x,Y_t),
\)
is unique in the stated approximation class.
\end{proof}

We now identify the function \(\tilde{\varrho}\) constructed in the preceding theorem with the
density of the unnormalised filter. For each \(\epsilon>0\), the
path-and-covariance smoothed problem is classical, and its solution
\(\tilde{\varrho}^\epsilon\) is the density of the corresponding unnormalised conditional
measure. Hence, for every \(\phi\in C_b^\infty(\mathbb R)\),
\[
    \int_{\mathbb R}\phi(x)\tilde{\varrho}^\epsilon(t,x)\,dx
    =
    \tilde\pi_t^\epsilon(\phi).
\]
By the backward duality argument, this identity is stable under the
path-and-covariance approximation. Passing to the limit as \(\epsilon\to0\)
therefore yields
\[
    \int_{\mathbb R}\phi(x)u(t,x)\,dx
    =
    \tilde\pi_t(\phi),
    \qquad
    \phi\in C_b^\infty(\mathbb R),\quad t\in[0,T].
\]
It follows that the unnormalised conditional measure \(\tilde\pi_t\) is
absolutely continuous with respect to Lebesgue measure and has density
\(\tilde{\varrho}(t,\cdot)\).

We begin this identification by introducing the operator \(A\), dual to
\(A^*\) on test functions, and by defining the class of weak solutions which
characterises the unnormalised robust filter. By the identification above, this is equivalent to applying
the argument to the original filter therefore, for notational convenience we use  \(\tilde\pi_t\)  to identify the filter in this section.
For \(y\in\mathbb R\), define the operator acting on test functions by
\[
    A_y\phi(x)
    :=
    \sigma(x,y)^2\partial_{xx}\phi(x)
    +
    \sigma(x,y)\partial_x\sigma(x,y)\partial_x\phi(x),
    \qquad
    \phi\in C_c^\infty(\mathbb R).
\]
Equivalently, \(A_y\) is the formal adjoint of \(A_y^*\), in the sense that
\[
    \int_{\mathbb R}\phi(x)A_y^* g(x)\,dx
    =
    \int_{\mathbb R}A_y\phi(x) g(x)\,dx .
\]

\begin{definition}[Rough weak solution]\label{def:weak_measure_rough_zakai_approx}
A finite measure-valued path
\(
    \varkappa:[0,T]\to\mathcal M_f(\mathbb R)
\)
is called a rough weak solution of the rough Zakai equation if for every \(\phi\in C_c^\infty(\mathbb R)\), the path
\(\varkappa^\epsilon\), constructed using smooth approximating paths of $\hat J$, satisfies the weak Zakai equation
\begin{equation}\label{weak_smooth_zakai}
\begin{aligned}
    \varkappa_t^\epsilon(\phi)
    &=
    \varkappa_0^\epsilon(\phi)
    +
    \int_0^t
        \varkappa_s^\epsilon\bigl(A_{Y_s^\epsilon}\phi\bigr)\,
        dc^\epsilon(s)        \\
    &\quad
    +
    \int_0^t
        \varkappa_s^\epsilon
        \bigl(\mathscr f(\cdot,Y_s^\epsilon)\phi\bigr)\,
        d\hat Y_s^\epsilon
    -
    \frac12
    \int_0^t
        \varkappa_s^\epsilon
        \bigl(f(\cdot,Y_s^\epsilon)^2\phi\bigr)\,ds 
\end{aligned}
\end{equation}
and for every \(\phi\in C_c^\infty(\mathbb R)\),
\[
    \sup_{t\in[0,T]}
    \left|
        \varkappa_t^\epsilon(\phi)-\varkappa_t(\phi)
    \right|
    \to 0 .
\]
\end{definition}
We are now ready to state and prove the results which implement the backward
uniqueness argument for the rough weak equation and allow us to identify the
limiting viscosity solution with the unnormalised conditional density.
\begin{lemma}[Backward dual equation]\label{lem:backward_dual_equation}
Assume the hypotheses of Theorem
\ref{thm:lifted_rough_mz_wellposedness} hold. Let
\(t\in[0,T]\) and \(\psi\in C_c^\infty(\mathbb R)\). Then there exists a bounded
backward rough test function
\[
    \varphi^{t,\psi}:[0,t]\times\mathbb R\to\mathbb R
\]
with terminal condition
\(
    \varphi_t^{t,\psi}(x)=\psi(x),
\)
that solves 
\[
\begin{aligned}
    d\varphi_s^{t,\psi}(x)
    &=
    -A_{Y_s}\varphi_s^{t,\psi}(x)\,dr(s)
    -
    \mathscr f(x,Y_s)\varphi_s^{t,\psi}(x)\,d\hat{\mathbf Y}_s        \\
    &\quad
    +
    \frac12 f(x,Y_s)^2\varphi_s^{t,\psi}(x)\,ds,
    \qquad 0\leq s\leq t.
\end{aligned}
\]
in the sense of rough viscosity solutions i.e. for every sequence of smooth approximating paths
    \((\hat J^\epsilon)_{\epsilon>0}\) that satisfies, the unique solutions
    \(\varphi^{\epsilon, t, \psi}\) to the associated smooth
    problems 
    \[
\begin{aligned}
    d\varphi_s^{\epsilon, t,\psi}(x)
    &=
    -A_{Y_s}\varphi_s^{\epsilon,t,\psi}(x)\,dr(s)
    -
    f(x,Y_s)\varphi_s^{\epsilon,t,\psi}(x)\,dW^\epsilon_s        \\
    &\quad
    +
    \frac12 f(x,Y_s)^2\varphi_s^{\epsilon,t,\psi}(x)\,ds,
    \qquad 0\leq s\leq t.
\end{aligned}
\]
with terminal condition \(
    \varphi_t^{\epsilon, t,\psi}(x)=\psi(x),
\)
    converge to
    \(\varphi^{t, \psi}\) locally uniformly on
    \([0,T]\times\mathbb R\).
    
\end{lemma}

\begin{proof}
The equation is constructed by lifting once more in the \(y\)-variable. More
precisely, let \(\Psi^{t,\psi}(s,x,y)\) solve the lifted backward rough PDE
\[
\begin{aligned}
    d\Psi_s^{t,\psi}(x,y)
    &=
    -A_y\Psi_s^{t,\psi}(x,y)\,dr(s)
    +
    \frac12 f(x,y)^2\Psi_s^{t,\psi}(x,y)\,ds        \\
    &\quad
    -
    \partial_y\Psi_s^{t,\psi}(x,y)\,d\mathbf Y_s
    -
    \mathscr f(x,y)\Psi_s^{t,\psi}(x,y)\,d\hat{\mathbf{Y}}_s,
\end{aligned}
\]
with terminal condition
\(
    \Psi_t^{t,\psi}(x,y)=\psi(x).
\)
This backward problem is reduced to a forward RPDE by reversing time on
\([0,t]\). The resulting forward equation has the same structural form as the
lifted equation in Theorem \ref{thm:lifted_rough_mz_wellposedness}; the
combined-clock argument applies exactly as before. 

Define
\[
    \varphi_s^{t,\psi}(x)
    :=
    \Psi_s^{t,\psi}(x,Y_s).
\]
For smooth approximations, the ordinary chain rule gives
\[
\begin{aligned}
    d\Psi_s^{\epsilon,t,\psi}(x,Y_s^\epsilon)
    &=
    d\Psi_s^{\epsilon,t,\psi}(x,y)
        \big|_{y=Y_s^\epsilon}
    +
    \partial_y\Psi_s^{\epsilon,t,\psi}
        (x,Y_s^\epsilon)\,dY_s^\epsilon .
\end{aligned}
\]
The transport term
\[
    -\partial_y\Psi_s^{\epsilon,t,\psi}(x,y)\,dY_s^\epsilon
\]
in the lifted equation cancels the chain-rule term. Therefore the trace
\[
    \varphi_s^{\epsilon,t,\psi}(x)
    :=
    \Psi_s^{\epsilon,t,\psi}(x,Y_s^\epsilon)
\]
solves
\[
\begin{aligned}
    d\varphi_s^{\epsilon,t,\psi}(x)
    &=
    -A_{Y_s^\epsilon}\varphi_s^{\epsilon,t,\psi}(x)\,
    dr^\epsilon(s)
    -
    f(x,Y_s^\epsilon)\varphi_s^{\epsilon,t,\psi}(x)\,
    dW_s^\epsilon        \\
    &\quad
    +
    \frac12
    f(x,Y_s^\epsilon)^2
    \varphi_s^{\epsilon,t,\psi}(x)\,ds .
\end{aligned}
\]
Passing to the limit along the smooth approximations gives the claimed rough
backward equation and the locally uniform convergence of
\(\varphi^{\epsilon,t,\psi}\) to \(\varphi^{t,\psi}\).
\end{proof}

\begin{proposition}[Smooth duality]\label{prop:smooth_duality}
Assume the Conditions of Theorem \ref{lem:backward_dual_equation} hold. Let \(\hat{J}^\epsilon\) be a sequence in as in Definition \ref{approximation_class_rpde}, and let
\(\varkappa^\epsilon\) be a finite measure-valued solution of the smooth weak
Zakai equation
\[
\begin{aligned}
    \varkappa_t^\epsilon(\phi)
    &=
    \varkappa_0^\epsilon(\phi)
    +
    \int_0^t
        \varkappa_s^\epsilon
        \bigl(A_{Y_s^\epsilon}\phi\bigr)\,dc^\epsilon_s        \\
    &\quad
    +
    \int_0^t
        \varkappa_s^\epsilon
        \bigl(f(\cdot,Y_s^\epsilon)\phi\bigr)\,dW_s^\epsilon
    -
    \frac12
    \int_0^t
        \varkappa_s^\epsilon
        \bigl(f(\cdot,Y_s^\epsilon)^2\phi\bigr)\,ds .
\end{aligned}
\]
Let \(\varphi^{\epsilon,t,\psi}\) be the smooth backward solution with
terminal condition
\(
    \varphi_t^{\epsilon,t,\psi}=\psi.
\)
Then
\(
    \varkappa_t^\epsilon(\psi)
    =
    \varkappa_0^\epsilon(\varphi_0^{\epsilon,t,\psi}).
\)
\end{proposition}
\begin{proof}
Since all paths are smooth, we may test the forward weak equation against the
time-dependent function \(\varphi_s^{\epsilon,t,\psi}\). Using the
smooth weak Zakai equation and the backward equation, we compute
\[
\begin{aligned}
    d\varkappa_s^\epsilon
    \bigl(\varphi_s^{\epsilon,t,\psi}\bigr)
    &=
    \varkappa_s^\epsilon
    \bigl(
        A_{Y_s^\epsilon}\varphi_s^{\epsilon,t,\psi}
    \bigr)\,dc^\epsilon_s        \\
    &\quad
    +
    \varkappa_s^\epsilon
    \bigl(
        f(\cdot,Y_s^\epsilon)
        \varphi_s^{\epsilon,t,\psi}
    \bigr)\,dW_s^\epsilon        \\
    &\quad
    -
    \frac12
    \varkappa_s^\epsilon
    \bigl(
        f(\cdot,Y_s^\epsilon)^2
        \varphi_s^{\epsilon,t,\psi}
    \bigr)\,ds        \\
    &\quad
    +
    \varkappa_s^\epsilon
    \bigl(
        d\varphi_s^{\epsilon,t,\psi}
    \bigr).
\end{aligned}
\]
But
\[
\begin{aligned}
    d\varphi_s^{\epsilon,t,\psi}
    &=
    -A_{Y_s^\epsilon}
        \varphi_s^{\epsilon,t,\psi}\,dc^\epsilon_s
    -
    f(\cdot,Y_s^\epsilon)
        \varphi_s^{\epsilon,t,\psi}\,dW_s^\epsilon        \\
    &\quad
    +
    \frac12
    f(\cdot,Y_s^\epsilon)^2
        \varphi_s^{\epsilon,t,\psi}\,ds .
\end{aligned}
\]
Substituting this expression into the previous equation, all terms cancel and
we get
\(
    d\varkappa_s^\epsilon
    \bigl(\varphi_s^{\epsilon,t,\psi}\bigr)=0.
\)
Therefore
\(
    \varkappa_t^\epsilon
    \bigl(\varphi_t^{\epsilon,t,\psi}\bigr)
    =
    \varkappa_0^\epsilon
    \bigl(\varphi_0^{\epsilon,t,\psi}\bigr).
\)
Since \(\varphi_t^{\epsilon,t,\psi}=\psi\), this gives
\(
    \varkappa_t^\epsilon(\psi)
    =
    \varkappa_0^\epsilon(\varphi_0^{\epsilon,t,\psi}).
\)
\end{proof}

\begin{proposition}[Rough duality]
\label{prop:rough_duality_by_approximation}
In the context of the previous proposition let \(\varkappa\) be a finite measure-valued weak solution of the rough Zakai
equation and let \(\varkappa^\epsilon\) be a sequence of smooth approximations of $\varkappa$. Assume that the convergence is compatible with the backward test
functions, namely, for every \(t\in[0,T]\) and
\(\psi\in C_c^\infty(\mathbb R)\),
\(
    \varkappa_0^\epsilon
    \bigl(\varphi_0^{\epsilon,t,\psi}\bigr)
    \to
    \varkappa_0
    \bigl(\varphi_0^{t,\psi}\bigr),
\)
and
\(
    \varkappa_t^\epsilon(\psi)
    \to
    \varkappa_t(\psi).
\)
Then
\[
    \varkappa_t(\psi)
    =
    \varkappa_0(\varphi_0^{t,\psi}),
    \qquad
    \psi\in C_c^\infty(\mathbb R),
    \quad t\in[0,T].
\]
\end{proposition}

\begin{proof}
For every \(\epsilon>0\), Proposition \ref{prop:smooth_duality} gives
\(
    \varkappa_t^\epsilon(\psi)
    =
    \varkappa_0^\epsilon
    \bigl(\varphi_0^{\epsilon,t,\psi}\bigr).
\)
Passing to the limit \(\epsilon\to0\), using the assumed convergence on both
sides, yields
\(
    \varkappa_t(\psi)
    =
    \varkappa_0(\varphi_0^{t,\psi}).
\)
\end{proof}

\begin{corollary}[Weak uniqueness]\label{cor:weak_uniqueness_rough_zakai}
Let \(\varkappa\) and \(\nu\) be two weak solutions of the rough Zakai equation
belonging to the approximation class of Proposition
\ref{prop:rough_duality_by_approximation}. Assume that
\(
    \varkappa_0=\nu_0.
\)
Then
\[
    \varkappa_t=\nu_t,
    \qquad t\in[0,T].
\]
\end{corollary}

\begin{proof}
Fix \(t\in[0,T]\) and \(\psi\in C_c^\infty(\mathbb R)\). By Proposition
\ref{prop:rough_duality_by_approximation},
\[
    \varkappa_t(\psi)
    =
    \varkappa_0(\varphi_0^{t,\psi}),
    \qquad
    \nu_t(\psi)
    =
    \nu_0(\varphi_0^{t,\psi}).
\]
Since \(\varkappa_0=\nu_0\), we get
\(
    \varkappa_t(\psi)=\nu_t(\psi).
\)
Because \(\psi\in C_c^\infty(\mathbb R)\) was arbitrary, it follows that
\(
    \varkappa_t=\nu_t.
\)
Since \(t\in[0,T]\) was arbitrary, the proof is complete.
\end{proof}
\begin{proposition}\label{prop:analytic_candidate_duality}
Assume the hypotheses of Theorem
\ref{thm:lifted_rough_mz_wellposedness}.
Define
\(
    \varkappa_t^u(dx):=\tilde{\varrho}(t,x)\,dx.
\)
Then \(\varkappa^u\)
satisfies the duality formula
\[
    \varkappa_t^{\tilde{\varrho}}(\psi)
    =
    \varkappa_0^{\tilde{\varrho}}(\varphi_0^{t,\psi}),
    \qquad
    \psi\in C_c^\infty(\mathbb R).
\]
\end{proposition}
\begin{proof}
Let \(\tilde{\varrho}^\epsilon\) be the smooth analytic solution associated with the
smooth approximation \(\hat J^\epsilon\), and set
\[
    \varkappa_t^{\tilde{\varrho},\epsilon}(dx):=\tilde{\varrho}^\epsilon(t,x)\,dx.
\]
By Theorem \ref{thm:lifted_rough_mz_wellposedness},
\[
    \tilde{\varrho}^\epsilon
    \to
    \tilde{\varrho}
    \qquad
    \text{locally uniformly on }[0,T]\times\mathbb R.
\]
For every \(\epsilon>0\), the smooth analytic solution satisfies the smooth
weak Zakai equation. Hence Proposition \ref{prop:smooth_duality} gives
\(
    \varkappa_t^{\tilde{\varrho},\epsilon}(\psi)
    =
    \varkappa_0^{\tilde{\varrho},\epsilon}
    \bigl(\varphi_0^{\epsilon,t,\psi}\bigr).
\)
Since
\(
    \varkappa_0^{\tilde{\varrho},\epsilon}(dx)=\tilde{\varrho}_0(x)\,dx
\)
and since
\(
    \varphi_0^{\epsilon,t,\psi}
    \to
    \varphi_0^{t,\psi}
\)
locally uniformly, with uniform boundedness, we obtain
\(
    \varkappa_0^{u,\epsilon}
    \bigl(\varphi_0^{\epsilon,t,\psi}\bigr)
    \to
    \varkappa_0^u
    \bigl(\varphi_0^{t,\psi}\bigr).
\)
Similarly, the locally uniform convergence of \(u^\epsilon\) to \(u\) gives
\(
    \varkappa_t^{u,\epsilon}(\psi)
    \to
    \varkappa_t^u(\psi).
\)
Passing to the limit in the smooth duality identity yields
\(
    \varkappa_t^u(\psi)
    =
    \varkappa_0^u(\varphi_0^{t,\psi}).
\)
\end{proof}

We now apply the backward uniqueness argument to the rough weak equation
satisfied by the unnormalised filter. The key point is that, when the test
function is evolved backwards according to the adjoint equation, its pairing
with the forward filter is conserved in time. This gives the following duality
identity.

\begin{proposition}\label{prop:filter_duality}
Assume Conditions \ref{condition_covariance_rho_var}, \ref{cameron_martin_condition}, \ref{condition_filtration}, \ref{condition_7} and \ref{cond:vector_field_reg_zakai} hold. 
Let $\tilde \pi$ be the unnormalized filter, then 
\[
    \tilde\pi_t(\psi)
    =
    \tilde\pi_0(\varphi_0^{t,\psi}),
    \qquad
    \psi\in C_c^\infty(\mathbb R).
\]
\end{proposition}
\begin{proof}

The first step in the proof consists in showing that for any smooth approximation $\hat Y^\epsilon$, the corresponding unnormalized filter satisfies equation \eqref{weak_smooth_zakai}, then by showing that \(
    \tilde\xi_0^\epsilon
    \bigl(\varphi_0^{\epsilon,t,\psi}\bigr)
    \to
    \tilde\xi_0
    \bigl(\varphi_0^{t,\psi}\bigr).
\)
allows to conclude by Proposition \ref{prop:smooth_duality}.
Fix \(\epsilon>0\) and let \(\phi\in C_c^\infty(\mathbb R)\). We write
\begin{equation}\label{representation_tilde_pi_eps}
    \tilde \pi_t^\epsilon(\phi)
    :=
    \mathbb E^{\overline{\mathbb Q}}
    \left[
        \phi(X_t^\epsilon)\exp(\Xi_t^\epsilon)
    \right].
\end{equation}
By the robust representation of the approximating unnormalised filter, this
is precisely the value of the approximating measure \(\tilde \pi_t^\epsilon\)
on the test function \(\phi\).

Since the approximating driver is smooth, the dynamics are classical and we
may apply the ordinary chain rule to
\[
    t\to
    \phi(X_t^\epsilon)\exp(\Xi_t^\epsilon).
\]
Thus, for every \(t\in[0,T]\),
\[
\begin{aligned}
    \tilde \pi_t^\epsilon(\phi)-\tilde \pi_0^\epsilon(\phi)
    &=
    {\mathbb E ^{\overline{\mathbb{Q}}}}
    \left[ \int_0^t \partial_x\phi(X_s^\epsilon)  \sigma(X_s^\epsilon,Y_s^\epsilon)
            \exp(\Xi_s^\epsilon)
            \,dB_s^\epsilon \right]\\
    &\quad
    +
    \int_0^t
        {\mathbb E^{\overline{\mathbb{Q}}}}
        \left[
            f(X_s^\epsilon,Y_s^\epsilon)
            \phi(X_s^\epsilon)
            \exp(\Xi_s^\epsilon)
        \right]
        d\hat Y_s^\epsilon \\
    &\quad
    -
    \frac12
    \int_0^t
        {\mathbb E^{\overline{\mathbb{Q}}}}
        \left[
            f(X_s^\epsilon,Y_s^\epsilon)^2
            \phi(X_s^\epsilon)
            \exp(\Xi_s^\epsilon)
        \right]ds .
\end{aligned}
\]
Equivalently,
\[
\begin{aligned}
    \tilde \pi_t^\epsilon(\phi)-\tilde \pi_0^\epsilon(\phi)
    &=
    {\mathbb E^{\overline{\mathbb{Q}}}}
    \left[
        \int_0^t
            \partial_x\phi(X_s^\epsilon)
            \sigma(X_s^\epsilon,Y_s^\epsilon)
            \exp(\Xi_s^\epsilon)
            \,dB_s^\epsilon
    \right] \\
    &\quad
    +
    \int_0^t
        \tilde \pi_s^\epsilon
        \bigl(f(\cdot,Y_s^\epsilon)\phi\bigr)
        d\hat Y_s^\epsilon
    -
    \frac12
    \int_0^t
        \tilde \pi _s^\epsilon
        \bigl(f(\cdot,Y_s^\epsilon)^2\phi\bigr)ds .
\end{aligned}
\]

It remains to identify the expectation of the integral with respect to
\(B^\epsilon\). Set
\[
    S_y\phi(x):=
    \sigma(x,y)\partial_x\phi(x).
    \]
With this notation,
\(
    \partial_x\phi(X_s^\epsilon)
    \sigma(X_s^\epsilon,Y_s^\epsilon)
    =
    S_{Y_s^\epsilon}\phi(X_s^\epsilon).
\)
By the Gaussian integration-by-parts formula for the smooth Gaussian
approximation \(B^\epsilon\) (see Section 1.3 in \cite{nualart2009malliavin}), together with the differentiability of the
solution map we obtain
\[
\begin{aligned}
    {\mathbb E^{\overline{\mathbb{Q}}}}
    \left[
        \int_0^t
            S_{Y_s^\epsilon}\phi(X_s^\epsilon)
            \exp(\Xi_s^\epsilon)
            \,dB_s^\epsilon
    \right]
    &=
    \int_0^t
        {\mathbb E^{\overline{\mathbb{Q}}}}
        \left[
            A_{Y_s^\epsilon}\phi(X_s^\epsilon)
            \exp(\Xi_s^\epsilon)
        \right]
        dc^\epsilon_s.
\end{aligned}
\]
Indeed, the Frechet derivative of
\(S_{Y_s^\epsilon}\phi(X_s^\epsilon)\) in the
\(B^\epsilon\)-direction is transported by the Jacobian.
Since the signal vector field is one-dimensional, it commutes with itself,
and the transported derivative reduces on the diagonal to
\[
    S_{Y_s^\epsilon}
    \bigl(S_{Y_s^\epsilon}\phi\bigr)
    (X_s^\epsilon)
    =
    A_{Y_s^\epsilon}\phi(X_s^\epsilon).
\]
Consequently,
\[
\begin{aligned}
    {\mathbb E^{\overline{\mathbb{Q}}}}
    \left[
        \int_0^t
            \partial_x\phi(X_s^\epsilon)
            \sigma(X_s^\epsilon,Y_s^\epsilon)
            \exp(\Xi_s^\epsilon)
            \,dB_s^\epsilon
    \right]
    &=
    \int_0^t
        \tilde \pi_s^\epsilon
        \bigl(A_{Y_s^\epsilon}\phi\bigr)
        dc^\epsilon_s.
\end{aligned}
\]

Substituting this identity into the previous expression for
\(\tilde \pi_t^\epsilon(\phi)-\tilde \pi_0^\epsilon(\phi)\), we obtain
\[
\begin{aligned}
    \tilde \pi_t^\epsilon(\phi)
    &=
    \tilde \pi_0^\epsilon(\phi)
    +
    \int_0^t
        \tilde \pi_s^\epsilon
        \bigl(A_{Y_s^\epsilon}\phi\bigr)
        dc^\epsilon_s \\
    &\quad
    +
    \int_0^t
        \tilde \pi_s^\epsilon
        \bigl(f(\cdot,Y_s^\epsilon)\phi\bigr)
        d\hat Y_s^\epsilon
    -
    \frac12
    \int_0^t
        \tilde \pi_s^\epsilon
        \bigl(f(\cdot,Y_s^\epsilon)^2\phi\bigr)ds .
\end{aligned}
\]
This proves that each smooth approximating measure
\(\tilde \pi^\epsilon\) satisfies the classical weak Zakai equation. 
Using the representation \eqref{representation_tilde_pi_eps} and passing to 
the limit in the sense of Definition
\ref{def:weak_measure_rough_zakai_approx} gives the desired weak solution
of the rough Zakai equation by using an argument analogous to the one employed in Section \ref{section_robustness} to show robustness of the filter.
For the second part of the proof, we begin by recalling that from Proposition \ref{lem:backward_dual_equation},  \(\varphi^{\epsilon,t,\psi}_0\to\varphi^{t,\psi}_0\) locally uniformly on \(\mathbb R\) and 
\[
    C:=
    \sup_{\epsilon>0}
    \|\varphi^{\epsilon,t,\psi}_0\|_\infty
    +
    \|\varphi^{t,\psi}_0\|_\infty < \infty.
\]
We have
\(
    \pi_0
    \left(
        \left|
            \varphi^{\epsilon,t,\psi}_0
            -
            \varphi^{t,\psi}_0
        \right|
    \right)
    \to 0 .
\)
Indeed, for every \(M>0\),
\[
\begin{aligned}
    \pi_0
    \left(
        \left|
            \varphi^{\epsilon,t,\psi}_0
            -
            \varphi^{t,\psi}_0
        \right|
    \right)
    &\leq
    \pi_0(\mathbb{1}_{\mathbb{R}})
    \sup_{|x|\leq R}
    \left|
        \varphi^{\epsilon,t,\psi}_0(x)
        -
        \varphi^{t,\psi}_0(x)
    \right| \\
    &\quad
    +
    C\,\pi_0\bigl(\mathbb{1}_{\{|x|>M\}}\bigr).
\end{aligned}
\]
The first term vanishes as \(\epsilon\to0\), for fixed \(M\), by local uniform convergence. The second term can be made arbitrarily small by choosing \(M\) sufficiently large, since \(\pi_0\) is a finite measure. Therefore
\(
    \pi_0(\varphi^{\epsilon,t,\psi}_0)
    \to
    \pi_0(\varphi^{t,\psi}_0).
\)
As we mentioned above we can now invoke Proposition \ref{prop:smooth_duality}, which, for each \(\epsilon>0\), gives
\(
    \tilde\pi_t^\epsilon(\psi)
    =
    \tilde\pi_0^\epsilon
    \bigl(\varphi_0^{\epsilon,t,\psi}\bigr).
\)
Passing to the limit \(\epsilon\to0\), using the assumed convergence on both
sides, gives
\(
    \tilde\pi_t(\psi)
    =
    \tilde\pi_0(\varphi_0^{t,\psi}).
\)
\end{proof}

\begin{remark}
This is the point at which the one-dimensional nature of the signal is used. In
dimension one, the integration-by-parts term closes as
\[
    S_y(S_y\phi)(x)=A_y\phi(x),
\]
and therefore produces exactly the second-order operator appearing in the weak
equation.

In higher dimension, the same computation generates additional commutator terms
involving iterated integrals of the driving signal. These terms do not collapse
to the local operator \(A_y\) driven only by the diagonal covariance clock.
Thus the integration-by-parts argument no longer yields the weak equation in
the form considered here. Extending the result beyond the one-dimensional
signal case would therefore require a different framework, capable of retaining
these additional rough second-order contributions.
\end{remark}
We can now complete the identification of the analytic solution constructed in
Theorem \ref{thm:lifted_rough_mz_wellposedness}. That theorem provides a
continuous function \(\tilde{\varrho}\), but does not by itself identify \(\tilde{\varrho}(t,\cdot)\) with the conditional
density of the filter. The backward duality relation proved above provides this
identification: it shows that the measure \(\tilde{\varrho}(t,x)\,dx\) and the
unnormalised filter have the same action on test functions. Hence the analytic
RPDE solution coincides with the density of the filtering problem.

\begin{theorem}\label{thm:continuous_density_by_duality}
Assume that Conditions \ref{cameron_martin_condition} and the hypotheses of Theorem
\ref{thm:lifted_rough_mz_wellposedness} hold. Let
\(
    \tilde{\varrho}(t,x)
\)
be the rough viscosity solution of equation \eqref{rough_mortensen_zakai}. Assume that
\[
    \tilde\pi_0(dx)=\tilde{\varrho}_0(x)\,dx.
\]
Then
\[
    \tilde\pi_t(dx)=\tilde{\varrho}(t,x)\,dx,
    \qquad t\in[0,T].
\]
In particular, the unnormalised filter admits a continuous density $\tilde{\varrho} = u$.
\end{theorem}
\begin{proof}
Define
\[
    \varkappa_t^{\tilde{\varrho}}(dx):=\tilde{\varrho}(t,x)\,dx.
\]
By Proposition \ref{prop:analytic_candidate_duality},
\(
    \varkappa_t^{\tilde{\varrho}}(\psi)
    =
    \varkappa_0^{\tilde{\varrho}}(\varphi_0^{t,\psi}),
    \) with \(
    \psi\in C_c^\infty(\mathbb R).
\)
Moreover, by Proposition \ref{prop:filter_duality},
\(
    \tilde\pi_t(\psi)
    =
    \tilde\pi_0(\varphi_0^{t,\psi}).
\)
Since
\(
    \tilde\pi_0(dx)=\tilde{\varrho}_0(x)\,dx=\varkappa_0^{\tilde{\varrho}}(dx),
\)
we have
\(
    \tilde\pi_0(\varphi_0^{t,\psi})
    =
    \varkappa_0^u(\varphi_0^{t,\psi}).
\)
Therefore
\(
    \tilde\pi_t(\psi)
    =
    \varkappa_t^{\tilde{\varrho}}(\psi).
\)
Since \(C_c^\infty(\mathbb R)\) separates finite measures, it follows that
\(
    \tilde\pi_t=\varkappa_t^{\tilde{\varrho}}.
\)
Hence
\(
    \tilde\pi_t(dx)=\tilde{\varrho}(t,x)\,dx.
\)

Finally, defining $P$ as in Theorem \ref{thm:lifted_rough_mz_wellposedness}
\(
    \tilde{\varrho}(t,x)=P(t,x,Y_t),
\)
where \(P\) is continuous and \(Y\) is continuous. Hence
\(
    (t,x)\to \tilde{\varrho}(t,x)
\)
is continuous on \([0,T]\times\mathbb R\). Therefore the unnormalised filter
admits a continuous density.
\end{proof}

\section{Analysis of the density of the conditional distribution} \label{section_density}
In this section, we continue the study of the conditional law of the signal $X_t$ given $\mathcal{Y}_t$  by establishing its smoothness with respect to the Lebesgue measure. Our approach leverages partial Malliavin calculus, more precisely we use the framework introduced in \cite{nualart1989partial}, which extends the Bouleau and Hirsch smoothness criterion of \cite{bouleau1986proprietes} from the classical Malliavin calculus setting to the partial Malliavin framework.  

We begin with an overview of the classical and  partial Malliavin calculus framework, introducing the necessary tools used in the main part of this section.  

\subsection{Preliminaries on Malliavin calculus}

Let
\(
    I=\{I(h):h\in\mathfrak H_1^d\}
\)
be the real-valued isonormal Gaussian process over \(\mathfrak H_1^d\). Let \(\mathbb S( \mathbb R^m )\) denote the set of \(\mathbb R^m\)-valued smooth cylindrical random variables of the form
\[
\Phi
=
\sum_{j=1}^{N}
\phi_j\bigl(I(h_{j,1}),\ldots,I(h_{j,n_j})\bigr)v_j,
\]
where \(\phi_j \in C^\infty(\mathbb R^{n_j},\mathbb R)\) and \(\phi_j\), together with all its derivatives, has at most polynomial growth, \(h_{j,\ell}\in \mathfrak H_1^d\), and \(v_j\in \mathbb R^m\).

For \(k\geq 1\), the \(k\)-th Malliavin derivative of \(\Phi\) is the random variable with values in
\(
\mathbb R^m\otimes (\mathfrak H_1^d)^{\otimes k}
\)
defined by
\[
\mathcal D^k\Phi
:=
\sum_{j=1}^{N}
\sum_{i_1,\ldots,i_k=1}^{n_j}
\frac{\partial^k \phi_j}{\partial x_{i_1}\cdots \partial x_{i_k}}
\bigl(I(h_{j,1}),\ldots,I(h_{j,n_j})\bigr)
\,v_j\otimes h_{j,i_1}\otimes\cdots\otimes h_{j,i_k}.
\]
The operator \(\mathcal D^k\) is closable from \(\mathbb S(\mathbb R^m)\) to
\(
L^r\bigl(\Omega,\mathbb R^m\otimes(\mathfrak H_1^d)^{\otimes k}\bigr)
\)
for every \(r\geq 1\).

We denote by \(\mathbb D^{k,r}(\mathbb R^m)\) the completion of \(\mathbb S(\mathbb R^m)\) with respect to the norm
\[
\|\Phi\|_{\mathbb D^{k,r}}^r
:=
\|\Phi\|_{L^r(\Omega,\mathbb R^m)}^r
+
\sum_{i=1}^{k}
\|\mathcal D^i\Phi\|_{L^r(\Omega,\mathbb R^m \otimes(\mathfrak H_1^d)^{\otimes i})}^r .
\]
We call these spaces Malliavin-Sobolev spaces.
We also define
\[
\mathbb D^\infty(\mathbb R^m)
:=
\bigcap_{k\geq 1}\bigcap_{r\geq 1}\mathbb D^{k,r}(\mathbb R^m).
\]

Let \(\operatorname{Dom}(\delta^k)\) be the set of all
\(
u\in L^2\bigl(\Omega,(\mathfrak H_1^d)^{\otimes k}\bigr)
\)
such that there exists a constant \(C>0\) satisfying
\[
\left|
\mathbb E^\mu\left[
\left\langle \mathcal D^k\Phi,u\right\rangle_{(\mathfrak H_1^d)^{\otimes k}}
\right]
\right|
\leq
C\|\Phi\|_{L^2(\Omega)}
\]
for every \(\Phi\in \mathbb S(\mathbb R)\). By the Riesz representation theorem, for every \(u\in\operatorname{Dom}(\delta^k)\) there exists a unique element \(\delta^k(u)\in L^2(\Omega)\), called the divergence of \(u\), such that
\begin{equation}\label{classic_integration_by_parts}
\mathbb E^\mu\left[\delta^k(u)\Phi\right]
=
\mathbb E^\mu\left[
\left\langle \mathcal D^k\Phi,u\right\rangle_{(\mathfrak H_1^d)^{\otimes k}}
\right]
\end{equation}
for every \(\Phi\in\mathbb S(\mathbb R)\). When \(k=1\), we write \(\mathcal D\Phi\) and \(\delta(u)\) instead of \(\mathcal D^1\Phi\) and \(\delta^1(u)\). In this case, for \(h\in\mathfrak H_1^d\), we also write
\[
\mathcal D_h\Phi:=\langle \mathcal D\Phi,h\rangle_{\mathfrak H_1^d}.
\]

\subsection{Partial Malliavin calculus}

Let
\(
\mathcal H=\{\mathscr{L}(\omega):\omega\in\Omega\}
\)
be a collection of closed subspaces of \(\mathfrak H_1^d\). We assume that the associated orthogonal projection is measurable, in the sense that, for every \(h\in\mathfrak H_1^d\), the map
\(
\omega\to \Pi_{\mathscr{L}(\omega)}h
\)
is measurable as an \(\mathfrak H_1^d\)-valued random variable. Equivalently, if \((e_i)_{i\geq 1}\) is an orthonormal basis of \(\mathfrak H_1^d\), then for every \(\mathfrak H_1^d\)-valued random variable \(\Phi\) the random variable
\[
\Pi_{\mathcal H}\Phi(\omega)
:=
\Pi_{\mathscr{L}(\omega)}\Phi(\omega)
=
\sum_{i\geq 1}
\langle \Phi(\omega),e_i\rangle_{\mathfrak H_1^d}
\Pi_{\mathscr{L}(\omega)}e_i
\]
is measurable.

For \(\Phi \in\mathbb D^{k,r}(\mathbb R)\), we define the \(k\)-th partial Malliavin derivative of \(F\) by
\[
\mathcal D_{\mathcal H}^k\Phi
:=
(\Pi_{\mathcal H})^{\otimes k}D^k\Phi.
\]
That is,
\(
\mathcal D_{\mathcal H}^k\Phi(\omega)
=
(\Pi_{\mathscr{L}(\omega)})^{\otimes k}\mathcal D^k\Phi(\omega).
\)

The corresponding \(k\)-th partial divergence operator is defined by
\[
\operatorname{Dom}(\delta_{\mathcal H}^k)
:=
\left\{
u\in L^2\bigl(\Omega,(\mathfrak H_1^d)^{\otimes k}\bigr)
:
(\Pi_{\mathcal H})^{\otimes k}u
\in \operatorname{Dom}(\delta^k)
\right\},
\]
and
\(
\delta_{\mathcal H}^k(u)
:=
\delta^k\bigl((\Pi_{\mathcal H})^{\otimes k}u\bigr).
\)
Hence, for every \(\Phi \in\mathcal S(\mathbb R)\),
\begin{align*}
\mathbb E\left[\delta_{\mathcal H}^k(u)\Phi\right]
&=
\mathbb E\left[
\delta^k\bigl((\Pi_{\mathcal H})^{\otimes k}u\bigr)\Phi
\right] \\
&=
\mathbb E\left[
\left\langle
\mathcal D^k\Phi,
(\Pi_{\mathcal H})^{\otimes k}u
\right\rangle_{(\mathfrak H_1^d)^{\otimes k}}
\right] \\
&=
\mathbb E\left[
\left\langle
\mathcal D_{\mathcal H}^k\Phi,
u
\right\rangle_{(\mathfrak H_1^d)^{\otimes k}}
\right].
\end{align*}

The next definition introduces the condition which allows one to associate such a family of projections with a sub-\(\sigma\)-field.

\begin{definition}
A sub-\(\sigma\)-field \(\mathcal G\subseteq\mathcal F\) is said to be countably smoothly generated if there exists a sequence of random variables \((G_i)_{i\geq 1}\), with \(G_i\in\mathbb D^{2,1}\), such that
\(
\mathcal G=\sigma(G_i:i\geq 1).
\)
\end{definition}

Given a countably smoothly generated sub-\(\sigma\)-field \(\mathcal G=\sigma(G_i:i\geq 1)\), we associate to it the family
\[
\mathscr{L}(\omega)
:=
\left(
\overline{\operatorname{span}\{\mathcal DG_i(\omega):i\geq 1\}}
\right)^\perp
\subseteq \mathfrak H_1^d.
\]
Equivalently, for \(h\in\mathfrak H_1^d\), the projection \(\Pi_{\mathscr{L}(\omega)}h\) is given by the limit
\[
\Pi_{\mathscr{L}(\omega)}h
=
\lim_{n\to\infty}
\left(
h-
\Pi_{\operatorname{span}\{\mathcal DG_i(\omega):1\leq i\leq n\}}h
\right),
\]
whenever the generating sequence \((G_i)_{i\geq 1}\) is fixed.

\begin{example}\label{c_s_gen_example}
Let \(B\) be a centered \(d\)-dimensional  Volterra Gaussian process with i.i.d components and let
\[
\mathcal F_t^B:=\sigma(B_s:0\leq s\leq t).
\]
Then for every $t$,\(\mathcal F_t^B\) is countably smoothly generated. Indeed, it is generated by the countable family
\[
\left\{
B^i(\mathbbm 1_{[0,q]})
:
i\in[1:d],
\ q\in[0,t]\cap\mathbb Q
\right\}.
\]
Moreover,
\(
\mathcal D B^i(\mathbbm 1_{[0,q]})
=
\mathbbm 1_{[0,q]}^i.
\)
Therefore the associated closed subspace is
\[
\mathfrak H_t^d
:=
\overline{
\operatorname{span}
\left\{
\mathbbm 1_{[0,q]}^i
:
i \in [1:d],
\ q\in[0,t]\cap\mathbb Q
\right\}
},
\]
and the partial Malliavin subspace associated with \(\mathcal F_t^B\) is
\(
(\mathfrak H_t^d)^\perp.
\)
In particular, the corresponding projection is the deterministic orthogonal projection onto \((\mathfrak H_t^d)^\perp\).
\end{example}

\subsection{Existence of smooth density for the unnormalized filter}

Within the framework described above, we prove that, under suitable non-degeneracy assumptions, the conditional distribution of \(X_t\) given \(\mathcal Y\),
\[
\pi_t(\phi)
=
\mathbb E^{\mathbb P}\left[\phi(X_t) \mid \mathcal Y\right],
\]
admits a smooth density with respect to the Lebesgue measure. More precisely, we will show that there exists a version
\[
 \varrho:\Omega\times (0, T]\times\mathbb R^{d_X}\to\mathbb R_+
\]
such that
\[
\pi_t(\phi)
=
\int_{\mathbb R^{d_X}}\phi(x) \varrho(t, x)\,dx,
\qquad
\phi\in C_b(\mathbb R^{d_X}),
\]
and, for every \(\omega\in\Omega, t \in (0, T]\), the map \(x\to \varrho(\omega,t,x)\) belongs to \(C^\infty(\mathbb R^{d_X})\).
We note that, throughout this section, $X$, $Y$, and $B$ are no longer assumed to be one-dimensional.

The proof combines the estimates of \cite{cass2015smoothness} with the partial Bouleau-Hirsch criterion of \cite[Theorem 5.7]{nualart1989partial}. More precisely, since Example~\ref{c_s_gen_example} shows that \(\mathcal Y\) is countably smoothly generated, it remains to verify the following conditions:
\begin{enumerate}
    \item $Z_t \in \mathbb{D}^\infty(\mathbb{R}^{d_Z})$ and the Radon--Nikodym process satisfies
    \(
    \Lambda_t\in\mathbb D^\infty(\mathbb R)
    \) for every $t\in[0,T]$.
    \item The orthogonal projection \(\Pi_{\mathcal H}\) maps
\(\mathbb{D}^{\infty}(\mathfrak{H}_1^{d_J})\) into itself, where $d_J = d_I$.
    \item The partial Malliavin matrix of \(X_t\),
    \[
    \Upsilon_{\mathcal H,t}^{ij}
    :=
    \left\langle
    \mathcal D_{\mathcal H}X_t^i,
    \mathcal D_{\mathcal H}X_t^j
    \right\rangle_{\mathfrak H_1^{d_I}},
    \qquad i,j \in [1:d_X],
    \]
    is almost surely invertible and satisfies
    \[
    (\Upsilon_{\mathcal H,t}^{-1})^{ij}
    \in
    \bigcap_{r\geq 1}L^r(\Omega),
    \qquad i,j \in [1:d_X].
    \]
\end{enumerate}
Once these properties are established, the integration-by-parts formula for the partial Malliavin derivative can be applied with \(F=X_t\) and \(L=\Lambda_t\). The criterion of \cite{nualart1989partial} then yields a \(C^\infty\)-version of the density of the unnormalised conditional distribution.

The additional conditions that we will use to verify these conditions are:
\begin{condition}[Partial H\"ormander condition]\label{condition_3}
Assume that
\(
    \sigma_1,\ldots,\sigma_{d_B}, f
    \in C_b^\infty
    \bigl(
        \mathbb R^{d_X}\times \mathbb R^{d_Y},
        \mathbb R^{d_X}
    \bigr).
\)
Let
\(
    \mathcal S_0
    :=
    \{\sigma_1,\ldots,\sigma_{d_B}\}.
\)
For \(n\geq 0\), define recursively
\[
\begin{aligned}
    \mathcal S_{n+1}
    :=
    \mathcal S_n
    &\cup
    \left\{
        [\sigma_i,S]_x
        :
        S\in\mathcal S_n,\ i=1,\ldots,d_B
    \right\}  \\
    &\cup
    \left\{
        \partial_{y_j}S
        :
        S\in\mathcal S_n,\ j=1,\ldots,d_Y
    \right\}.
\end{aligned}
\]
Here, for smooth vector fields
\(S,\tilde S :\mathbb R^{d_X}\times\mathbb R^{d_Y}\to\mathbb R^{d_X}\), we write
\[
    [S,\tilde S]_x
    :=
    \partial _x\tilde S\,S-\partial_xS\,\tilde S.
\]
We assume that there exists a deterministic \(N\in\mathbb N\) such that
\[
    \operatorname{span}
    \left\{
        S(x_0,0)
        :
        S\in\bigcup_{n=0}^{N}\mathcal S_n
    \right\}
    =
    \mathbb R^{d_X},
    \qquad \mathbb P\text{-a.s.}
\]
\end{condition}

\begin{condition}\label{condition_4}
We assume that one of the following two conditions holds.

\smallskip

\noindent\emph{(A) Local non-determinism and non-negative conditional covariance.}
For each \(i=1,\ldots,d_B\), let
\(
    \mathcal F^i_{a,b}
    :=
    \sigma\left(J^i_r:a\leq r\leq b\right).
\)
There exist constants \(\alpha>0\) and \(c>0\) such that, for every \(i=1,\ldots,d_I\),
\[
    \inf_{0\leq s<t\leq T}
    \frac{
        \operatorname{Var}
        \left(
            J^i_{s,t}
            \mid
            \mathcal F^i_{0,s}\vee \mathcal F^i_{t,T}
        \right)
    }{
        |t-s|^\alpha
    }
    \geq c.
\]
Moreover, for every
\(0\leq s<t\leq u<v\leq T\),
\[
    \operatorname{Cov}
    \left(
        J^i_{s,t},
        J^i_{u,v}
        \mid
        \mathcal F^i_{0,s}\vee \mathcal F^i_{t,T}
    \right)
    \geq 0.
\]

\smallskip

\noindent\emph{(B) Interpolated functional covariance coercivity.}
We assume that there exist constants
\[
    \gamma>0,
    \qquad
    \vartheta\in(0,1],
    \qquad
    \eta \geq 0,
    \qquad
    C>0,
\]
with $ \eta/\vartheta +\gamma(1-\vartheta)/\vartheta < 2/p \wedge 1$, 
such that, for every interval \([s,t]\subset[0,T]\) and every \(\varphi\in C^\gamma([s,t])\),
\[
    \|\varphi\|_{\infty;[s,t]}
    \leq
    C |t-s|^{-\eta}
    \bigl\|K^*(\varphi \mathbbm 1_{[s,t]})\bigr\|_{L^2([0,T])}^{\vartheta}
    \|\varphi\|_{C^\gamma([s,t])}^{1-\vartheta}.
\]
\end{condition}

Condition \ref{condition_4} is satisfied by Type I and Type II fractional Brownian motion, see Appendix \ref{app:examples}.

The next goal is to establish the Malliavin smoothness of the likelihood
process \(\Lambda_t\). We do this by first studying its logarithm. Under the
standing assumptions, Theorem~1.2 in \cite{inahama2014malliavin} applies to
the dynamics of \(Z\), and gives \(Z_t\in\mathbb D^\infty\). The remaining
steps are then to verify that the deterministic and It\^o integral terms in
\(\pi_t\) preserve Malliavin smoothness, and that the boundedness of
the observation coefficient provides the exponential moment estimates needed
to prove the result for  \(\Lambda_t=:\exp(\chi_t)\).

\begin{proposition}\label{prop:skorohod_brownian_smooth}
Assume Condition \ref{condition_covariance_rho_var}, \ref{condition_filtration} and \ref{condition_3} hold. Let \(S\in C_b^\infty(\mathbb R^{d_X+d_Y}, \mathbb R)\), let
\(i\in [1:d_Y]\), and let \(t\in[0,T]\). Then
\(
    \delta^{W,i}\bigl(S(Z)\mathbbm 1_{[0,t]}\bigr)
    \in \mathbb D^\infty(\mathbb R).
\)
Here $\delta^{W}$ denotes the divergence operator with respect to the Brownian motion $W$.
\end{proposition}
\begin{proof}

Fix \(k\geq 1\) and \(r>1\), and let
\(
    \Phi\in\mathbb D^{k,r}\bigl(L^2([0,T])\bigr).
\)
By Condition~\ref{condition_filtration}, the operator
\((K^*)^{-1}\colon L^2([0,T])\to\mathfrak H_1\) is well defined and
bounded. Therefore,
\(
    (K^*)^{-1}\Phi
    \in
    \mathbb D^{k,r}(\mathfrak H_1).
\)
In particular, both divergence operators below are well defined, and we can write
\begin{equation}\label{transfer_identity_divergence}
    \delta^{W,i}(\Phi)
    =
    \delta^i\bigl((K^*)^{-1}\Phi\bigr).
\end{equation}
This is enough to apply Proposition~1.5.7 in \cite{nualart2009malliavin} and obtain
\[
\begin{aligned}
    \left\|
        \delta^{W,i}(\Phi)
    \right\|_{\mathbb D^{k-1,r}}
    &=
    \left\|
        \delta^i\bigl((K^*)^{-1}\Phi\bigr)
    \right\|_{\mathbb D^{k-1,r}} \\
    &\leq
    C_{k,r}
    \left\|
        (K^*)^{-1}\Phi
    \right\|_{\mathbb D^{k,r}(\mathfrak H_1)} .
\end{aligned}
\]
Therefore it is enough to prove that
$    (K^*)^{-1}\bigl(S(Z)\mathbbm 1_{[0,t]}\bigr)
    \in
    \mathbb D^\infty(\mathfrak H_1).
$
Let
\(    \mathcal P=\{0=t_0<t_1<\cdots<t_n=t\}
\)
be a dyadic partition of \([0,t]\). Define the left-point approximation
\[
    S_{\mathcal P}(s)
    :=
    \sum_{j=1}^n
    S(Z_{t_{j-1}})
    \mathbbm 1_{(t_{j-1},t_j]}(s),
    \qquad s\in[0,T],
\]
where \(S_{\mathcal P}(s)=0\) for \(s\notin[0,t]\). Set
\[
    \widetilde S_{\mathcal P}
    :=
    (K^*)^{-1}S_{\mathcal P}
    =
    \sum_{j=1}^n
    S(Z_{t_{j-1}})
    (K^*)^{-1}\mathbbm 1_{(t_{j-1},t_j]} .
\]
For every fixed partition \(\mathcal P\), as already noted before, Theorem 1.2 in
\cite{inahama2014malliavin} gives
\[
    Z_{t_{j-1}}\in
    \mathbb D^\infty(\mathbb R^{d_Z}),
    \qquad j=1,\ldots,n.
\]
Since \(S\in C_b^\infty(\mathbb R^{d_X+d_Y},\mathbb R)\), it follows from
the chain rule that
\[
    S(Z_{t_{j-1}})\in \mathbb D^\infty(\mathbb R),
    \qquad j=1,\ldots,n.
\]
Moreover,
\(
    (K^*)^{-1}\mathbbm 1_{(t_{j-1},t_j]}
    \in \mathfrak H_1.
\)
Hence
\(
    \widetilde S_{\mathcal P}
    \in
    \mathbb D^\infty(\mathfrak H_1).
\)
We now prove that \(\widetilde S_{\mathcal P}\) converges in the Malliavin
Sobolev norms. By Theorem 1.2 in \cite{inahama2014malliavin}, the map
\(
    s\to Z_s
\)
is continuous from \([0,T]\) into
\(\mathbb D^\infty(\mathbb R^{d_Z})\). Since \(S\in C_b^\infty\), the map
\(
    s\to S(Z_s)
\)
is continuous from \([0,T]\) into \(\mathbb D^\infty(\mathbb R)\). In
particular, for every \(q>1\) and every \(m\geq 0\),
\[
    \sup_{|u-v|\leq |\mathcal P|}
    \left\|
        S(Z_u)-S(Z_v)
    \right\|_{\mathbb D^{m,q}}
    \to 0
    \qquad
    \text{as }|\mathcal P|\to 0 .
\]
Let
\[
    \operatorname{proj}_{\mathcal P}(s)=t_{j-1}
    \qquad
    \text{whenever }s\in(t_{j-1},t_j].
\]
Then
\(
    S_{\mathcal P}(s)
    =
    S(Z_{\operatorname{proj}_{\mathcal P}(s)})
    \mathbbm 1_{[0,t]}(s).
\)
Choose \(q\geq 2\) with \(q\geq r\). For \(m\in\{0,\ldots,k\}\), by
Minkowski's inequality,
\[
\begin{aligned}
    &\left\|
        \mathcal D^m\left(
            S_{\mathcal P}
            -
            S(Z)\mathbbm 1_{[0,t]}
        \right)
    \right\|_{L^q(\Omega,\mathfrak H_1^{\otimes m}\otimes L^2([0,T]))}
    \\
    &\qquad
    =
    \left\|
        \left(
            \int_0^t
            \left\|
                \mathcal D^m\left(
                    S(Z_{\operatorname{proj}_{\mathcal P}(s)})-S(Z_s)
                \right)
            \right\|_{\mathfrak H_1^{\otimes m}}^2
            ds
        \right)^{1/2}
    \right\|_{L^q(\Omega)}
    \\
    &\qquad
    \leq
    \left(
        \int_0^t
        \left\|
            \mathcal D^m\left(
                S(Z_{\operatorname{proj}_{\mathcal P}(s)})-S(Z_s)
            \right)
        \right\|_{L^q(\Omega,\mathfrak H_1^{\otimes m})}^2
        ds
    \right)^{1/2}
    \\
    &\qquad
    \leq
    t^{1/2}
    \sup_{|u-v|\leq |\mathcal P|}
    \left\|
        S(Z_u)-S(Z_v)
    \right\|_{\mathbb D^{k,q}} .
\end{aligned}
\]
The right-hand side converges to \(0\) as \(|\mathcal P|\to 0\). Therefore
\[
    S_{\mathcal P}
    \to
    S(Z)\mathbbm 1_{[0,t]}
    \qquad
    \text{in }
    \mathbb D^{k,q}(L^2([0,T])) .
\]
Since \(q\geq r\), this also implies convergence in
\(\mathbb D^{k,r}(L^2([0,T]))\).
Using that \(K^*\colon \mathfrak H_1\to \mathfrak H_2\) is an isometry, we obtain
\[
    \left\|
        (K^*)^{-1}h
    \right\|_{\mathfrak H_1}
    =
    \|h\|_{L^2}
\]
for every \(h\) in the range of \(K^*\). Since \((K^*)^{-1}\) is deterministic
and linear, it commutes with Malliavin differentiation. Hence
\[
    \widetilde S_{\mathcal P}
    =
    (K^*)^{-1}S_{\mathcal P}
    \to
    (K^*)^{-1}\bigl(S(Z)\mathbbm 1_{[0,t]}\bigr)
    \qquad
    \text{in }
    \mathbb D^{k,r}(\mathfrak H_1).
\]
Because
\(
    \widetilde S_{\mathcal P}
    \in
    \mathbb D^\infty(\mathfrak H_1)
\)
for every partition \(\mathcal P\), and because
\(\mathbb D^{k,r}(\mathfrak H_1)\) is complete by construction, we conclude that
\(
    (K^*)^{-1}\bigl(S(Z)\mathbbm 1_{[0,t]}\bigr)
    \in
    \mathbb D^{k,r}(\mathfrak H_1).
\)
Since \(k\geq 1\) and \(r>1\) were arbitrary, this proves that
\[
    (K^*)^{-1}\bigl(S(Z)\mathbbm 1_{[0,t]}\bigr)
    \in
    \mathbb D^\infty(\mathfrak H_1).
\]
Applying Meyer's inequality and the identity \eqref{transfer_identity_divergence} once more, we obtain
\[
\begin{aligned}
    \left\|
        \delta^{W,i}\bigl(S(Z)\mathbbm 1_{[0,t]}\bigr)
    \right\|_{\mathbb D^{k-1,r}}
    &=
    \left\|
        \delta^i\left(
            (K^*)^{-1}\bigl(S(Z)\mathbbm 1_{[0,t]}\bigr)
        \right)
    \right\|_{\mathbb D^{k-1,r}}
    \\
    &\le C_{k, r}
    \left\|
        (K^*)^{-1}\bigl(S(Z)\mathbbm 1_{[0,t]}\bigr)
    \right\|_{\mathbb D^{k,r}(\mathfrak H_1)}
    <
    \infty .
\end{aligned}
\]
Thus
\(
    \delta^{W,i}\bigl(S(Z)\mathbbm 1_{[0,t]}\bigr)
    \in
    \mathbb D^{k-1,r}(\mathbb R)
\)
for every \(k\geq 1\) and \(r>1\). Hence
\(
    \delta^{W,i}\bigl(S(Z)\mathbbm 1_{[0,t]}\bigr)
    \in
    \mathbb D^\infty(\mathbb R).
\)
\end{proof}
\begin{proposition}\label{lambda_d_infty}
Assume Condition \ref{condition_covariance_rho_var}, \ref{condition_filtration} and \ref{condition_3} hold. For every \(t\in[0,T]\), we have
\(
    \Lambda_t \in \mathbb{D}^{\infty}(\mathbb{R}).
\)
\end{proposition}
\begin{proof}
The proof follows the argument of Lemma 18 in
\cite{crisan2020high}. Recall that
\(
    \Lambda_t = \exp(\chi_t).
\)
By the preceding proposition, together with the fact that
\[
    \int_0^t S(Z_s)\,ds \in \mathbb{D}^{\infty}(\mathbb{R}),
    \qquad
    S\in C_b^\infty(\mathbb{R}^{d_Z},\mathbb{R}),
\]
we have
\(
    \chi_t \in \mathbb{D}^{\infty}(\mathbb{R}).
\)
Moreover, since the coefficients appearing in the definition of
\(\chi_t\) are bounded, the usual exponential martingale estimates imply
that, for every \(q\geq 1\),
\(
    \mathbb{E}\bigl[\exp(q|\chi_t|)\bigr] < \infty .
\)

We now prove that \(\Lambda_t\in\mathbb{D}^{k,r}(\mathbb{R})\) for
arbitrary \(k\geq 1\) and \(r>1\). By the Faà di Bruno formula, for every
\(n\in\{1,\ldots,k\}\),
\[
    \mathcal D^n \Lambda_t
    =
    \mathcal D^n \exp(\chi_t)
    =
    \exp(\chi_t)
    \sum_{\substack{m_1,\ldots,m_n\geq 0\\
                    \sum_{\ell=1}^n \ell m_\ell=n}}
    \frac{n!}{m_1!\cdots m_n!}
    \prod_{\ell=1}^n
    \left(
        \frac{\mathcal D^\ell \chi_t}{\ell!}
    \right)^{\otimes m_\ell}.
\]
Hence, for every \(r>1\),
\[
\begin{aligned}
    \bigl\|\mathcal D^n\Lambda_t\bigr\|_{L^r(\Omega,\mathfrak{H}_1^{\otimes n})}
    &\leq
    \sum_{\substack{m_1,\ldots,m_n\geq 0\\
                    \sum_{\ell=1}^n \ell m_\ell=n}}
    \frac{n!}{m_1!\cdots m_n!}
    \left\|
        \exp(\chi_t)
        \prod_{\ell=1}^n
        \left(
            \frac{\mathcal D^\ell \chi_t}{\ell!}
        \right)^{\otimes m_\ell}
    \right\|_{L^r(\Omega,\mathfrak{H}_1^{\otimes n})}.
\end{aligned}
\]
Applying H\"older's inequality, using the exponential integrability of
\(\chi_t\), and the fact that
\[
   \mathcal  D^\ell \chi_t \in L^p(\Omega,\mathfrak{H}_1^{\otimes \ell})
    \qquad
    \text{for all } p>1,\ \ell\geq 1,
\]
we obtain
\(
    \bigl\|\mathcal D^n\Lambda_t\bigr\|_{L^r(\Omega,\mathfrak{H}_1^{\otimes n})}
    < \infty .
\)
Since this holds for every \(n\in\{1,\ldots,k\}\), every \(k\geq 1\), and
every \(r>1\), it follows that
\[
    \Lambda_t \in \mathbb{D}^{k,r}(\mathbb{R})
    \qquad
    \text{for all } k\geq 1,\ r>1.
\]
Therefore
\(
    \Lambda_t \in \mathbb{D}^{\infty}(\mathbb{R}),
\)
as claimed.
\end{proof}
As a simple consequence of the example \ref{c_s_gen_example} we obtain the result in the next Proposition which guarantees that the second condition towards proving the smoothness of the density is guaranteed.

\begin{proposition}\label{prop:projection_d_infty_elements} Let $J$ satisfy Condition \ref{condition_covariance_rho_var}. Then the orthogonal projection \(\Pi_{\mathcal H}\) maps
\(\mathbb{D}^{\infty}(\mathfrak{H}_1^{d_J})\) into itself. More precisely, if
\(\Phi \in \mathbb{D}^{\infty}(\mathfrak{H}_1^{d_J})\), then
\(
    \Pi_{\mathcal H}\Phi \in \mathbb{D}^{\infty}(\mathfrak{H}_1^{d_J}).
\)
\end{proposition}

\begin{proof}
By Example~\ref{c_s_gen_example}, the projection associated with
\(\mathcal H\) is deterministic. More precisely, for every \(\omega\),
the operator \(\mathscr{L}(\omega)\) coincides with the orthogonal projection
\(\Pi_{\mathcal H}\) of \(\mathfrak{H}_1^{d_J}\) onto
\[
    \bigoplus_{k=d_B+1}^{d_J}\mathfrak{H}_1^{(k)}.
\]
In particular, \(\Pi_{\mathcal H}\) is a deterministic bounded linear
operator on \(\mathfrak{H}_1^{d_J}\), with operator norm at most one.

Let \(\Phi \in\mathbb{D}^{\infty}(\mathfrak{H}_1^{d_J})\). Since
\(\Pi_{\mathcal H}\) is deterministic and bounded, for every \(n\geq 1\)
we have
\[
    \mathcal D^n(\Pi_{\mathcal H}\Phi)
    =
    \bigl(\operatorname{Id}_{(\mathfrak{H}_1^{d_J})^{\otimes n}}\otimes \Pi_{\mathcal H}\bigr)
    \mathcal D^n \Phi.
\]
Here $\operatorname{Id}_{(\mathfrak{H}_1^{d_J})^{\otimes n}}$ is the identity operator in $(\mathfrak{H}_1^{d_J})^{\otimes n}$.
Therefore, for every \(r>1\),
\[
\begin{aligned}
    \bigl\|\mathcal D^n(\Pi_{\mathcal H}\Phi)\bigr\|_
    {L^r(\Omega,\mathfrak{H}_1^{\otimes n}\otimes\mathfrak{H}_1^{d_J})}
    &=
    \bigl\|
        \bigl(\operatorname{Id}_{(\mathfrak{H}_1^{d_J})^{\otimes n}}\otimes \Pi_{\mathcal H}\bigr)
        \mathcal D^n \Phi
    \bigr\|_
    {L^r(\Omega,(\mathfrak{H}_1^{d_J})^{\otimes n}\otimes\mathfrak{H}_1^{d_J})}
    \\
    &\leq
    \bigl\|\mathcal D^n \Phi\bigr\|_
    {L^r(\Omega,(\mathfrak{H}_1^{d_J})^{\otimes n}\otimes \mathfrak{H}_1^{d_J})}.
\end{aligned}
\]
Moreover,
\[
    \|\Pi_{\mathcal H}\Phi\|_{L^r(\Omega,\mathfrak{H}_1^{d_J})}
    \leq
    \|\Phi\|_{L^r(\Omega,\mathfrak{H}_1^{d_J})}.
\]
Combining the preceding estimates gives
\[
    \|\Pi_{\mathcal H}\Phi\|_{\mathbb{D}^{n,r}(\mathfrak{H}_1^{d_J})}
    \leq
    \|\Phi\|_{\mathbb{D}^{n,r}(\mathfrak{H}_1^{d_J})}
\]
for every \(n\geq 1\) and every \(r>1\). Since
\(\Phi\in\mathbb{D}^{\infty}(\mathfrak{H}_1^{d_J})\), it follows that
\[
    \Pi_{\mathcal H}\Phi\in\mathbb{D}^{n,r}(\mathfrak{H}_1^{d_J})
    \qquad
    \text{for all } n\geq 1,\ r>1.
\]
Hence
\(
    \Pi_{\mathcal H}\Phi\in\mathbb{D}^{\infty}(\mathfrak{H}_1^{d_J}),
\)
as required.
\end{proof}

In order to study the invertibility of the partial Malliavin matrix, we first
recall the representation of the Malliavin derivative of the solution. Since
\(Z=(X,Y)\) takes values in \(\mathbb R^{d_Z}\), the Jacobian of the solution is a
\(d_Z\times d_Z\) matrix. By Section 11 of
\cite{friz2010multidimensional}, for every
\(h=(h^1,\ldots,h^{d_B+d_Y})\in\mathfrak H_1^{d_J}\) and every
\(t\in[0,T]\), one has
\[
    \mathcal D_h Z_t^{(0)}
    =
    \mathscr J_t
    \int_0^t
        \sum_{k=1}^{d_J}
        \left(\mathscr J_s\right)^{-1}
        g_k(Z_s)
        \,d\overline{\psi}(h^k)_s .
\]
Here \(\mathscr J\) denotes the Jacobian of the solution,
which satisfies the linear RDE
\[
    d\mathscr J_t
    =
    Dg(Z_t)\mathscr J_t
    \,d\mathbf J_t,
    \qquad
    \mathscr J_0
    =
    \mathbbm{1}_{d_Z,d_Z}.
\]

Since we only differentiate with respect to the signal-driving directions and
we are interested in the signal component \(X\), we introduce the
\(d_X\times d_X\) signal block of the Jacobian by
\[
    \left(\mathscr J_t^{X}\right)^{ij}
    :=
    \left(\mathscr J_t\right)^{ij},
    \qquad
    i,j=1,\ldots,d_X .
\]
The matrix \(\mathscr J_t^{X}\) is invertible for every \(t\in[0,T]\).
Indeed, \(\mathscr J_t^{\tilde{\mathbf B}}\) takes values in the space of
block triangular matrices with invertible diagonal blocks. Moreover,
\[
    g_{ij}=0,
    \qquad
    i\leq d_X,\quad j>d_X,
\]
so that the partial derivative of the signal component is obtained by
restricting the previous representation to the first \(d_B\) driving
directions.

Additionally, following Hairer and Pillai \cite{hairer2013regularity}, for a smooth bounded
vector field \(S:\mathbb R^{d_Z}\to\mathbb R^{d_X}\), we define
\[
    \mathcal Z_t^S
    :=
    \left(\mathscr J^X_t\right)^{-1}
    S(Z_t),
    \qquad
    t\in[0,T].
\]
Therefore, for \(h=(h^1,\ldots,h_1^{d_B})\in\mathfrak H^{d_B}\), the partial
Malliavin derivative of the signal component can be written as
\[
    \mathcal D_{\mathcal H,h}X_t^{(0)}
    :=
    \left\langle
        \mathcal D_{\mathcal H}X_t^{(0)},h
    \right\rangle_{\mathfrak H_1^{d_B}}
    =
    \mathscr J_t^{X}
    \int_0^t
        \sum_{k=1}^{d_B}
        \mathcal Z_s^{\sigma_k}
        \,d\overline{\psi}(h^k)_s .
\]
Equivalently,
\[
    \mathcal D_{\mathcal H,h}X_t^{(0)}
    =
    \mathscr J^X_t
    \int_0^t
        \sum_{k=1}^{d_B}
        \left(\mathscr J_s^X\right)^{-1}
        \sigma_k(Z_s)
        \,d\overline{\psi}(h^k)_s .
\]

The corresponding partial Malliavin matrix is the \(d_X\times d_X\) matrix
\[
    \Upsilon_{\mathcal H,t}
    :=
    \left(
        \left\langle
            \mathcal D_{\mathcal H}X_t^{(0),i},
            \mathcal D_{\mathcal H}X_t^{(0),j}
        \right\rangle_{\mathfrak H_1^{d_B}}
    \right)_{i,j=1}^{d_X}.
\]
In particular, for every \(v\in\mathbb R^{d_X}\),
\[
    v^\top \Upsilon_{\mathcal H,t}v
    =
    \left\|
        \left\langle
            v,\mathcal D_{\mathcal H}X_t^{(0)}
        \right\rangle_{\mathbb R^{d_X}}
    \right\|_{\mathfrak H_1^{d_B}}^2 .
\]

To establish non-degeneracy of the Malliavin covariance matrix, we shall use
the deterministic Norris lemma of \cite{cass2015smoothness}, which requires a
quantitative lower bound on the directional oscillations of the driving path.
The following proposition introduces the corresponding notion of H\"older
roughness and verifies it, together with the required integrability of its
roughness modulus, for the Gaussian processes considered here.

\begin{lemma}[Hölder roughness and integrability of the roughness modulus]
\label{lem:holder_roughness_from_Kstar}
Let $I$ be a $d_I$-dimensional centered Volterra Gaussian process with i.i.d. components satisfying Conditions \ref{condition_covariance_rho_var} and  \ref{condition_4}. 
Then, there exists a $\varsigma > 0$ such that for every \(\theta\in( 2/p \wedge 1 -\varsigma,  2/p \wedge 1 )\), the process \(I\) is almost surely
\(\theta\)-Hölder rough. More precisely, if
\[
    L_\theta(I)
    :=
    \inf_{\substack{s\in[0,T],\,0<\epsilon\leq T/2\\
                    v\in\mathbb R^{d_I},\,|v|=1}}
    \frac{1}{\epsilon^\theta}
    \sup_{\substack{t\in[0,T]\\
                    \epsilon/2<|t-s|<\epsilon}}
    \bigl|\langle v,I_t-I_s\rangle\bigr|,
\]
then
\(
    L_\theta(I)>0 
\) almost surely.
Moreover, for every \(l>0\),
\(
    \mathbb E\left[L_\theta(I)^{-l}\right]<\infty.
\)
\end{lemma}
\begin{proof}
Deferred to Appendix \ref{app:proof_holder_roughness}.
\end{proof}

\begin{proposition}\label{prop:bracket_control_malliavin_covariance}
Under the conditions \ref{condition_covariance_rho_var},\ref{cameron_martin_condition}, \ref{condition_filtration}, \ref{condition_3} and \ref{condition_4}, for any \(k\in\mathbb N\), there exist constants
\(
    \varpi_k>0,
    C_{t,k}>0,
\)
and a random variable
\(
    H\in\bigcap_{r=1}^{\infty}L^r(\Omega)
\)
such that, for all \(S\in\mathcal S_k\) and all
\(v\in\mathbb R^{d_X}\) with \(|v|=1\), we have
\[
    \bigl\|\langle v,\mathcal Z^S\rangle\bigr\|_{\infty;[0,t]}
    \leq
    C_{t,k} H^{\varpi_k}
    \bigl(v^T\Upsilon_{\mathcal H,t}v\bigr)^{\varpi_k}.
\]
\end{proposition}

\begin{proof}
We first prove the estimate for the vector fields
\(\sigma_i\), \(i=1,\ldots,d_B\). Define for $t \in [0, T]$
\[
    \mathscr v^i_s
    :=
    \left\langle v,\mathcal Z^{\sigma_i}_s\right\rangle,
    \qquad s\in[0,t],
\]
and
\[
    \Gamma_t^i
    :=
    \int_0^t\int_0^t
        \mathscr v^i_u\mathscr v^i_s
    \,dC_B(u,s).
\]
Then, by definition of the partial Malliavin matrix,
\(
    v^T\Upsilon_{\mathcal H,t}v
    =
    \sum_{i=1}^{d_B}\Gamma_t^i.
\)

Assume first that Condition~\ref{condition_4}(A) holds. Then the assumptions
of Corollary~6.10 in \cite{cass2015smoothness} are satisfied. Hence, for
\(i=1,\ldots,d_B\), there exists a constant \(c>0\) such that
\begin{align*}
    \|\mathscr v^i\|_{\infty;[0,t]}
    &\leq
    2\max
    \left[
        \frac{|\Gamma_t^i|^{1/2}}
             {\mathbb E[(B_t^i)^2]^{1/2}},
        \frac{1}{\sqrt c}
        |\Gamma_t^i|^{1/(2+\alpha/p)}
        \|\mathscr v^i\|_{1/p\operatorname{-Hol};[0,t]}^{\alpha p/(2+\alpha p)}
    \right].
\end{align*}
Here $\|\cdot\|_{1/p\operatorname{-Hol}}$ denotes the $1/p$-H\"older seminorm.
By Condition~\ref{condition_4}(A), we have
\(
    \mathbb E[(B_t^i)^2]>0.
\) for every $t > 0$.
Moreover, using Lemma~\ref{regularity_Kf},
\[
    \Gamma_t^i
    =
    \|K^*\mathscr v^i\|_{L^2([0,t])}^2
    \leq
    C
    \left(
        |\mathscr v^i_0|
        +
        \|\mathscr v^i\|_{1/p\operatorname{-Hol};[0,t]}
    \right)^2.
\]
Therefore,
\begin{align*}
    \|\mathscr v^i\|_{\infty;[0,t]}
    &\leq
    C_1
    |\Gamma_t^i|^{1/(2+\alpha/p)}
    \max_{j=1,\ldots,d_B}
    \left(
        |\mathscr v^j_0|
        +
        \|\mathscr v^j\|_{1/p\operatorname{-Hol};[0,t]}^{\alpha p/(2+\alpha p)}
    \right)  \\
    &\leq
    C_1
    \bigl(v^T\Upsilon_{\mathcal H,t}v\bigr)^{1/(2+\alpha/p)}
    \max_{j=1,\ldots,d_B}
    \left(
        |\mathscr v^j_0|
        +
        \|\mathscr v^j\|_{1/p\operatorname{-Hol};[0,t]}^{\alpha p/(2+\alpha p)}
    \right).
\end{align*}
Assume now that Condition~\ref{condition_4} (B) holds in its interpolated
form, with \(\gamma=1/p\). Applying the condition on the interval \([0,t]\),
we obtain
\begin{align*}
    \|\mathscr v^i\|_{\infty;[0,t]}
    &\leq
    C t^{-\eta}
    \left\|K_i^*\bigl(\mathscr v^i\mathbbm 1_{[0,t]}\bigr)
    \right\|_{L^2([0,T])}^{\vartheta}
    \|\mathscr v^i\|_{C^{1/p};[0,t]}^{1-\vartheta}  \\
    &\leq
    C t^{-\eta}
    \bigl(v^T\Upsilon_{\mathcal H,t}v\bigr)^{\vartheta/2}
    \|\mathscr v^i\|_{C^{1/p};[0,t]}^{1-\vartheta}.
\end{align*}
Indeed,
\[
    \left\|K_i^*\bigl(\mathscr v^i\mathbbm 1_{[0,t]}\bigr)
    \right\|_{L^2([0,T])}^2
    \leq
    \sum_{j=1}^{d_B}
    \left\|K_j^*\bigl(\mathscr v^j\mathbbm 1_{[0,t]}\bigr)
    \right\|_{L^2([0,T])}^2
    =
    v^T\Upsilon_{\mathcal H,t}v .
\]
Using Corollary~8.2 in \cite{cass2015smoothness}, there exist
\(\nu>0\) and a random variable
\[
    H\in\bigcap_{r=1}^{\infty}L^r(\Omega)
\]
such that,
\[
    |\mathscr v^i_0|
    +
    \|\mathscr v^i\|_{C^{1/p};[0,t]}^{1-\vartheta}
    \leq
    C_2 H^\nu,
    \qquad
    i=1,\ldots,d_B.
\]
Consequently,
\[
    \|\mathscr v^i\|_{\infty;[0,t]}
    \leq
    C_2 t^{-\eta}
    H^\nu
    \bigl(v^T\Upsilon_{\mathcal H,t}v\bigr)^{\vartheta/2},
    \qquad
    i=1,\ldots,d_B.
\]
Since \(t>0\) is fixed, we may absorb \(t^{-\eta}\) into the deterministic
constant. Thus, after possibly changing the constants and the exponent, we get
\[
    \|\mathscr v^i\|_{\infty;[0,t]}
    \leq
    C_{2,t}
    H^\nu
    \bigl(v^T\Upsilon_{\mathcal H,t}v\bigr)^\nu,
    \qquad
    i=1,\ldots,d_B.
\]
This proves the desired estimate for \(\mathcal Z^{\sigma_i}\),
\(i=1,\ldots,d_B\).
We now pass from the first-order vector fields to the full bracket family.
Let
\(
    S:\mathbb R^{d_X}\times\mathbb R^{d_Y}\to\mathbb R^{d_X}
\)
be a $C^\infty_b$ vector field, then, for every \(v\in\mathbb R^{d_X}\) and every \(u\in[0,T]\), we have
\[
\begin{aligned}
    \left\langle v,\mathcal Z_u^S\right\rangle
    &=
    \left\langle v,S(x_0,0)\right\rangle
    +
    \sum_{i=1}^{d_B}
    \int_0^u
        \left\langle
            v,
            \mathcal Z_r^{[\sigma_i,S]_x}
        \right\rangle
    \,dB_r^i  \\
    &\quad
    +
    \sum_{j=1}^{d_Y}
    \int_0^u
        \left\langle
            v,
            \mathcal Z_r^{\partial_{y_j}S}
        \right\rangle
    \,dY_r^j .
\end{aligned}
\]
The role of Condition~\ref{condition_4} (A) or Condition~\ref{condition_4} (B)
is only to obtain the preceding first-order estimate for
\(\mathcal Z^{\sigma_i}\). The remainder of the argument is the deterministic
induction used in Theorem 5.6 in \cite{cass2015smoothness}. Indeed, the proof only requires
that the processes \(\langle v,\mathcal Z^S\rangle\) satisfy the recursive
rough integral identity above and that the driving rough path \(J\) is
\(\theta\)-Hölder rough for some \(\theta<2/p \wedge 1\). 
Both properties hold in the present setting (see Lemma \ref{lem:holder_roughness_from_Kstar}), with the joint driver
\((B,Y)\) and with the bracket family \((\mathcal S_k)_{k\geq0}\), this implies the existence of two constants
\(
    \varpi_k>0,
    C_{t,k}>0,
\)
and a random variable
\[
    H\in\bigcap_{r=1}^{\infty}L^r(\Omega)
\]
such that, for every \(S\in\mathcal S_k\) and every
\(v\in\mathbb R^{d_X}\) with \(|v|=1\),
\[
    \bigl\|\langle v,\mathcal Z^S\rangle\bigr\|_{\infty;[0,t]}
    \leq
    C_{t,k}H^{\varpi_k}
    \bigl(v^T\Upsilon_{\mathcal H,t}v\bigr)^{\varpi_k}.
\]
This concludes the proof.
\end{proof}

\begin{theorem}[Smooth conditional density]\label{thm:smooth_conditional_density}
Under assumptions \ref{condition_covariance_rho_var},\ref{cameron_martin_condition}, \ref{condition_filtration}, \ref{condition_3} and \ref{condition_4}, for every \(t>0\),
the partial Malliavin matrix \(\Upsilon_{\mathcal H,t}\) is
non-degenerate and satisfies
\[
    \bigl(\det\Upsilon_{\mathcal H,t}\bigr)^{-1}
    \in
    \bigcap_{r=1}^{\infty}L^r(\Omega).
\]
Consequently, the conditional law of \(X_t\) given \(\mathcal Y_t\) admits
a smooth density with respect to the Lebesgue measure for every $t \in (0, T]$.
\end{theorem}

\begin{proof}
We follow closely the proof of Theorem 3.5 in \cite{cass2015smoothness}.
By the partial H\"ormander condition, there exist \(k\in\mathbb N\) and
finitely many vector fields \(S_1,\ldots,S_m\in\mathcal S_k\) such that $\mathbb P$-a.s.
\[
    \operatorname{span}
    \left\{
        S_1(x_0,0),\ldots,S_m(x_0,0)
    \right\}
    =
    \mathbb R^{d_X}.
\]
Hence, by compactness of the unit sphere in \(\mathbb R^{d_X}\),
\[
    a
    :=
    \inf_{|v|=1}
    \sum_{\ell=1}^{m}
    \left|
        \left\langle v,S_\ell(x_0,0)\right\rangle
    \right|
    >
    0.
\]
For each \(\ell=1,\ldots,m\), we have
\[
    \mathcal Z^{S_\ell}_0
    =
    \mathscr J^X_{0}S_\ell(X_0,Y_0)
    =
    S_\ell(x_0,0).
\]
Therefore, for every \(v\in\mathbb R^{d_X}\) with \(|v|=1\),
\[
    \left|
        \left\langle v,S_\ell(x_0,y_0)\right\rangle
    \right|
    =
    \left|
        \left\langle v,\mathcal Z^{S_\ell}_0\right\rangle
    \right|
    \leq
    \bigl\|
        \langle v,\mathcal Z^{S_\ell}\rangle
    \bigr\|_{\infty;[0,t]}.
\]
Thus
\(
    a
    \leq
    \sum_{\ell=1}^{m}
    \bigl\|
        \langle v,\mathcal Z^{S_\ell}\rangle
    \bigr\|_{\infty;[0,t]}.
\)
Applying Proposition~\ref{prop:bracket_control_malliavin_covariance}, we
obtain
\(
    a
    \leq
    m C_{t,k}H^{\varpi_k}
    \bigl(v^T\Upsilon_{\mathcal H,t}v\bigr)^{\varpi_k}.
\)
Consequently,
\[
    v^T\Upsilon_{\mathcal H,t}v
    \geq
    \left(
        \frac{a}{mC_{t,k}}
    \right)^{1/\varpi_k}
    H^{-1}.
\]
Taking the infimum over all \(v\in\mathbb R^{d_X}\) with \(|v|=1\), we get
\(
    \lambda_{\min}(\Upsilon_{\mathcal H,t})
    \geq
    C_{t,k}^{-1}H^{-1},
\)
after changing the deterministic constant. Hence
\(
    \lambda_{\min}(\Upsilon_{\mathcal H,t})^{-1}
    \leq
    C_{t,k}H.
\)
Since
\[
    H\in\bigcap_{r=1}^{\infty}L^r(\Omega),
\]
it follows that
\[
    \lambda_{\min}(\Upsilon_{\mathcal H,t})^{-1}
    \in
    \bigcap_{r=1}^{\infty}L^r(\Omega).
\]
Moreover,
\(
    \det\Upsilon_{\mathcal H,t}
    \geq
    \lambda_{\min}(\Upsilon_{\mathcal H,t})^{d_X},
\)
and therefore
\(
    \bigl(\det\Upsilon_{\mathcal H,t}\bigr)^{-1}
    \leq
    \lambda_{\min}(\Upsilon_{\mathcal H,t})^{-d_X}.
\)
Thus
\[
    \bigl(\det\Upsilon_{\mathcal H,t}\bigr)^{-1}
    \in
    \bigcap_{r=1}^{\infty}L^r(\Omega).
\]
This together with Proposition \ref{lambda_d_infty} and Proposition \ref{prop:projection_d_infty_elements} alllows to use the criterion in \cite{nualart1989partial} presented at the beginning of this section to conclude that  the conditional law of \(X_t\) given \(\mathcal Y_t\) admits
a smooth density with respect to the Lebesgue measure for every $t \in (0, T]$.
\end{proof}

\appendix
\section{Verification of the Conditions for Type I and Type II Fractional Brownian Motion}\label{app:examples}

In this section we collect the proofs regarding the validity of the various conditions imposed through the paper for both type I and type II fractional Brownian motion (defined in Example \ref{ex:fbm} and \ref{ex:rl-fbm}). We recall that in both cases, whenever the parameter $H = 1/2$ those processes are standard Brownian motions. For simplicity we restrain to the case single dimensional case, the case of general case follows in identical fashion. Throughout this section, for ease of exposition, we assume that $I$ is real-valued. The multidimensional case follows componentwise without further difficulty.
\paragraph{Condition \ref{condition_covariance_rho_var}}
\begin{itemize}
    \item \emph{Type I fBm.}
    When $H > 1/2$, from Proposition 5.1.2 in \cite{nualart2009malliavin} we have 
    \[
    \mathbb E\left| I_t - I_s \right|^2 \leq C_H |t-s|^{2H},
    \] via Kolmogorov continuity criterion it is possible to conclude for $(A)$. For $(B)$ , whenever $ H > 1/4$, is verified in Proposition 14 in \cite{friz2010differential} shows that this process has covariance in $\mathcal C^\rho([0, T]^2, \mathbb{R})$ for every \(
    \rho>\frac1{2H},
\)  .
    \item \emph{Type II fBm.}
    When $H > 1/2$, Example 1 in \cite{decreusefond2003stochasticintegrationrespectvolterra} shows that 
    \[
    \mathbb E\left| I_t - I_s \right|^2 \leq C_H |t-s|^{2H},
    \] via Kolmogorov continuity criterion it is possible to conclude for $(A)$. \newline
    For \((B)\), take $H \in (1/4, 1/2]$ and set
\[
    \alpha:=H-\frac12\in(-1/4,0).
\]
The covariance of the type II fractional Brownian motion is
\[
    C_I(s,t)
    =
    C_H
    \int_0^{s\wedge t}
    (s-r)^\alpha(t-r)^\alpha\,dr .
\]
Let \(C\) denote the type I fractional Brownian covariance,
\[
    C(s,t)
    =
    \widetilde C_H
    \frac12\left(s^{2H}+t^{2H}-|t-s|^{2H}\right).
\]
By the Mandelbrot--Van Ness decomposition (Proposition 5.1.2 in \cite{nualart2009malliavin}), i.e. by splitting the two-sided
kernel into its positive-time and negative-time parts, there exist constants
\(a_H,b_H>0\), depending only on the chosen normalisations, such that
\[
    C(s,t)=a_H C_I(s,t)+b_H Q(s,t),
\]
where
\[
    Q(s,t)
    :=
    \int_0^\infty
    \left((s+x)^\alpha-x^\alpha\right)
    \left((t+x)^\alpha-x^\alpha\right)\,dx .
\]
Equivalently,
\[
    C_I(s,t)=a_H^{-1}C(s,t)-a_H^{-1}b_H Q(s,t).
\]
Since multiplication by constants does not affect finite two-dimensional
\(\rho\)-variation, and since the type I fractional Brownian covariance \(C\)
has finite two-dimensional \(\rho\)-variation for every
\(
    \rho>\frac1{2H},
\)
it remains only to show that \(Q\) has finite two-dimensional
\(1\)-variation.
For \(s,t>0\), differentiating under the integral gives
\[
    \partial_s\partial_t Q(s,t)
    =
    \alpha^2
    \int_0^\infty
    (s+x)^{\alpha-1}(t+x)^{\alpha-1}\,dx
    \geq 0.
\]
Moreover,
\[
\begin{aligned}
    \int_0^T\int_0^T \partial_s\partial_t Q(s,t)\,ds\,dt
    &=
    \alpha^2
    \int_0^\infty
    \left(
        \int_0^T (s+x)^{\alpha-1}\,ds
    \right)^2 dx                                      \\
    &=
    \int_0^\infty
    \left((T+x)^\alpha-x^\alpha\right)^2\,dx .
\end{aligned}
\]
This last integral is finite. Indeed, near \(x=0\) the integrand is of order
\(x^{2\alpha}\), which is integrable because \(2\alpha>-1\). Near infinity,
by the mean-value theorem, there is a $C > 0$, with
\[
    (T+x)^\alpha-x^\alpha\leq C x^{\alpha-1},
\]
so the integrand is of order \(x^{2\alpha-2}\), which is integrable.

Therefore
\[
    \int_0^T\int_0^T |\partial_s\partial_t Q(s,t)|\,ds\,dt<\infty.
\]
Since \(\partial_s\partial_t Q\geq0\), this implies that
\(
    Q \in \mathcal C^1([0, T]^2, \mathbb R) \subseteq C^\rho([0, T]^2, \mathbb R) \text{ for every }\rho\geq1.
\)
Combining this with
\[
    C_I=a_H^{-1}C-a_H^{-1}b_H Q
\]
is enough to conclude.
\end{itemize}
\paragraph{Condition \ref{cameron_martin_condition}}
We recall that the only case that is not immediate to verify is the case where the sample paths have finite $p$-variation for $p \geq 2$, that is when $H \leq 1/2$ (see the proof of Condition \ref{condition_covariance_rho_var} above) 
\begin{itemize}
    \item \emph{Type I fBm.} Example 2.8 in \cite{friz2016jain} verifies the embedding for the case $1/4 < H \leq 1/2$.
    \item \emph{Type II fBm.} Example 2.11 in \cite{friz2016jain} verifies the embedding for the case $1/4 < H \leq 1/2$. 
\end{itemize}

\paragraph{Condition \ref{condition_filtration}}
\begin{itemize}
\item \emph{Type I fBm.} 
For any $H \in (0,1)$ the proof of the fact that Condition \ref{condition_filtration} is satisfied can be found in Chapter 5 in \cite{nualart2009malliavin}.
\item \emph{Type II fBm.} 
Let $\alpha = H+1/2$, \(h\in L^2([0,T])\) be orthogonal to \(K(t, \cdot)\) for every \(t\in[0,T]\). Then
\[
    0
    =
    \langle h,K(t, \cdot) \rangle_{L^2}
    =
    c_H\int_0^t (t-r)^{\alpha-1}h(r)\,dr
    =
    I_{0+}^{\alpha}h(t)
\]
for every \(t\in[0,T]\). Since the Riemann--Liouville fractional integral
\(I_{0+}^{\alpha}\) is injective on \(L^2([0,T])\), it follows that \(h=0\).
By the construction of $K^*$ in \ref{sec:VGP}, we can extend $K^*$ to an isometry from $\mathfrak H_1$ to $L^2([0, T])$. Thus, in particular, \(
    \mathbf 1_{[0,t]}\in \operatorname{Range}(K^*) \)
for every \(t\in[0,T]\).
\end{itemize}

\paragraph{Condition \ref{condition_7}}
\begin{itemize}
    \item \emph{Type I fBm.} The covariance is given by
    \[
        C_I(s,t)
        =
        \frac{1}{2}
        \left(
            s^{2H}+t^{2H}-|t-s|^{2H}
        \right).
    \]
    In particular, \(C_I(t,t)=t^{2H}\), so that
    \[
        c_t=t^{2H},
        \qquad
        \dot c_t=2Ht^{2H-1},
    \]
    and therefore \(c\in C([0,T])\cap C^1((0,T])\). Fix \(\delta>0\) and
    take \(\epsilon_\delta<\delta\). For every
    \(a,b\in[0,\epsilon_\delta]\),
    \[
        C_I(t-a,t-b)
        =
        \frac{1}{2}
        \left(
            (t-a)^{2H}
            +(t-b)^{2H}
            -|a-b|^{2H}
        \right),
    \]
    and hence \(t\mapsto C_I(t-a,t-b)\) belongs to
    \(W^{1,1}([\delta,T])\), with
    \[
        \frac{d}{dt}C_I(t-a,t-b)
        =
        Hc_H
        \left(
            (t-a)^{2H-1}
            +(t-b)^{2H-1}
        \right).
    \]
    Since \(t-a\) and \(t-b\) are uniformly bounded away from zero on
    \([\delta,T]\), this derivative converges uniformly, and therefore in
    \(L^1([\delta,T])\), to
    \(2Ht^{2H-1}=\dot c_t\) as \((a,b)\to(0,0)\).

    \item \emph{Type II fBm.} The covariance is given by
    \[
        C_{I}(s,t)
        =
        2Hc_H
        \int_0^{s\wedge t}
            (s-r)^{H-\frac12}
            (t-r)^{H-\frac12}
        \,dr.
    \]
    In particular, \(C_{I}(t,t)=c_Ht^{2H}\), and hence
    \(c\in C([0,T])\cap C^1((0,T])\), with
    \[
        \dot c_t=2Hc_Ht^{2H-1}.
    \]
    Fix \(\delta>0\) and take \(\epsilon_\delta<\delta\). For every
    \(a,b\in[0,\epsilon_\delta]\), a change of variables gives
    \[
        C_{I}(t-a,t-b)
        =
        2Hc_H
        \int_{a\vee b}^{t}
            (u-a)^{H-\frac12}
            (u-b)^{H-\frac12}
        \,du.
    \]
    Therefore \(t\mapsto C_{I}(t-a,t-b)\) belongs to the Sobolev space
    \(W^{1,1}([\delta,T])\), with
    \[
        \frac{d}{dt}C_{I}(t-a,t-b)
        =
        2Hc_H
        (t-a)^{H-\frac12}
        (t-b)^{H-\frac12}.
    \]
    Since \(t-a\) and \(t-b\) are uniformly bounded away from zero on
    \([\delta,T]\), this derivative converges uniformly, and therefore in
    \(L^1([\delta,T])\), to
    \(2Hc_Ht^{2H-1}=\dot c_t\) as \((a,b)\to(0,0)\).
\end{itemize}

It remains to verify Condition~\ref{condition_4}. We begin by recording the fractional estimate that will be used below.

For $\beta,\gamma \in (0, \infty]$ and $\alpha\in\mathbb{R}$, we denote by
$B_{\beta,\gamma}^{\alpha}(\mathbb R)$ the corresponding Besov space on
$\mathbb R$ and by $B_{\beta,\gamma}^{\alpha}([0, 1])$ its restriction to $(0,1)$. Let $\alpha>1/2$. Then there exists a constant $C_\alpha>0$
such that, for every $\Phi\in L^2([0,1])$,
\[
    \|\Phi\|_{B^{- \alpha}_{2,2}([0,1])}
    \leq
    C_\alpha
    \|I_{1-}^{\alpha}\Phi\|_{L^2([0,1])}.
\]
Indeed, since the Besov space $B^{- \alpha}_{2,2}([0,1])$ coincides with
the negative-order Sobolev space $H^{- \alpha}([0,1])$, with equivalent
norms, it suffices to establish the corresponding estimate in
$H^{- \alpha}([0,1])$. Set
\(
    G:=I_{1-}^{\alpha}\Phi.
\)
Then, in the sense of distributions,
\(
    \Phi=D_{1-}^{\alpha}G.
\)
Thus, for every \(\phi\in C_c^\infty(0,1)\), the fractional integration
by parts formula (Corollary 2 in \cite{kilbas1993fractional}) gives
\[
    \langle \Phi,\phi\rangle
    =
    \langle D_{1-}^{\alpha}G,\phi\rangle
    =
    \langle G,D_{0+}^{\alpha}\phi\rangle .
\]
Consequently,
\[
    |\langle \Phi,\phi\rangle|
    \leq
    \|G\|_{L^2([0,1])}
    \|D_{0+}^{\alpha}\phi\|_{L^2([0,1])}.
\]
Introduce \(H_0^\alpha([0,1])\) as the closure of
\(C_c^\infty(0,1)\) in \(H^\alpha([0,1])\).
Since by Theorem 2.1 in \cite{Jin2015} \(D_{0+}^{\alpha}\) is bounded from \(H^\alpha_0([0,1])\) to
\(L^2([0,1])\), we have
\[
    \|D_{0+}^{\alpha}\phi\|_{L^2([0,1])}
    \leq
    C_\alpha \|\phi\|_{H^\alpha([0,1])}.
\]
Therefore
\[
    |\langle \Phi,\phi\rangle|
    \leq
    C_\alpha
    \|I_{1-}^{\alpha}\Phi\|_{L^2([0,1])}
    \|\phi\|_{H^\alpha([0,1])}.
\]
Taking the
supremum over all \(\phi\in C_c^\infty(0,1)\) satisfying
\(\|\phi\|_{H^\alpha([0,1])}\leq 1\), and using the density of
\(C_c^\infty(0,1)\) in \(H_0^\alpha([0,1])\), gives
\[
    \|\Phi\|_{B^{-\alpha}_{2,2}([0,1])} = \|\Phi\|_{H^{-\alpha}([0,1])}
    \leq
    C_\alpha \|I_{1-}^{\alpha}\Phi\|_{L^2([0,1])}.
\]

We also derive an interpolation estimate needed for both Type I and Type II fBm. Let
$a>0$ and $\vartheta\in(0,1)$, and set
\[
    \frac{1}{\beta}=\frac{\vartheta}{2},
    \qquad
    \sigma=(1-\vartheta)\gamma-\vartheta a.
\]
Real interpolation between
$B^{-a}_{2,2}([0,1])$ and
$C^\gamma([0,1])=B^\gamma_{\infty,\infty}([0,1])$ gives
\[
    \|g\|_{B^\sigma_{\beta,\infty}([0,1])}
    \leq
    C
    \|g\|_{B^{-a}_{2,2}([0,1])}^{\vartheta}
    \|g\|_{C^\gamma([0,1])}^{1-\vartheta}.
\]
Moreover,
\(
    \vartheta
    <
    \frac{\gamma}{\gamma+a+\frac12}
\)
is equivalent to $\rho>1/\beta$. The Besov embedding
\[
    B^\sigma_{\beta,\infty}([0,1])
    \hookrightarrow
    L^\infty([0,1])
\]
therefore yields
\begin{equation}
    \label{eq:unit-besov-interpolation}
    \|g\|_{\infty;[0,1]}
    \leq
    C
    \|g\|_{B^{-a}_{2,2}([0,1])}^{\vartheta}
    \|g\|_{C^\gamma([0,1])}^{1-\vartheta}.
\end{equation}
\paragraph{Condition \ref{condition_4}}
\begin{itemize}
    \item \emph{Type I fBm.}
    When $1/4<H\leq 1/2$, Section 4 in
    \cite{cass2015smoothness} verifies that
    Condition~\ref{condition_4} is satisfied. We therefore consider
    $H>1/2$. By Equation 5.5 in Chapter 5 of
    \cite{nualart2009malliavin},
    \[
        \|K^*(f\mathbbm 1_{\mathcal I})\|_{L^2([0,T])}^2
        =
        c_H
        \int_{\mathcal I}\int_{\mathcal I}
        f(u)f(v)|u-v|^{2H-2}\,du\,dv
    \]
    for every interval $\mathcal I\subset[0,T]$. Setting
    \(
        \delta:=H-\frac12,
    \)
    this quadratic form controls the negative Besov norm on the
    unit interval:
    \[
        \|g\|_{B^{-\delta}_{2,2}([0,1])}
        \leq
        C
        \left(
            c_H
            \int_0^1\int_0^1
            g(x)g(y)|x-y|^{2H-2}\,dx\,dy
        \right)^{1/2}.
    \]

    Let $\mathcal I=[s,t]$, set $\ell:=t-s$, and define
    \[
        \Phi(x):=f(s+\ell x),
        \qquad x\in[0,1].
    \]
    A change of variables gives
    \[
        c_H
        \int_0^1\int_0^1
        \Phi(x)\Phi(y)|x-y|^{2H-2}\,dx\,dy
        =
        \ell^{-2H}
        \|K^*(f\mathbbm 1_{\mathcal I})\|_{L^2([0,T])}^2.
    \]
    Furthermore,
    \[
        \|\Phi\|_{\infty;[0,1]}
        =
        \|f\|_{\infty;\mathcal I},
        \qquad
        \|\Phi\|_{C^\gamma([0,1])}
        \leq
        C_T\|f\|_{C^\gamma(\mathcal I)}.
    \]
    Applying \eqref{eq:unit-besov-interpolation} with
    $a=\delta$, for any
    \[
        \vartheta
        <
        \frac{\gamma}{\gamma+\delta+\frac12}
        =
        \frac{\gamma}{\gamma+H},
    \]
    gives
    \[
        \|f\|_{\infty;\mathcal I}
        \leq
        C|\mathcal I|^{-\vartheta H}
        \|K^*(f\mathbbm 1_{\mathcal I})\|_{L^2([0,T])}^{\vartheta}
        \|f\|_{C^\gamma(\mathcal I)}^{1-\vartheta}.
    \]
    Therefore Condition~\ref{condition_4}(B) holds with
    \[
        \eta=\vartheta H,
        \qquad
        \vartheta<\frac{\gamma}{\gamma+H}.
    \]

    \item \emph{Type II fBm.}
    In this case,
    \(
        K^*=I_{T-}^{H+\frac12}.
    \)
    Set
    \(
        \alpha:=H+ 1/2.
    \)
    The fractional estimate proved above gives
    \[
        \|\Phi\|_{B^{-\alpha}_{2,2}([0,1])}
        \leq
        C
        \|I_{1-}^{\alpha}\Phi\|_{L^2([0,1])}.
    \]
    Applying \eqref{eq:unit-besov-interpolation} with
    $a=\alpha$, for any
    \[
        \vartheta
        <
        \frac{\gamma}{\gamma+\alpha+\frac12}
        =
        \frac{\gamma}{\gamma+H+1},
    \]
    yields
    \[
        \|\Phi\|_{\infty;[0,1]}
        \leq
        C
        \|I_{1-}^{\alpha}\Phi\|_{L^2([0,1])}^{\vartheta}
        \|\Phi\|_{C^\gamma([0,1])}^{1-\vartheta}.
    \]

    Now let $\mathcal I=[s,t]$, set $\ell:=t-s$, and define
    \[
        \Phi(x):=f(s+\ell x),
        \qquad x\in[0,1].
    \]
    Then
    \[
        \|\Phi\|_{\infty;[0,1]}
        =
        \|f\|_{\infty;\mathcal I},
        \qquad
        \|\Phi\|_{C^\gamma([0,1])}
        \leq
        C_T\|f\|_{C^\gamma(\mathcal I)}.
    \]
    By the scaling of the right-sided fractional integral,
    \[
        I_{1-}^{\alpha}\Phi(x)
        =
        \ell^{-\alpha}
        I_{t-}^{\alpha}f(s+\ell x),
    \]
    and consequently
    \[
        \|I_{1-}^{\alpha}\Phi\|_{L^2([0,1])}
        =
        \ell^{-\alpha-\frac12}
        \|I_{t-}^{\alpha}f\|_{L^2(\mathcal I)}.
    \]
    Since $\alpha+\frac12=H+1$ and
    \[
        I_{t-}^{\alpha}f
        =
        I_{T-}^{\alpha}(f\mathbbm 1_{\mathcal I})
        \quad\text{on }\mathcal I,
    \]
    it follows that
    \[
        \|I_{1-}^{\alpha}\Phi\|_{L^2([0,1])}
        \leq
        \ell^{-(H+1)}
        \|K^*(f\mathbbm 1_{\mathcal I})\|_{L^2([0,T])}.
    \]
    Substitution into the unit-interval estimate gives
    \[
        \|f\|_{\infty;\mathcal I}
        \leq
        C|\mathcal I|^{-\vartheta(H+1)}
        \|K^*(f\mathbbm 1_{\mathcal I})\|_{L^2([0,T])}^{\vartheta}
        \|f\|_{C^\gamma(\mathcal I)}^{1-\vartheta}.
    \]
    Therefore Condition~\ref{condition_4}(B) holds with
    \[
        \eta=\vartheta(H+1),
        \qquad
        \vartheta<\frac{\gamma}{\gamma+H+1}.
    \]
\end{itemize}

The next example shows how to construct a smooth approximating sequence under Condition \ref{approximation_class_rpde}

\begin{example}
\label{ex:pathwise_variance_regularisation}

Let \(\rho_0\in C_c^\infty((0,1))\) be non-negative and satisfy
\(
    \int_0^1\rho_0(r)\,dr=1.
\)
Extend \(\hat J\) to \((-\infty,T]\) by setting
\[
    \hat J_t= \hat J_0=0,
    \qquad t\leq0.
\]
Let \(\chi\in C^\infty([0,\infty),[0,1])\) satisfy
\[
    \chi(r)=0
    \quad\text{for }r\in[0,1/2],
    \qquad
    \chi(r)=1
    \quad\text{for }r\geq1,
\]
and define
\[
    \lambda_\epsilon(t)
    :=
    \chi\left(\frac{t}{\epsilon}\right),
    \qquad
    t\in[0,T].
\]
Set
\[
    \widetilde J_t^\epsilon
    :=
    \lambda_\epsilon(t)
    \int_0^1
        \rho_0(r)\hat J_{t-\epsilon r}
    \,dr.
\]
Then 
\(
    \widetilde J^\epsilon
    \in C^\infty([0,T],\mathbb R^{d_I})
    \text{ almost surely}.
\)
Moreover, by continuity of the paths of \(J\),
\(
    \|\widetilde J^\epsilon-\hat J\|_{\infty;[0,T]}
    \to0
    \text{ almost surely}.
\)
Moreover, since for every $\epsilon > 0$ the covariance of $(\hat J, \tilde J ^\epsilon)$ has finite $\rho$-variation and $C_{\tilde J^\epsilon}$ converges in the sup norm to $C_J$, then by Theorem 37 in \cite{friz2010differential} and convergence of the It\^o integrals, for every \(\ell<\infty\),
\[
    \mathbb{E}\left[
        \left\|
            \mathscr S_\kappa(\widetilde J^\epsilon) - 
            \hat{\mathbf J}
        \right|^\ell
    \right]
    \to0.
\]
To make the approximating variance strictly positive, let \(z\) be a
standard normal random variable independent of \(\hat J\), and let
\(\eta_\epsilon>0\) satisfy \(\eta_\epsilon\to0\). Define
\[
    \hat J_t^\epsilon
    :=
    \widetilde J_t^\epsilon+\eta_\epsilon tz.
\]
Then \(\hat J^\epsilon\) is centred Gaussian and
\(
    \hat J^\epsilon
    \in C^\infty([0,T],\mathbb R^{d_I})
    \text{ almost surely}.
\)
Its variance is
\[
    c_t^\epsilon
    :=
    \mathbb E\left[
        \left(\hat J_t^\epsilon\right)^2
    \right]
    =
    \mathbb E\left[
        \left(\widetilde J_t^\epsilon\right)^2
    \right]
    +
    \eta_\epsilon^2t^2.
\]
Consequently,
\[
    c^\epsilon\in C^\infty([0,T]),
    \qquad
    c_0^\epsilon=0,
    \qquad
    c_t^\epsilon>0
    \quad\text{for }t>0.
\]
Moreover,
\[
    c^\epsilon\to c
    \qquad\text{uniformly on }[0,T].
\]
Indeed, the uniform convergence of \(\widetilde J^\epsilon\) to \(\hat J\),
together with Gaussian moment bounds, gives uniform convergence of the
corresponding variances, while
\[
    \sup_{t\in[0,T]}\eta_\epsilon^2t^2
    \leq
    \eta_\epsilon^2T^2
    \to0.
\]
By Condition~\ref{condition_7}, for every \(\delta>0\) there exists
\(\epsilon_\delta>0\) such that, whenever
\(
    0<\epsilon<\min\{\delta,\epsilon_\delta\},
\)
the map
\[
    t\to C_{\hat J}(t-\epsilon r,t-\epsilon q)
\]
belongs to \(W^{1,1}([\delta,T])\) for every \(r,s\in[0,1]\). Since
\(\lambda_\epsilon(t)=1\) on \([\delta,T]\), it follows that, for
almost every \(t\in[\delta,T]\),
\[
    \dot c_t^\epsilon
    =
    \int_0^1\int_0^1
        \rho_0(r)\rho_0(s)
        \frac{d}{dt}
        C_{\hat J}(t-\epsilon r,t-\epsilon s)
    \,dr\,ds
    +
    2\eta_\epsilon^2t.
\]
Consequently,
\begin{align*}
    \int_\delta^T
        |\dot c_t^\epsilon-\dot c_t|
    \,dt
    &\leq
    \int_0^1\int_0^1
        \rho_0(r)\rho_0(s)
        \left\|
            \frac{d}{dt}
            C_{\hat J}(\,\cdot-\epsilon r,\cdot-\epsilon s\,)
            -
            \dot c
        \right\|_{L^1([\delta,T])}
    \,dr\,ds
    \\
    &\quad
    +
    \eta_\epsilon^2(T^2-\delta^2).
\end{align*}
By Condition~\ref{condition_7}, the first term converges to zero as
\(\epsilon\to0\), while the second converges to zero since
\(\eta_\epsilon\to0\). Hence, for every \(\delta>0\),
\[
    \dot c^\epsilon
    \to
    \dot c
    \qquad\text{in }L^1([\delta,T]).
\]
Finally, the additional path \(t\mapsto\eta_\epsilon tz\) has finite
variation and
\[
    \|\eta_\epsilon(\cdot)z\|_{1\text{-var};[0,T]}
    =
    \eta_\epsilon T|z|
    \to0
\]
in every \(L^\ell\). Since \(
    \mathscr S(\widetilde J^\epsilon)
    \to
    \hat{\mathbf J}
\)
in \(L^\ell\) with respect to the \(p\)-variation rough-path metric,
continuity of Young translation yields
\[
    \mathbb{E}\left[
        \left\|
            \mathscr{S}_{\kappa }(\hat J^\epsilon) - 
            \hat{\mathbf J}
        \right\|_{p}^r
    \right]
    \to0
\]
for every \(r<\infty\).
\end{example}

\section{Proof of the well-posedness of the signal-observation system}\label{app:proofs}

This appendix develops the well-posedness theory underlying
Theorem~\ref{weaK_solution_existence}. The analysis is divided into two parts.
We first study \eqref{sensor_equation} pathwise, establishing local existence
and global uniqueness by means of a ball-preserving argument in a suitable
space of controlled rough paths. We then address the probabilistic aspects of
the construction and adapt the classical Yamada--Watanabe argument to the
rough-path setting in order to obtain an adapted solution defined outside a
single null set. This second step is only needed in the genuinely rough
regime. When the driving signal has finite \(p\)-variation for some \(p<2\),
so that the equation can be interpreted in the Young sense, the required
adapted solution follows directly from the pathwise construction.
\subsection{Local existence and global uniqueness of RDEs with Volterra drift}

In order to prove well-posedness of the signal--observation system considered in this work, we first isolate a general class of rough Volterra equations. The signal--observation dynamics considered later will be obtained as a particular instance of this class, for a suitable choice of state variable, drift and diffusion coefficients. At this level of generality, we prove local existence and global uniqueness in a pathwise sense. In the subsequent subsection, we then strengthen this result in the particular case of the filtering dynamics, where the additional structure of the equation allows us to obtain global well-posedness.

Let \((\Omega,\mathcal F,\mathbb P)\) be a probability space, and let \(I\) be a \(d\)-dimensional Volterra Gaussian process with kernel \(K\), satisfying Condition \ref{condition_covariance_rho_var} and \ref{cameron_martin_condition}. In particular there exists a set of full measure on which \(I\) admits the rough path lift \(\mathbf I \in \Omega G_p(\mathbb{R}^d)\), $p \in [1, 4)$.  On this set of full measure set we consider, for each fixed realisation \(\omega\), the equation
\begin{align}\label{proxy_eqtn}
    Z_t
    =
    z_0
    +
    Kb\bigl(Z^{(0)}\bigr)(t)
    +
    \int_0^t g(Z)_s\,d\mathbf I_s,
\end{align}
where \(z_0\in\mathbb R^{d_Z}\) and the precise regularity assumptions on \(b\) and \(g\) will be specified below. The unknown \(Z\) will be interpreted as a controlled rough path above \(\mathbf I\), and the integral in \eqref{proxy_eqtn} is understood in the corresponding rough-path sense. For notational convenience, we set $\kappa := \floor{p}$, moreover, we shall write \(Kb(Z)\) instead of \(Kb\bigl(Z^{(0)}\bigr)\) whenever no confusion can arise.

Recall that, in Section~\ref{section_setup}, in order to regard \(Z\) as a
controlled path in the sense of Definition~\ref{operator_drift_solution}, we imposed
Condition~\ref{cameron_martin_condition} to ensure sufficient regularity of the
deterministic Volterra term \(Kb\). The following proposition makes this claim
precise.

\begin{lemma}\label{regularity_Kf}
Assume Condition \ref{cameron_martin_condition} holds. Then, for every
\(b \in L^2([0,T],\mathbb R^d)\),
\(
    Kb \in \mathcal C^{q}([0,T],\mathbb R^d).
\)
\end{lemma}

\begin{proof}
We first consider the case \(d=1\). The general case follows componentwise.

Let \(h=Kb\). We show that \(h\in \mathfrak H\). Let
\(\mathcal E\) be as in \eqref{span_e}. Define \(L_h:\mathcal E\to \mathbb R\) by
\[
    L_h(\varphi)
    :=
    \sum_i a_i h(t_i),
    \qquad
    \varphi=\sum_i a_i C_{t_i}.
\]
Using identities \eqref{K(t)_star} and \eqref{cameron_martin_condition}
we have
\[
    \|\varphi\|_{\mathfrak H}^2
    =
    \sum_{i,j} a_i a_j C_B(t_i,t_j)
    =
    \left\|
        \sum_i a_i K(t_i,\cdot)
    \right\|_{L^2}^2 .
\]
Moreover, since \(h=Kb\),
\(
    h(t_i)
    =
    \langle K(t_i,\cdot),b\rangle_{L^2}.
\)
Therefore, by Cauchy--Schwarz,
\begin{align*}
    |L_h(\varphi)|
    &=
    \left|
        \sum_i a_i h(t_i)
    \right|  \\
    &=
    \left|
        \left\langle
            \sum_i a_i K(t_i,\cdot),
            b
        \right\rangle_{L^2}
    \right|  \\
    &\leq
    \|b\|_{L^2}
    \left\|
        \sum_i a_i K(t_i,\cdot)
    \right\|_{L^2}  \\
    &=
    \|b\|_{L^2}\|\varphi\|_{\mathfrak H}.
\end{align*}
Hence \(L_h\) is bounded on \(\mathcal E\), and therefore extends uniquely to
a bounded linear functional on \(\mathfrak H\). By the Riesz representation
theorem, there exists \(g\in\mathfrak H\) such that
\[
    L_h(\varphi)
    =
    \langle g,\varphi\rangle_{\mathfrak H},
    \qquad
    \varphi\in \mathfrak H.
\]
In particular, taking \(\varphi=C_t\), the reproducing property gives
\[
    h(t)
    =
    L_h(C_t)
    =
    \langle g,C_t\rangle_{\mathfrak H}
    =
    g(t),
    \qquad t\in[0,T].
\]
Thus \(h=g\) pointwise, and hence \(h\in\mathfrak H\).

By Condition \ref{cameron_martin_condition}, the embedding
\[
    \mathfrak H
    \hookrightarrow
    \mathcal C^{q\text{-var}}([0,T],\mathbb R)
\]
is continuous. Therefore \(h=Kb\in \mathcal C^{q\text{-var}}([0,T],\mathbb R)\).

For \(d>1\), the same argument applies to each component of \(b\), yielding
\[
    Kb\in \mathcal C^{q\text{-var}}([0,T],\mathbb R^d).
\]
\end{proof}
We begin by proving that, under appropriate regularity assumptions on the
coefficients \(b\) and \(g\), equation \eqref{proxy_eqtn} is locally well posed.
The proof is based on a contraction argument on a suitable ball, following a
strategy similar to that of \cite{deya2009rough}. The main difference is that,
in the present setting, the Volterra structure enters through the drift term
\(Kb\), rather than through the coefficient of the rough signal.

\begin{proposition}\label{local_solution_sensor_equation}
Suppose Conditions \ref{condition_covariance_rho_var} and \ref{cameron_martin_condition} hold, suppose in addition that
\(g \in C_b^{\kappa+1}(\mathbb  R^{d_Z},\newline
\mathrm{Hom}(\mathbb{R}^{d_I},\mathbb  R^{d_Z}))\) and
\(b\in C_b^2(\mathbb  R^{d_Z},\mathbb  R^{d_Z})\). Then equation
\eqref{proxy_eqtn} admits a unique local solution
\(
    Z \in \mathscr{D}_{\mathbf I,p_*}(\mathbb  R^{d_Z})
\)
define an interval $[0, T_0]$ for some \(T_0\in(0,T]\), depending on \(z_0\), \(f\), \(g\), and
\(\mathbf I\).
\end{proposition}
\begin{proof}
We will show the result for the case where $p \in [3, 4)$, the case $p \in (1, 3)$ will then follow from similar calculations.
Fix $T_0 \in (0, T]$ and a tuple of positive numbers $\nu = (\nu_0, \nu_1, \nu_2)$, both quantities will be tuned later in the proof. Let $p_* < \check{p} < 4$ and consider $\mathbf{I} \in \Omega G_{\check p}([0, T_0], \mathbb R^{d_I})$.
For $z_0 \in \mathbb R^{d_Z}$, define the ball 
\[
\mathcal{B}^{\nu}_{T_0, z_0} := \left\{ Z \in \mathscr{D}_{\mathbf{I}, \check{p}}(\mathbb  R^{d_Z}) : Z_0 = (z_0, \sigma(z_0), D\sigma(z_0) \sigma(z_0)),  \| R^{(i)} \|_{\check p / (\kappa-i); [0, T_0]} \leq \nu_i,\, i \in [0:\kappa] \right\}.
\]
Consider the map
\[
\mathcal{M}_t(Z)
=
\left(
    z_0
    +
    Kb(Z^{(0)})(t)
    +
    \int_0^t g(Z)_u\,d\mathbf I_u,
    \,
    g(Z^{(0)}_t),
    \,
    Dg(Z^{(0)}_t)Z^{(1)}_t
\right),
\qquad
Z\in \mathcal{B}^{\nu}_{T_0, z_0}.
\]
We first observe that \(\mathcal M(Z)\) takes values in
\(\mathscr{D}_{\mathbf I,\check p}(\mathbb R^{d_Z})\). Indeed, by Lemma
\ref{regularity_Kf}, since $Kb(Z) \in \mathcal C^q$, 
\(
    \bigl(Kb(Z),0,0\bigr)
\)
defines a \(\mathbf I\)-controlled path. On the other hand, by the standard
stability of controlled paths under composition and rough integration (see Section 2 in \cite{cass2022combinatorial}),
\[
\left(
    z_0+\int_0^\cdot g(Z)_u\,d\mathbf I_u,
    \,
    g(Z^{(0)}),
    \,
    Dg(Z^{(0)})Z^{(1)}
\right)
\in
\mathscr{D}_{\mathbf I,\check p}(\mathbb R^{d_Z}).
\]
Therefore we may decompose
\[
\mathcal M(Z)
=
\left(
    z_0+\int_0^\cdot g(Z)_u\,d\mathbf I_u,
    \,
    g(Z^{(0)}),
    \,
    Dg(Z^{(0)})Z^{(1)}
\right)
+
\bigl(Kb(Z),0,0\bigr).
\]
Since \(\mathscr{D}_{\mathbf I,\check p}(\mathbb R^{d_Z})\) is a vector space, and
in fact a Banach space under the controlled-path norm, it follows that
\(
    \mathcal M(Z)\in \mathscr{D}_{\mathbf I,\check p}(\mathbb R^{d_Z}).
\)
We now show that, after choosing the radii appropriately and possibly
decreasing \(T_0\), the map \(\mathcal M\) maps \(\mathcal{B}^{\nu}_{T_0, z_0}\) into itself.
Throughout this argument \(C\) denotes a constant depending only on
\(\check p,g,f,K,T\), and may change from line to line. Set
\[
    \alpha(T_0)
    :=
    \sup_{Z\in \mathcal{B}^{\nu}_{T_0, z_0}}
    \|Kb(Z)\|_{\check p/\kappa;[0,T_0]} .
\]
By Lemma \ref{regularity_Kf} and the bounds defining
\(\mathcal{B}^{\nu}_{T_0, z_0}\), the family
\[
    \mathcal K:=\{Kb(Z):Z\in \mathcal{B}^{\nu}_{T_0, z_0}\}
\]
is uniformly equicontinuous and bounded in
\(\mathcal C^{q}\), for some \(q<\check p/\kappa\). Hence
\(\mathcal K\) is relatively compact in
\(\mathcal C^{\check p/\kappa}\) (Proposition 1.15 in \cite{lyons2007differential}), and the variation seminorms are
uniformly small on shrinking intervals. Therefore
\[
    \alpha(T_0):=
    \sup_{Z\in \mathcal{B}^{\nu}_{T_0, z_0}}
    \|Kb(Z)\|_{\check p/\kappa;[0,T_0]}
    \to0
    \qquad \text{as }T_0\downarrow0.
\]
Together with the continuity of the rough path variation control,
\[
    \|\mathbf I\|_{\check p;[0,T_0]}\to 0.
\]
Let \(Z\in \mathcal{B}^{\nu}_{T_0, z_0}\). 
From calculations in Section 4 in  \cite{friz2018differential}, the remainder $R^{\mathcal{M}(Z)}$ satisfies:
\begin{align*}
\|R^{\mathcal{M}(Z), (1)}\|_{\check p /2; [0, T_0]}
&\leq
C
\biggl(
    \|R^{Z, (0)}\|_{\check p/\kappa; [0, T_0]}
    + \|\mathbf I^{(1)}\|_{\check p; [0, T_0]}
\\
&\hspace{4cm}
    + \Bigl(
        \|R^{Z, (0)}\|_{\check p/\kappa; [0, T_0]}
        + \|\mathbf I\|_{\check p; [0, T_0]}
    \Bigr)^2
\biggr)\\
&\leq
C\left(
    \nu_0+\|\mathbf I\|_{\check p;[0,T_0]}
    +(\nu_0+\|\mathbf I\|_{\check p;[0,T_0]})^2
\right), 
\\[0.4em]
\|R^{\mathcal{M}(Z), (2)}\|_{\check p; [0, T_0]}
&\leq
C
\biggl(
    \|R^{Z, (0)}\|_{\check p/\kappa; [0, T_0]}
    + \|\mathbf I^{(1)}\|_{\check p; [0, T_0]}
\\
&\hspace{4cm}
    + \Bigl(
        \|R^{Z, (0)}\|_{\check p/\kappa; [0, T_0]}
        + \|\mathbf I\|_{\check p; [0, T_0]}
    \Bigr)^2
\\
&\hspace{4cm}
    + \|\mathbf I^{(1)}\|_{\check p; [0, T_0]}
    + \|R^{Z, (1)}\|_{\check p/(\kappa-1); [0, T_0]}
\biggr)\\
&\leq
C\left(
    \nu_0+\nu_1+\|\mathbf I\|_{\check p;[0,T_0]}
    +(\nu_0+\|\mathbf I\|_{\check p;[0,T_0]})^2
\right).
\end{align*}
It remains to control the zeroth remainder. By definition,
\[
\begin{aligned}
R^{\mathcal M(Z),(0)}_{s,t}
&=
Kb(Z)(t)-Kb(Z)(s)
+
\int_s^t g(Z)_u\,d\mathbf B_u
-
g(Z_s)\mathbf I^{(1)}_{s,t}
-
Dg(Z_s)Z^{(1)}_s\mathbf I^{(2)}_{s,t}.
\end{aligned}
\]
Using the \enquote{sewing estimate} for the rough integral (Section 2 in \cite{cass2022combinatorial}) and the expansion
of \(g(Z)\) against $\mathbf I$, one obtains
\[
\begin{aligned}
\|R^{\mathcal M(Z),(0)}\|_{\check p/\kappa;[0,T_0]}
&\leq
\alpha(T_0)
+
C\|\mathbf I\|_{\check p; [0, T_0]}
\|R^{g(Z),(0)}\|_{\check p/\kappa;[0,T_0]}
\\
&\quad
+
C\|\mathbf I\|_{\check p; [0, T_0]}
\|R^{g(Z),(1)}\|_{\check p/2;[0,T_0]}
+
C\|\mathbf I\|_{\check p; [0, T_0]}
\|g(Z)^{(2)}\|_{\check p;[0,T_0]} .
\end{aligned}
\]
Since
\[
    R^{g(Z),(0)}=R^{\mathcal M(Z),(0)},
    \qquad
    R^{g(Z),(1)}=R^{\mathcal M(Z),(1)},
\]
and
\[
    g(Z)^{(2)}
    =
    D^2g(Z)(Z^{(1)})^{\otimes 2}
    +
    Dg(Z)Z^{(2)},
\]
the boundedness assumptions defining \(\mathcal{B}^{\nu}_{T_0, z_0}\), together with the
boundedness of \(g\) and its derivatives, imply
\[
    \|g(Z)^{(2)}\|_{\check p;[0,T_0]}\leq C.
\]
Therefore,
\[
\begin{aligned}
\|R^{\mathcal M(Z),(0)}\|_{\check p/\kappa;[0,T_0]}
&\leq
\alpha(T_0)
+
C\|\mathbf I\|_{\check p; [0, T_0]}
\|R^{\mathcal M(Z),(0)}\|_{\check p/\kappa;[0,T_0]}
\\
&\quad
+
C\|\mathbf I\|_{\check p; [0, T_0]}
\left[
    1+\nu_0+\|\mathbf I\|_{\check p; [0, T_0]}+(\nu_0+\|\mathbf I\|_{\check p; [0, T_0]})^2
\right].
\end{aligned}
\]
Choose \(T_0\) small enough that
\[
    C\|\mathbf I\|_{\check p; [0, T_0]}\leq \frac12 .
\]
Then the second term on the right-hand side can be absorbed into the left-hand
side, giving
\[
\begin{aligned}
\|R^{\mathcal M(Z),(0)}\|_{\check p/\kappa;[0,T_0]}
&\leq
2\alpha(T_0)
+
C\|\mathbf I\|_{\check p; [0, T_0]}
\left[
    1+\nu_0+\|\mathbf I\|_{\check p; [0, T_0]}+(\nu_0+\|\mathbf I\|_{\check p; [0, T_0]})^2
\right].
\end{aligned}
\]
Since the right-hand side tends to \(0\) as \(T_0\downarrow0\), we may
decrease \(T_0\) further so that
\[
    \|R^{\mathcal M(Z),(0)}\|_{\check p/\kappa;[0,T_0]}
    \leq \nu_0 .
\]
We now choose the remaining radii successively. First choose \(\nu_1\) such
that
\[
    \nu_1>C(\nu_0+\nu_0^2).
\]
Then, by decreasing \(T_0\) if necessary, the first estimate above yields
\[
    \|R^{\mathcal M(Z),(1)}\|_{\check p/2;[0,T_0]}
    \leq \nu_1 .
\]
Next choose \(\nu_2\) such that
\[
    \nu_2>C(\nu_0+\nu_1+\nu_0^2).
\]
Again decreasing \(T_0\), if necessary, the second estimate yields
\[
    \|R^{\mathcal M(Z),(2)}\|_{\check p;[0,T_0]}
    \leq \nu_2 .
\]
It remains to verify the uniform bounds defining the ball. Since
\[
    \mathcal M(Z)_t
    =
    z_0+Kb(Z)(t)+\int_0^t g(Z)_u\,d\mathbf I_u ,
\]
the rough integral estimate (Section 2 in \cite{cass2022combinatorial}) gives
\[
    \|\mathcal M(Z)-z_0\|_{\infty;[0,T_0]}
    \leq
    \|Kb(Z)\|_{\infty;[0,T_0]}
    +
    C\|\mathbf I\|_{\check p; [0, T_0]}.
\]
The right-hand side can be made arbitrarily small by decreasing \(T_0\).
Moreover,
\[
    \mathcal M(Z)^{(1)}=g(Z),
    \qquad
    \mathcal M(Z)^{(2)}=Dg(Z)Z^{(1)} ,
\]
and these quantities are uniformly bounded by the boundedness of \(g\),
\(Dg\), and the bounds defining \(\mathcal{B}^{\nu}_{T_0, z_0}\). Thus, choosing \(M\) large
enough and then \(T_0\) sufficiently small, all the defining bounds of
\(\mathcal{B}^{\nu}_{T_0, z_0}\) are preserved.

Consequently,
\[
    \mathcal M(\mathcal{B}^{\nu}_{T_0, z_0})\subseteq \mathcal{B}^{\nu}_{T_0, z_0}.
\]
For the stability part we 
begin by defining the equivalent controlled path norm
\[
\|Z\|^{(\varsigma)}_{\mathscr{D}_{\mathbf I, \check p}; [s,t]} = |Z_s| + \sum_{j=0}^{\kappa -1 } \varsigma_j \|R^{(j)}\|_{\check p / (\kappa - j); [s,t]},
\]
where $\varsigma = (\varsigma_0, \varsigma_1, \varsigma_2)$  will be tuned later in the proof.

Let $Z, \tilde Z \in B_{T_0, z_0}$.  We invoke the estimates in Lemma  4.11 of \cite{friz2018differential} to write for some positive constant $C$  depending only on
\(\check p,g,f,K,T, \nu\), and may change from line to line
\begin{align*}
&\left\|R^{\mathcal M (Z),(0)} - R^{\mathcal M (\tilde Z),(0)}\right\|_{\check p /\kappa ;[0,T_0]} \\
&\qquad \leq \|Kb(Z) - Kb(\tilde Z)\|_{\check p /k; [0; T_0]}  +
C\|\mathbf I\|_{\check p;[0,T_0]}
\sum_{j =0 }^{\kappa -1}
\left\|R^{Z,(j)} - R^{\tilde Z,(j)}\right\|_{\check p/ (\kappa - j);[0,T_0]}, \\
&\left\|R^{\mathcal M(Z),(k)}\ -R^{\mathcal M(\tilde Z),(k)} \right\|_{\check p / (\kappa -k) ;[0,T_0]}\\
&\quad \leq
C\left( \sum_{j = 0}^{k-1}
\left\|R^{Z,(j)} - R^{\tilde Z,(j)}\right\|_{\check p /(\kappa - j+1);[0,T_0]} + 
\|\mathbf I\|_{\check p;[0,T_0]}
\sum_{j = k}^{\kappa}
\left\| R^{Z,(j)} - R^{\tilde Z,(j)}\right\|_{\check p / (\kappa - j);[0,T_0]}
\right),
\end{align*}
where $k = 1,2$.

We now choose the weights. Normalise \(\varsigma_0=1\). Choose first
\(\varsigma_1>0\) and then \(\varsigma_2>0\) sufficiently small so that
\[
    C(\varsigma_1+\varsigma_2)\leq \frac18,
    \qquad
    C\varsigma_2\leq \frac{\varsigma_1}{8}.
\]
With these weights fixed, decrease \(T_0\) so that
\[
    \alpha(T_0)
    +
    C\| \mathbf I \|_{\check p; [0, T_0]}\left(1+\varsigma_1^{-1}+\varsigma_2^{-1}\right)
    \leq \frac18 .
\]
Multiplying the three estimates above by
\(\varsigma_0,\varsigma_1,\varsigma_2\), respectively, and summing, we obtain
\[
    \|\mathcal M(Z) - \mathcal M(\tilde Z)\|_{\mathscr{D}_{\mathbf I, \tilde p}; [0,T_0]}
    \leq
    \frac12 \|Z - \tilde Z\|_{\mathscr{D}_{\mathbf I, \tilde p}; [0,T_0]}.
\]
Thus, after the above choice of weights and after possibly decreasing
\(T_0\), the map \(\mathcal M\) is a contraction on \(\mathcal{B}^{\nu}_{T_0, z_0}\). Since \(\mathbf B\) has finite \(p_*\)-variation, the solution constructed
above may be viewed a posteriori at the regularity of the driver. Indeed,
from the equation
\[
    Z_t
    =
    z_0+Kb(Z)(t)+\int_0^t g(Z)_u\,d\mathbf I_u ,
\]
the controlled expansion of the rough integral is governed by the
\(p_*\)-variation of \(\mathbf B\), while Lemma
\ref{regularity_Kf} gives the required regularity of the drift term
\(Kb(Z)\). Therefore the same remainder estimates, now computed with
\(p_*\) in place of \(\check p\), show that
\[
    Z\in \mathscr D_{\mathbf I,p_*}(\mathbb R^{d_Z}).
\]
This is enough to conclude
\end{proof}

\begin{remark}\label{rem:local_not_global_general_rough}
The preceding argument gives local existence and uniqueness, but it does
not, in general, give global existence when \(p \in [2,4)\). The obstruction is
the lack of a uniform control on the radii of the balls used in the successive
local fixed-point arguments. Indeed, in the genuinely rough regime, the norm of
the image of a controlled path under the solution map grows polynomially in the
controlled-path norm. For instance, when \(2<p<3\), if \(Z^{[T_i]}\) denotes the
solution constructed on \([0,T_i]\), then on a further interval
\([T_i,T_i+\zeta]\) one obtains an estimate of the form
\[
    \|\mathcal M(Z)\|_{\mathscr D_{\mathbf I,\check p};[0,T_i+\zeta]}
    \leq C\left(
    \|Z^{[T_i]}\|_{\mathscr D_{\mathbf I,\check p};[0,T_i]}
    +
    1
    +
    \|\mathbf I\|_{p;[T_i,T_i+\zeta]}^{\lambda}
    \|Z\|_{\mathscr D_{\mathbf I,\check p};[0,T_i+\zeta]}^{2}\right),
\]
for some \(C, \lambda>0\), whenever \(Z|_{[0,T_i]}=Z^{[T_i]}\). Thus the radius of
the ball on which the next contraction is performed depends quadratically on the
previous controlled-path norm. In the regime \(3<p<4\), the corresponding
dependence is cubic. Consequently, iterating the local construction does not by
itself preclude the successive radii from diverging in finite time.

This is a limitation of the fixed-point argument, rather than a statement that
global well-posedness fails. A detailed discussion of this phenomenon, including
worked-out examples, can be found in Section~5 of \cite{deya2009rough}.

In the Young regime \(p\in[1,2)\), this obstruction disappears. In that case the
controlled-path structure contains no higher rough levels, and the estimates for
the solution map grow only linearly in the relevant path norm. The local
construction can therefore be iterated on a finite partition of \([0,T]\), with
the radii controlled by a Gronwall-type argument, yielding global well-posedness.
\end{remark}

\begin{proposition}\label{global_well_posedness_young}
Assume that $p \in [1,2)$ and  in addition that
\(g \in C_b^{2}(\mathbb  R^{d_Z},
\mathrm{Hom}(\mathbb{R}^{d_I},\mathbb  R^{d_Z}))\) and
\(b\in C_b^2(\mathbb  R^{d_Z},\mathbb  R^{d_Z})\). Then under Conditions \ref{condition_covariance_rho_var} and \ref{cameron_martin_condition}, equation
\eqref{proxy_eqtn} admits a unique global solution
\(
    Z \in \mathscr{D}_{\mathbf I,p_*}(\mathbb  R^{d_Z}).
\)

\end{proposition}
\begin{proof}
Since \(p<2\), the integral against \(I\) is a Young integral. Thus no higher
rough path levels are needed in the continuation argument.

For \(0\leq s<t\leq T\), set
\[
    \varpi_K(s,t)
    :=
    \sup_{\|h\|_{\infty;[0,T]}\leq 1}
    \|Kh\|_{p\text{-var};[s,t]},
    \qquad
    \varpi(s,t)
    :=
    \|I\|_{p\text{-var};[s,t]}+\varpi_K(s,t).
\]
By Lemma~\ref{regularity_Kf}, together with the continuity of the
\(p\)-variation control of \(I\), we have
\[
    \varpi(s,t)\to 0
    \qquad
    \text{uniformly as } |t-s|\downarrow0 .
\]

We first prove that the local existence time can be chosen independently of
the current value of the solution. Fix \(0\leq s<t\leq T\), an initial value
\(z\in\mathbb R^{d_Z}\), and define
\[
    \mathcal B_{s,t}(z)
    :=
    \left\{
        Z\in \mathcal C^{p}([s,t];\mathbb R^{d_Z})
        :
        Z_s=z,\,
        \|Z-z\|_{p;[s,t]}\leq 1
    \right\}.
\]
On this ball consider the map
\[
    \mathcal M(Z)_u
    :=
    z
    +
    \bigl(Kb(Z)\bigr)_{s,u}
    +
    \int_s^u g(Z_r)\,dI_r,
    \qquad u\in[s,t],
\]
where, if \(s>0\), the path \(Z\) is concatenated with the already constructed
solution on \([0,s]\).

By the Young estimate and the boundedness of \(b\), \(g\), and \(Dg\), there
exists a constant \(C\), depending only on \(p,b,g,K,T\), such that for every
\(Z\in\mathcal B_{s,t}(z)\),
\[
\begin{aligned}
    \|\mathcal M(Z)-z\|_{p\text{-var};[s,t]}
    &\leq
    \|Kb(Z)\|_{p\text{-var};[s,t]}
    +
    C
    \bigl(
        \|g\|_\infty
        +
        \|Dg\|_\infty
        \|Z\|_{p\text{-var};[s,t]}
    \bigr)
    \|I\|_{p\text{-var};[s,t]}
    \\
    &\leq
    C\,\omega(s,t).
\end{aligned}
\]
Hence, if \(\omega(s,t)\leq \delta\) with \(\delta>0\) sufficiently small,
then
\[
    \mathcal M(\mathcal B_{s,t}(z))
    \subseteq
    \mathcal B_{s,t}(z).
\]

Similarly, for \(Z,\widetilde Z\in\mathcal B_{s,t}(z)\), using the Lipschitz
boundedness of \(b\) and \(g\), the regularising estimate for \(K\), and again
the Young estimate, we obtain
\[
\begin{aligned}
    \|\mathcal M(Z)-\mathcal M(\widetilde Z)\|_{p;[s,t]}
    &\leq
    \|K(b(Z)-b(\widetilde Z))\|_{p;[s,t]}
    \\
    &\quad
    +
    C
    \|g(Z)-g(\widetilde Z)\|_{p;[s,t]}
    \|I\|_{p;[s,t]}
    \\
    &\leq
    C\,\omega(s,t)
    \|Z-\widetilde Z\|_{p;[s,t]} .
\end{aligned}
\]
Decreasing \(\delta\) if necessary, we may assume \(C\delta\leq 1/2\). Thus,
whenever \(\varpi(s,t)\leq\delta\), the map \(\mathcal M\) is a contraction on
\(\mathcal B_{s,t}(z)\). Therefore the equation has a unique solution on
\([s,t]\), and the constants involved do not depend on \(x\).

Since \(\varpi(s,t)\to0\) uniformly on shrinking intervals, there exists a
finite partition
\[
    0=t_0<t_1<\cdots<t_N=T
\]
such that
\[
    \varpi(t_i,t_{i+1})\leq\delta,
    \qquad i=0,\ldots,N-1.
\]
Starting from \(Z_{t_0}=z_0\), we solve successively on each interval
\([t_i,t_{i+1}]\), using \(Z_{t_i}\) as the initial condition for the next
step. Since the local existence radius and the contraction constants are
independent of \(Z_{t_i}\), this construction reaches \(T\) after finitely
many steps.

Moreover,
\[
    \|Z\|_{p;[0,T]}^p
    \leq
    \sum_{i=0}^{N-1}
    \|Z\|_{p;[t_i,t_{i+1}]}^p
    <\infty,
\]
so the concatenated path belongs to \(\mathcal C^{p}([0,T];\mathbb R^{d_Z})\).
Uniqueness follows by applying the same local contraction estimate on each
interval of the partition. Hence the local solution extends uniquely to the
whole interval \([0,T]\).
\end{proof}

However it is still possible to show that uniqueness holds even for the case $ p \in [2, 4)$.
\begin{proposition}\label{pathwise_uniqueness_proxy_equation}
Assume that the conditions of Proposition \ref{local_solution_sensor_equation} hold. Then
equation \eqref{proxy_eqtn} admits at most one solution on \([0,T]\) in
\(
    \mathscr{D}_{\mathbf I, p_*}(\mathbb R^{d_Z}).
\)
\end{proposition}

\begin{proof}
The case \(p\in(1,2]\) is obtained from the global well-posedness result above.
We give the argument in the case \(p\in[3,4)\). The case $p \in [2, 3)$ 
is obtained by the same argument, with fewer controlled-path levels.

Let \(Z\) and \(\widetilde Z\) be two solutions of \eqref{proxy_eqtn}.
From Proposition \ref{local_solution_sensor_equation}  \(Z\) and \(\widetilde Z\) agree on \([0,T_0]\). We show that there exists
\(\zeta>0\) such that they also agree on \([0,T_0+\zeta]\). Iterating this
argument then gives uniqueness on the whole interval \([0,T]\).

Let \(Z^{[T_0]}\) denote the common restriction of the two solutions to
\([0,T_0]\). Choose \(p_*<\check p<4\). For \(\zeta>0\), define the class of
extensions
\[
\mathcal B_{T_0,\zeta}
:=
\left\{
z\in \mathscr{D}_{\mathbf I,\check p}([0,T_0+\zeta];\mathbb R^{d_Z})
:
z=Z^{[T_0]}\ \text{on }[0,T_0],
\ 
\|z\|_{\mathscr{D}_{\mathbf I,\check p};[T_0,T_0+\zeta]}
\le M
\right\}.
\]
Here equality on \([0,T_0]\) is understood in the controlled-path sense,
i.e. including the controlled-path derivatives.

On this set define the map \(\mathcal M\) by
\[
\mathcal M(z)_t
=
\begin{cases}
\left(
Z^{[T_0]}_t,
g(Z^{[T_0], (0)}_t),
Dg(Z^{[T_0], (0)}_t)Z^{[T_0],(1)}_t
\right),
& t\in[0,T_0],
\\[0.4em]
\left(
Z_0+Kb(z)(t)+\displaystyle\int_0^t g(z^{(0)})_u\,d\mathbf I_u,
g(z_t^{(0)}),
Dg(z_t)z^{(1)}_t
\right),
& t\in(T_0,T_0+\zeta].
\end{cases}
\]
By taking \(M\) large enough, both \(Z\) and \(\widetilde Z\) belong to
\(\mathcal B_{T_0,\zeta}\) for \(\zeta\) sufficiently small.

We now show that \(\mathcal M\) is a contraction on
\(\mathcal B_{T_0,\zeta}\) for \(\zeta\) small enough. Since any two
elements of \(\mathcal B_{T_0,\zeta}\) agree on \([0,T_0]\), all their
controlled-path remainders agree on \([0,T_0]\). Therefore the variation of
their difference on \([0,T_0+\zeta]\) is the same as the variation of their
difference on \([T_0,T_0+\zeta]\).

Let \(z,\widetilde z\in\mathcal B_{T_0,\zeta}\). By the local estimates for
the \(Kb\)-term from Lemma \ref{regularity_Kf}, there exists a
modulus \(\omega_K\), with \(\omega_K(\zeta)\to0\) as \(\zeta\downarrow0\),
such that
\[
\|Kb(z)-Kb(\widetilde z)\|_{\mathscr{D}_{\mathbf I,\check p};
[T_0,T_0+\zeta]}
\le
\omega_K(\zeta)
\|z-\widetilde z\|_{\mathscr{D}_{\mathbf I,\check p};
[T_0,T_0+\zeta]}.
\]
Moreover, by the local Lipschitz estimates for the rough integral (see Proposition \ref{local_solution_sensor_equation}), there is
a constant \(C\), depending only on \(\check p,f,g,K,T\) and \(M\), but not
on \(T_0\) or \(\zeta\), such that
\[
\left\|
\int_0^\cdot g(z)_u\,d\mathbf I_u
-
\int_0^\cdot g(\widetilde z)_u\,d\mathbf I_u
\right\|_{\mathscr{D}_{\mathbf I,\check p};
[T_0,T_0+\zeta]}
\le
C\|\mathbf I\|_{\check p;[T_0,T_0+\zeta]}
\|z-\widetilde z\|_{\mathscr{D}_{\mathbf I,\check p};
[T_0,T_0+\zeta]}.
\]
Combining these estimates gives
\[
\|\mathcal M(z)-\mathcal M(\widetilde z)\|_{\mathscr{D}_{\mathbf I,\check p};
[T_0,T_0+\zeta]}
\le
\left(
\omega_K(\zeta)
+
C\|\mathbf I\|_{\check p;[T_0,T_0+\zeta]}
\right)
\|z-\widetilde z\|_{\mathscr{D}_{\mathbf I,\check p};
[T_0,T_0+\zeta]}.
\]
Since \(\mathbf I\) has finite \(\check p\)-variation, we have
\[
\|\mathbf I\|_{\check p;[T_0,T_0+\zeta]}\to0
\qquad
\text{as }\zeta\downarrow0.
\]
Thus, choosing \(\zeta>0\) sufficiently small, we obtain
\[
\omega_K(\zeta)
+
C\|\mathbf I\|_{\check p;[T_0,T_0+\zeta]}
<1.
\]
Hence \(\mathcal M\) is a contraction on \(\mathcal B_{T_0,\zeta}\).

Since both \(Z\) and \(\widetilde Z\) solve \eqref{proxy_eqtn}, they are
fixed points of \(\mathcal M\) on \([0,T_0+\zeta]\). By uniqueness of the
fixed point, \(Z=\widetilde Z\) on \([0,T_0+\zeta]\).

Starting from the local uniqueness interval and repeating the argument over
successive subintervals, we cover the whole interval \([0,T]\). This proves
uniqueness on \([0,T]\). Existence on \([0,T]\) follows by the same
continuation argument applied to the local solution.
\end{proof}
\begin{remark}
The assumption that the driving signal is a Volterra Gaussian process is used
only to identify the associated deterministic operator \(K\). The arguments in
this section are otherwise pathwise. Consequently, the results extend directly
to any deterministic kernel \(K\) satisfying the analytic estimates stated in
Lemma \ref{regularity_Kf}, provided that the driving signal is enhanced
to a rough path \(\mathbf I\).
\end{remark}
The Volterra Gaussian assumption will play a role again in the next section.
The preceding results provide a pathwise local theory, valid beyond the
Gaussian setting whenever the corresponding assumptions on \(K\) and
\(\mathbf I\) are satisfied. To obtain global existence when $p \in [2, 4)$, however, we shall use
the additional probabilistic structure coming from the Gaussianity of the
driving signal.

\subsection{Global existence of Gaussian RDEs with Volterra drift}
By the pathwise local well-posedness result established in the previous
section, equation \eqref{sensor_equation} admits a local solution in $\mathscr D_{\mathbf I, p_*}$ on some
interval \([0,T_0]\), with \(T_0\in(0,T]\), where $p_*$ is defined as above. Moreover, whenever a global solution
exists, it is unique path by path in the following sense.

\begin{definition}[Path-by-path uniqueness]\label{definition_path_by_path_uniqueness}
We say that equation \eqref{sensor_equation} satisfies path-by-path uniqueness
if there exists a set \(\Omega_1\) of full measure such that, for every \(\omega\in\Omega_1\), any two
solutions
\[
    Z(\omega),\widetilde Z(\omega)
    \in
    \mathscr D_{\mathbf I(\omega),p_*}(\mathbb R^{d_Z})
\]
of the deterministic equation driven by the realised rough path
\(\mathbf I(\omega)\), with the same initial condition \(z_0(\omega)\), must
coincide. That is,
\(
    Z(\omega)=\widetilde Z(\omega).
\)
\end{definition}
In what follows we consider the setup introduced in Section \ref{section_setup}, additionally, to keep the notation light, we will keep omitting the explicit dependence on
\(\omega\) of both the driving rough path and the corresponding solution. The
remainder of this section is devoted to proving global well-posedness for
equation \eqref{sensor_equation}.

Having established path-by-path uniqueness, our strategy is to formulate a
suitable analogue of the Yamada--Watanabe theorem for
\eqref{sensor_equation}. As in the classical setting, the aim is to show that
weak existence, combined with path-by-path uniqueness, implies the existence of
a \enquote{strong} solution.

More precisely, we shall prove that, under an appropriate notion of weak
solution, there exists a full-measure set \(\Omega_2 \subseteq\Omega\) and an
\((\mathcal F_t)\)-adapted process
\(
    Z\in \mathscr D_{\mathbf I,p_*}(\mathbb R^{d_Z})
\)
such that, for every \(\omega\in\Omega_2\), \(Z(\omega)\) solves
\eqref{sensor_equation} in the sense of Definition
\ref{operator_drift_solution} (see Theorem \ref{weaK_solution_existence}). Moreover, this solution is unique path by path,
and hence its law is uniquely determined. We shall refer to such a solution as a
strong solution of \eqref{sensor_equation}.

The next definition clarifies the notion of weak solution used below.
\begin{definition}[Weak solution]\label{definition_weak_solution}
A weak solution of equation \eqref{sensor_equation} consists of a filtered
probability space
\[
    (\Omega,\mathcal F,(\mathcal F_t)_{t\in[0,T]},\mathbb P),
\]
an \(\mathbb R^{d_X}\)-valued random variable \(x_0\), a
\((d_B+d_Y)\)-dimensional Volterra Gaussian process \(I\), and a process
\[
    Z\in \mathscr D_{\mathbf I,p_*}
    (\mathbb R^{d_Z}),
\]
such that the following properties hold:
\begin{enumerate}
    \item the filtration \((\mathcal F_t)_{t\in[0,T]}\) satisfies the usual
    conditions;

    \item \(I\) is adapted to \((\mathcal F_t)_{t\in[0,T]}\), has Volterra
    kernel \(K\), and admits a measurable rough path enhancement \(\mathbf I\);

    \item \(x_0\) is \(\mathcal F_0\)-measurable, has law \(\mu_0\), and is
    independent of \(I\);

    \item \(Z\) is adapted to \((\mathcal F_t)_{t\in[0,T]}\);

    \item for \(\mathbb P\)-almost every \(\omega\in\Omega\), the path
    \(Z(\omega)\) solves for every $t \in [0, T]$
    \[
        Z_t
        =
        z_0 + K b(Z)(t)
        +
        \int_0^t g(Z)_u\,d\mathbf I_u
    \]
    in the sense of Definition \ref{operator_drift_solution}, where
    \(z_0=(x_0,0)\).
\end{enumerate}
\end{definition}

Before turning to the weak formulation, we recall Condition
\ref{condition_filtration}. The purpose of this condition is to ensure that
the Volterra Gaussian process carries exactly the same information as the
Brownian motion appearing in its Volterra representation. Intuitively, the
Volterra property gives one inclusion, since the process is constructed from
the Brownian motion up to the same time. The condition below gives the reverse
inclusion: it allows one to reconstruct the Brownian increment \(V_t\) from the
Volterra process observed up to time \(t\).

This identification will be used repeatedly in the sequel since it allows us
to work with the filtration generated by the Volterra Gaussian process as if it
were the Brownian filtration. This is needed  in the
Yamada--Watanabe-type argument below.

As we mentioned above, since \(I\) generates the same filtration as the underlying Brownian motion
\(V\), the $(\mathcal F_t)_{t \in [0, T]}$ is equivalently the usual augmentation of the
filtration generated by \(x_0\) and \(V\).

We recall the notions of pathwise uniqueness and uniqueness in law that will be
used below. In the following definitions, uniqueness is stated at the level of the trace
\(Z^{(0)}\). This is sufficient for solutions in the sense of Definition
\ref{operator_drift_solution}, since the higher controlled-path components are
then prescribed recursively by
\[
    Z^{(k)} = g(Z)^{(k-1)}.
\]
Thus, once the trace of a solution is fixed, the corresponding
controlled-path structure is fixed as well. 

\begin{definition}[Pathwise uniqueness]\label{definition_pathwise_uniqueness}
We say that pathwise uniqueness holds for equation \eqref{sensor_equation} if
the following property is satisfied. Let
\[
    (Z,\mathbf I),\quad
    (\Omega,\mathcal F,(\mathcal F_t)_{t\in[0,T]},\mathbb P)
\]
and
\[
    (\widetilde Z,\mathbf I),\quad
    (\Omega,\mathcal F,(\widetilde{\mathcal F}_t)_{t\in[0,T]},\mathbb P)
\]
be two weak solutions of \eqref{sensor_equation} on the same probability space,
driven by the same enhanced Volterra Gaussian process \(\mathbf I\), and with
the same initial condition \(z_0\). Then
\[
    Z^{(0)}_t=\widetilde Z^{(0)}_t
    \qquad
    \text{for all }t\in[0,T],
    \quad \mathbb P\text{-a.s.}
\]
Equivalently, the solution paths \(Z^{(0)}\) and \(\widetilde Z^{(0)}\) are
indistinguishable.
\end{definition}

\begin{definition}[Uniqueness in law]\label{definition_uniqueness_in_law}
We say that uniqueness in law holds for equation \eqref{sensor_equation} if the
following property is satisfied. Let
\[
    (Z,\mathbf I),\quad
    (\Omega,\mathcal F,(\mathcal F_t)_{t\in[0,T]},\mathbb P)
\]
and
\[
    (\widetilde Z,\widetilde{\mathbf I}),\quad
    (\widetilde\Omega,\widetilde{\mathcal F},
    (\widetilde{\mathcal F}_t)_{t\in[0,T]},\widetilde{\mathbb P})
\]
be two weak solutions of \eqref{sensor_equation}, possibly defined on different
probability spaces, with the same initial distribution. Then the laws of the
solution paths for the trace agree:
\(
    \mathcal L_{\mathbb P}(Z^{(0)})
    =
    \mathcal L_{\widetilde{\mathbb P}}(\widetilde Z^{(0)})
\)
as probability measures on the path space
\(
    C([0,T],\mathbb R^{d_Z}).
\)
\end{definition}

It is immediate that path-by-path uniqueness implies pathwise uniqueness.
We are now ready to state one of the main results of this section
\begin{proposition}[Yamada--Watanabe for Volterra RDEs]\label{Yamada_Watanabe}
Assume that pathwise uniqueness holds for equation \eqref{sensor_equation}.
Then uniqueness in law holds.
\end{proposition}

\begin{proof}
The proof is inspired by the classical argument; see, for example, Proposition 5.3.20
in \cite{karatzas1991brownian}. 

Let
\(
    (Z,\mathbf I),
    (\Omega,\mathcal F,(\mathcal F_t)_{t\in[0,T]},\mathbb P)
\)
and
\(
    (\widetilde Z,\widetilde{\mathbf I}),
    (\widetilde\Omega,\widetilde{\mathcal F},
    (\widetilde{\mathcal F}_t)_{t\in[0,T]},\widetilde{\mathbb P})
\)
be two weak solutions of \eqref{sensor_equation} with the same initial
distribution. Here \(\mathbf I\) and \(\widetilde{\mathbf I}\) denote the
enhanced Volterra Gaussian driving signals. 
Set
\[
    L_t := Z^{(0)}_t- Z^{(0)}_0,
    \qquad
    \widetilde L_t :=
    \widetilde Z^{(0)}_t-\widetilde Z^{(0)}_0,
    \qquad 0\leq t\leq T.
\]
The two weak solutions can be represented through the triples
\(
    (Z^{(0)}_0,\mathbf I,L)\) and \(
    (\widetilde Z^{(0)}_0,\widetilde{\mathbf I},\widetilde L).
\)
Let
\[
    \Theta
    :=
    \mathbb R^{d_Z}
    \times
    \Omega G_p([0, T], \mathbb R ^I)
    \times
    C_0([0,T],\mathbb R^{d_Z}),
\]
endowed with its Borel \(\sigma\)-field. Define probability measures \(Q\) and
\(\widetilde Q\) on \(\Theta\) by
\(
    Q(A)
    :=
    \mathbb P\left[
        (Z^{(0)}_0,\mathbf I,L)\in A
    \right]
\)
and
\(
    \widetilde Q(A)
    :=
    \widetilde{\mathbb P}\left[
        (\widetilde Z^{(0)}_0,
        \widetilde{\mathbf I},\widetilde L)\in A
    \right].
\)
The first marginal of both \(Q\) and \(\widetilde Q\) is the common initial law
\(\mu_0\). The second marginal is the law of the enhancement of
the Volterra Gaussian driving signal to a rough path, which we denote by \(\mathbf P_I\).
Since the initial condition is independent of the driving signal in both weak
solutions, the joint law of the first two coordinates is
\(
    \lambda := \mu_0\otimes \mathbf P_I
\)
for both \(Q\) and \(\widetilde Q\).

Since the space of geometric rough path is Polish (see Proposition 5.38 in \cite{friz2010multidimensional}), by disintegration, there exist regular conditional probabilities
\(
    q(z,\mathbf i;d\ell),
    \) and  \(
    \widetilde q(z,\mathbf i;d\widetilde\ell),
\)
such that
\(
    Q(dz,d\mathbf i,d\ell)
    =
    \lambda(dz,d\mathbf i)\,q(z,\mathbf i;d\ell),
\)
and
\(
    \widetilde Q(dz,d\mathbf i,d\widetilde\ell)
    =
    \lambda(dz,d\mathbf i)\,
    \widetilde q(z,\mathbf i;d\widetilde\ell).
\)
Now define a probability measure \(\overline Q\) on the product space with
coordinates \((z,\mathbf i,\ell,\widetilde\ell)\) by
\[
    \overline Q(dz,d\mathbf i,d\ell,d\widetilde\ell)
    :=
    \lambda(dz,d\mathbf i)\,
    q(z,\mathbf i;d\ell)\,
    \widetilde q(z,\mathbf i;d\widetilde\ell).
\]
On this space define
\[
   \mathcal Z_t := z+\ell_t,
    \qquad
    \mathcal {\widetilde Z}_t := z+\widetilde\ell_t,
    \qquad
    \overline{\mathbf I}:=\mathbf i.
\]
By construction, \(\mathcal Z\) and \(\mathcal {\widetilde Z}\) have the same
initial condition and are driven by the same enhanced rough path
\(\overline{\mathbf I}\). Moreover, their respective laws coincide with those
of the traces of the original weak solutions. Therefore they can be made into weak
solutions of \eqref{sensor_equation} on a common probability space.

By pathwise uniqueness,
\[
    \mathcal Z_t=\mathcal {\widetilde Z}_t
    \qquad
    \text{for all }t\in[0,T],
    \quad \overline Q\text{-a.s.}
\]
Equivalently,
\(
    \ell=\widetilde\ell
    \)  \(
    \overline Q\text{-a.s.}
\)
It follows that
\(
    q(z,\mathbf i;\cdot)
    =
    \widetilde q(z,\mathbf i;\cdot)
\)
for \(\lambda\)-almost every \((x,\mathbf i)\). Hence
\(
    Q=\widetilde Q.
\)
In particular, the laws of
\(
    Z^{(0)}
    =
    Z^{(0)}_0+L
\)
and
\(
    \widetilde Z^{(0)}
    =
    \widetilde Z^{(0)}_0+\widetilde L
\)
coincide on \(C([0,T],\mathbb R^{d_Z})\). Thus uniqueness in law holds.
\end{proof}

\begin{remark}
For the filtering problem developed below, the essential well-posedness input is
existence of a weak solution together with uniqueness in law. Indeed, the
filter is determined by the joint law of the signal-observation process solving
the prescribed dynamics, and therefore does not depend on the particular
probability space on which this process is realised. The stronger pathwise
formulation that we derive below is not needed merely to identify the filter, but it is useful
for the subsequent, since it realises the
solution as a measurable non-anticipative functional of the initial condition
and the driving rough path.
\end{remark}

Having established that pathwise uniqueness implies uniqueness in law, we
continue with the Yamada--Watanabe argument by showing that the solution is, in
fact, a measurable functional of the initial condition and the driving rough
path. To this end, given a weak solution \((Z,I)\), we disintegrate the law of
\((Z_0,I,L)\), where \(L_t:=Z_t-Z_0\), with respect to \((Z_0,I)\). The first
result below shows that the resulting conditional law is adapted: events
depending on the solution only up to time \(t\) have conditional probabilities
which depend only on the initial condition and on the stopped driver. The
second result uses pathwise uniqueness to show that this conditional law is
almost surely a Dirac mass. This yields a measurable, non-anticipative map
\(\zeta\) such that \(Z=Z_0+\zeta(Z_0,I)\), and hence gives the desired strong
realisation of the solution.
Before proceeding with the statements, we introduce the notation \[
\mathscr B_t\bigl(C_0([0,T],\mathcal V)\bigr) := \sigma\bigl(\gamma \to \gamma_s:\,0\le s\le t\bigr). 
\]

\begin{proposition}[Adapted from Problem 5.3.21 in \cite{karatzas1991brownian}]
\label{Problem_5.3.21}
Let
\[
    (Z,\mathbf I),
    (\Omega,\mathcal F,(\mathcal F_t)_{t\in[0,T]},\mathbb P)
\]
be a weak solution of \eqref{sensor_equation}. Let $\mathbf P_I$ the law of $\mathbf{I}$ and \( Q\) be the  law of
\(
    (Z^{(0)}_0,\mathbf I,L),
\)
where
\(
    L_t := Z^{(0)}_t-Z^{(0)}_0.
\)
Let \(q(x,\mathbf i;\cdot)\) denote regular conditional probability
of the third coordinate, conditional on the first two coordinates.

Fix \(t\in[0,T]\). Then, for every
\(
    F\in
    \mathscr B_t\bigl(C_0([0,T],\mathbb R^{d_Z})\bigr),
\)
the map
\[
    (x,\mathbf i)\to q(x,\mathbf i;F),
\]
are \(\widehat{\mathscr B}_t\)-measurable, where
\(\widehat{\mathscr B}_t\) denotes the completion, with respect to
\(\mu_0(dx)\mathbf P_I(d\mathbf i)\), of
\[
    \mathscr B(\mathbb R^{d_Z})
    \otimes
    \mathscr B_t
    \left(
    \Omega G_{p}([0, T], \mathbb R^I)
    \right).
\]
\end{proposition}
\begin{proof}
Let \(V\) denote the Brownian motion underlying the Volterra Gaussian driver
\(I\). By Condition \ref{condition_filtration}, the completed filtration
generated by \(I\), equivalently by its enhancement \(\mathbf I\), coincides
with the completed filtration generated by \(V\). Therefore it is enough to
prove the corresponding measurability statement after pulling the kernel back
to the Brownian path space.

Let \(P_V\) be Wiener measure and let
\[
    \mathscr S_K:
    C_0([0,T],\mathbb R^{d_I})
    \to \Omega G_{p}([0, T], \mathbb R^I)
\]
denote the measurable map which sends the Brownian path to the canonical
Volterra Gaussian rough path enhancement. Thus
\(
    \mathbf P_I=(\mathscr S_K)_*P_I.
\)
Define
\[
    \mathcal Q(x,w;F):=q(x,\mathscr S_K(w);F).
\]
It suffices to show that \(\mathcal Q(\cdot,\cdot;F)\) is measurable with respect to
the completion of
\(
    \mathscr B(\mathbb R^{d_Z})
    \otimes
    \mathscr B_t\bigl(C_0([0,T],\mathbb R^{d_I})\bigr).
\) Let \(\varphi_t\) be the stopping map on Brownian paths,
\[
    (\varphi_t w)(s):=w_{s\wedge t},
    \qquad 0\le s\le T,
\]
and let \(\sigma_t\) be the map,
\[
    (\sigma_t w)(s):=w_{(t+s)\wedge T}-w_t,
    \qquad 0\le s\le T.
\]
Let \(\mathcal Q^t(x,w;F)\) be a regular conditional probability of \(L\), restricted
to events \(F\in \mathscr B_t(C_0([0,T],\mathbb R^{d_Z}))\), conditional
on \((Z^{(0)}_0,\varphi_t V)\). By construction,
\[
    (x,w)\to \mathcal Q^t(x,w;F)
\]
is measurable with respect to the completed
\[
    \mathscr B(\mathbb R^{d_Z})
    \otimes
    \mathscr B_t\bigl(C_0([0,T],\mathbb R^{d_I})\bigr).
\]
We now show that \(\mathcal Q=\mathcal Q^t\), \(\mu_0\otimes P_V\)-almost surely, on such
events \(F\).

It is enough to prove that, for every
\[
    G\in
    \mathscr B(\mathbb R^{d_Z})
    \otimes
    \mathscr B\bigl(C_0([0,T],\mathbb R^{d_I})\bigr),
\]
one has
\begin{equation}\label{goal_ex_5.3.21}
    \int_G \mathcal Q^t(x,w;F)\,\mu_0(dx)P_V(dw)
    =
    \mathbb P\bigl[(Z^{(0)}_0,V)\in G,\ L\in F\bigr].
\end{equation}
Indeed, since \(\mathcal Q\) is a regular conditional probability given
\((Z^{(0)}_0,V)\), the same identity with \(\mathcal Q\) in place of \(\mathcal Q^t\) holds by
definition. Hence \eqref{goal_ex_5.3.21} implies that \(\mathcal Q\) and \(\mathcal Q^t\) are
versions of the same conditional probability.

Fix \(F\in \mathscr B_t(C_0([0,T],\mathbb R^{d_Z}))\). Since the class of sets \(G\) for which
\eqref{goal_ex_5.3.21} holds is a Dynkin system, it is therefore enough to prove \eqref{goal_ex_5.3.21} on a \(\pi\)-system
generating
\[
    \mathscr B(\mathbb R^{d_Z})
    \otimes
    \mathscr B\bigl(C_0([0,T],\mathbb R^{d_I})\bigr).
\]
We take this \(\pi\)-system to consist of sets of the form
\[
    G
    =
    G_1
    \times
    \bigl(
        \varphi_t^{-1}G_2
        \cap
        \sigma_t^{-1}G_3
    \bigr),
\]
where
\[
    G_1\in\mathscr B(\mathbb R^{d_Z}),
    \qquad
    G_2,G_3\in
    \mathscr B\bigl(C_0([0,T],\mathbb R^{d_I})\bigr).
\]
For such \(G\), we have
\[
\begin{aligned}
&\int_G \mathcal Q^t(x,w;F)\,\mu_0(dx)P_V(dw)
\\
&\quad =
\mathbb E^{P_V}
\left[
    \mathbbm 1_{\varphi_t^{-1}G_2}(V)
    \mathbbm 1_{\sigma_t^{-1}G_3}(V)
    \int_{G_1} \mathcal Q^t(x,V;F)\,\mu_0(dx)
\right].
\end{aligned}
\]
The random variable
\[
    \mathbbm 1_{\varphi_t^{-1}G_2}(V)
    \int_{G_1} \mathcal Q^t(x,V;F)\,\mu_0(dx)
\]
is measurable with respect to the Brownian filtration up to time \(t\), while
\(\mathbbm 1_{\sigma_t^{-1}G_3}(V)\) depends only on the Brownian increments
after time \(t\). Hence, by independence of Brownian increments,
\[
\begin{aligned}
&\int_G \mathcal Q^t(x,w;F)\,\mu_0(dx)P_V(dw)
\\
&\quad =
P_V(\sigma_t^{-1}G_3)
\int_{G_1\times\varphi_t^{-1}G_2}
    \mathcal Q^t(x,w;F)\,\mu_0(dx)P_V(dw).
\end{aligned}
\]
By the definition of \(R^t\), the last integral is
\[
    \mathbb P
    \bigl[
        Z^{(0)}_0\in G_1,\,
        \varphi_t V\in G_2,\,
        L\in F
    \bigr].
\]
Since \(F\in\mathscr B_t(C_0([0,T],\mathbb R^{d_Z}))\) and \(Z\) is
adapted, the event
\[
    \{Z^{(0)}_0\in G_1,\ \varphi_t V\in G_2,\ L\in F\}
\]
belongs to \(\mathcal F_t\). The event
\(
    \{\sigma_t V\in G_3\}
\)
depends only on the future Brownian increments and is therefore independent of
\(\mathcal F_t\). Consequently,
\[
\begin{aligned}
& P_V(\sigma_t^{-1}G_3)
\mathbb P
\bigl[
    Z^{(0)}_0\in G_1,\,
    \varphi_t V\in G_2,\,
    L\in F
\bigr]
\\
&\quad =
\mathbb P
\bigl[
    Z^{(0)}_0\in G_1,\,
    \varphi_t V\in G_2,\,
    \sigma_t V\in G_3,\,
    L\in F
\bigr]
\\
&\quad =
\mathbb P
\bigl[
    (Z^{(0)}_0,V)\in G,\,
    L\in F
\bigr].
\end{aligned}
\]
Thus \eqref{goal_ex_5.3.21} holds for the generating \(\pi\)-system, and hence
for every \(G\) by the Dynkin-system theorem. Therefore
\[
    \mathcal Q(x,w;F)=\mathcal Q^t(x,w;F)
\]
for \(\mu_0\otimes P_V\)-almost every \((x,w)\).

It follows that \(\mathcal Q(\cdot,\cdot;F)\) is measurable with respect to the
completed past Brownian \(\sigma\)-field. Using again Condition
\ref{condition_filtration}, this completed \(\sigma\)-field coincides with the
completed \(\sigma\)-field generated by the stopped enhanced Volterra rough
path. Hence
\[
    (x,\mathbf i)\to q(x,\mathbf i;F)
\]
is \(\widehat{\mathscr B}_t\)-measurable.
\end{proof}

\begin{proposition}[Adapted from Problem 5.3.22 in \cite{karatzas1991brownian}]
\label{exercise_3_22_ks}
In the setting of Propositions \ref{Yamada_Watanabe} and
\ref{Problem_5.3.21}, there exists a measurable map
\[ \varsigma
    :
    \mathbb R^{d_Z}
    \times
    \Omega G_{p}([0, T], \mathbb{R}^{d_I})
    \to
    C_0([0,T],\mathbb R^{d_Z})
\]
such that, for \(\mu_0\otimes  \mathbf P_I\)-almost every \((z,w)\),
\begin{equation}\label{identity_Q}
    q(z,w;\{\varsigma(z,w)\})
    =
    \tilde q(x,w;\{\varsigma(x,w)\})
    =
    1.
\end{equation}
Moreover, \(\varsigma\) is
\[
    \mathscr B(\mathbb R^{d_Z})
    \otimes
    \mathscr B\left(
    \Omega G_{p}([0, T], \mathbb{R}^{d_I})\right)
\]
measurable. In addition, for every \(t\in[0,T]\), the stopped map
\[
    (z,w)\to \varsigma (z,w)_{\cdot\wedge t}
\]
is
\(
    \widehat{\mathscr B}_t
    \big/
    \mathscr B_t
    \left(C_0([0,T],\mathbb R^{d_Z})\right)
\)
measurable.

Finally, if \(\overline{Q}\) denotes the probability measure on
\[
    \mathbb R^{d_Z}
    \times
    \Omega G_{p}([0, T], \mathbb{R}^{d_I})
    \times
    C_0([0,T],\mathbb R^{d_Z})
    \times
    C_0([0,T],\mathbb R^{d_Z})
\]
defined by
\[
    \overline{Q}(dz,dw,d\ell,d \tilde \ell)
    :=
    \mu_0(dz) \mathbf P_{I}(d\mathbf i)q(x,w;d \tilde \ell)\tilde q(x,w; d \tilde \ell),
\]
then
\begin{equation}\label{pathwise_uniqueness_k}
    \overline{Q}
    \left[
        \ell=\tilde \ell=k(z,w)
    \right]
    =
    1.
\end{equation}
\end{proposition}
\begin{proof}
Set
\[
    \lambda(dz,dw):=\mu_0(dz)\mathbf P_I(dw).
\]
We first construct the product coupling of the two conditional laws. Define
\(\overline{Q}\) by
\[
    \overline{Q}(dz,dw,d\ell,d\tilde \ell)
    :=
    \lambda(dz,dw)q(z,w;d\ell)\tilde q(z,w;d\tilde \ell).
\]
Under this probability measure, by Proposition \ref{Yamada_Watanabe}, the two coordinate processes
\(
    z+\ell,
    \) and  \(
    z+ \tilde \ell
\)
have the laws of the traces of two weak solutions with the same initial
condition \(z\) and the same enhanced driving rough path \(w\). Therefore, by
pathwise uniqueness,
\[
    \overline{Q}[\ell=\tilde \ell]=1.
\]
Equivalently, if
\(
    \Delta
    :=
    \{(\ell,\tilde \ell):\ell=\tilde \ell\},
\)
then
\[
    \int
    \bigl(q(z,w;\cdot)\otimes \tilde q(z,w;\cdot)\bigr)(\Delta)
    \,\lambda(dz,dw)
    =
    1.
\]
Hence there exists a \(\lambda\)-null set \(N\) such that, for every
\((z,w)\notin N\),
\begin{equation}\label{diagonal_full_measure}
    \bigl(q(z,w;\cdot)\otimes \tilde q(z,w;\cdot)\bigr)(\Delta)=1.
\end{equation}

We now show that \eqref{diagonal_full_measure} forces both conditional laws to
be the same Dirac mass. Fix \((z,w)\notin N\), and write
\[
    \mathscr q:= q(z,w;\cdot),
    \qquad
    \tilde{\mathscr q}:= \tilde q(z,w;\cdot).
\]
Then
\[
    1
    =
    (\mathscr q\otimes \tilde{\mathscr q})(\Delta)
    =
    \int \mathscr q(\{\ell\})\,\tilde{\mathscr q}(d\ell).
\]
Since \(0\le \mathscr q(\{\ell\})\le 1\), this implies that
\(\mathscr q(\{\ell\})=1\) for \(\tilde{\mathscr q}\)-almost every \(\ell\). Thus there exists a point,
which we denote by \(\varsigma (z,w)\), such that
\(
    \mathscr q(\{\varsigma(z,w)\})=1.
\)
Substituting this into \eqref{diagonal_full_measure} gives
\(
    \tilde{\mathscr q}(\{\varsigma(z,w)\})=1.
\)
Therefore
\[
    q(z,w;\{\varsigma(z,w)\})
    =
    \tilde q(z,w;\{\varsigma(z,w)\})
    =
    1
\]
for every \((z,w)\notin N\). On \(N\), define \(\varsigma (z,w)\) arbitrarily, for
instance \(\varsigma (z,w)\equiv 0\). This proves \eqref{identity_Q}.

It remains to check measurability. 
For every \(A\in\mathscr B(C_0([0,T],\mathbb R^{d_Z}))\), outside the null set \(N\) we have
\[
    \varsigma (z,w)\in A
    \quad\Longleftrightarrow\quad
    q(z,w;A)=1.
\]
Since \((z,w)\to q(z,w;A)\) is measurable by the definition of regular
conditional probability, it follows that \(\varsigma\) is measurable with respect to
\(
    \mathscr B(\mathbb R^{d_Z})
    \otimes
    \mathscr B\left(
    \Omega G_{p}([0, T], \mathbb R^I)
    \right),
\)
up to modification on a \(\lambda\)-null set. After modifying \(\varsigma \) on \(N\),
we obtain a measurable version.

We now prove the non-anticipative measurability. Fix \(t\in[0,T]\), and let
\(
    A\in
    \mathscr B_t\left(C_0([0,T],\mathbb R^{d_Z})\right).
\)
Again, outside \(N\),
\[
    \varsigma (z,w)\in A
    \quad\Longleftrightarrow\quad
    q(z,w;A)=1.
\]
By Proposition \ref{Problem_5.3.21}, the map
\[
    (x,w)\to q(z,w;A)
\]
is \(\widehat{\mathscr B}_t\)-measurable. Hence
\[
    \{(z,w):\varsigma(z,w)\in A\}\in\widehat{\mathscr B}_t.
\]
This proves that the stopped map
\[
    (z,w)\to \varsigma(z,w)_{\cdot\wedge t}
\]
is
\(
    \widehat{\mathscr B}_t
    \big/
    \mathscr B_t\left(C_0([0,T],\mathbb R^{d_Z})\right)
\)
measurable.

Finally, from \eqref{identity_Q}, for \(\lambda\)-almost every \((z,w)\),
\[
    q(z,w;\{\varsigma (z,w)\})\tilde q(z,w;\{\varsigma (z,w)\})=1.
\]
Integrating this identity with respect to \(\lambda\), we obtain
\(
    \overline{Q}
    \left[
        \ell=\tilde \ell=\varsigma (z,w)
    \right]
    =
    1,
\)
which is \eqref{pathwise_uniqueness_k}.
\end{proof}
In order to establish weak existence, we follow the classical approach by applying Girsanov’s theorem. A version of this theorem, applicable to fractional Brownian motion, is presented in Theorem 4.9 in \cite{decreusefond1999stochastic}. However, the proof extends with minimal modifications to Volterra Gaussian processes.  

In particular, for a shifted \( d_I \)-dimensional $\mathbb{P}$ Volterra Gaussian process of the form  
\[
\tilde{I}^k_t = I^k_t + \int_0^t K(t,s) b_k(s) ds \quad k = 1, \dots, d_I,
\]  
where the integral term is an adapted process in \( L^2([0, T], \mathbb{R}^d) \), an equivalent probability measure \(\mathbb{Q}\) can be introduced via the density process  
\[
\tilde{\Lambda}_t = \exp\Bigg(-\sum_{k=1}^{d}\int_0^t b_k(s) dV^{k}_s - \frac{1}{2} \sum_{k=1}^{d} \int_0^t (b_k(s))^2 ds\Bigg).
\]  
If this process is a true martingale, meaning that \(\mathbb{E}[\tilde{\Lambda}_T] = 1\), then \(\mathbb{Q}\) is well-defined and ensures that \(\tilde{I}_t\) remains a Volterra Gaussian process with the same kernel as \(I\) and is adapted to $ (\mathcal F_t)_{t \in [0, T]}$. 

We will now show that there exists a weak solution to \eqref{sensor_equation}.
Let \(x_0\) be as before and consider the probability space 
\(
    (\Omega, \mathcal{F}, (\tilde{\mathcal{F}}_t)_{t \in [0, T]}, \tilde{\mathbb P})
\) which supports a Volterra Gaussian process $(\tilde{B}, \tilde U )$ with same kernel $K$ as $(B, U)$. The filtration $(\tilde{\mathcal{F}}_t)_{t \in [0, T]}$ is constructed 
following the same procedure as in Section \ref{section_setup}.

On this space, consider the RDE driven by $\tilde{\mathbf I} = \boldchar{(}\tilde{\mathbf B}, \tilde{\mathbf U}\boldchar{)}$  
\begin{equation}\label{weak_sol_equation}
\begin{pmatrix}
\tilde{X}_t \\ \tilde{Y}_t 
\end{pmatrix} = \begin{pmatrix} x_0 \\ 0 \end{pmatrix} + \int_0^ t\begin{pmatrix} \sigma(\tilde{X}, \tilde{Y})_u & 0 \\0 & \mathbbm{1} \end{pmatrix} d\tilde{\mathbf I}_u.
\end{equation}  

For notational convenience, we introduce the process \(\tilde{Z}_t = (\tilde{X}_t, \tilde{Y}_t)\), which satisfies  
\[
\tilde{Z}_t = z_0  + \int_0^t  g(\tilde{Z})_u d \mathbf{\tilde I}_u.
\]  
As recalled in the introduction, this RDE admits a unique global solution in the space $\mathscr{D}_{\tilde I, p}(\mathbb R^{d_Z})$.  

We will show how this formulation can be interpreted as a weak solution to equation \eqref{sensor_equation} by introducing a change of probability measure. Specifically, consider $W^U$ to be the Brownian motion associated to $U$ (see Section \ref{section_preliminaries}) and the process  
\[
\Lambda^K_t = \exp\Bigg(\sum_{j =1}^{ d_Y}\int_0^t f_j(\tilde{X}, \tilde{Y})_u dW^{U, k}_u - \frac{1}{2} \sum_{j= 1}^{d_Y} \int_0^t (f_j(\tilde{X}_u, \tilde{Y}_u))^2 du\Bigg),
\]  
which defines the new probability measure \(\mathbb{P}_K\) under which the process \[
(\tilde B, \overline{U}) := \left(
\widetilde B,\ldots,\widetilde B^{d_B},
\tilde U^1 - Kf_1(\widetilde Z),\ldots,
\tilde U^{d_Y} - Kf_{d_Y}(\widetilde Z)
\right)\] is a centered Volterra Gaussian process with kernel $K$.

With an appropriate adaptation of the Yamada-Watanabe theorem, together with the fact that $(\tilde{Z}, \boldchar{(}\tilde{\mathbf B}, \overline{\mathbf U}\boldchar{)}), (\Omega, \mathcal{F}, \mathbb{P}_K), (\tilde{\mathcal{F}}_t)_{t \in [0, T]}$ is a weak solution, we conclude that a strong solution to equation \eqref{sensor_equation} exists.

\begin{proof}[Proof of Theorem \ref{weaK_solution_existence}]
We prove the claim only for \(p\in[3,4)\); the case \(p\in[2,3)\) follows from
the same argument with the obvious simplifications.

Let \(\widetilde Z\) be the solution of \eqref{weak_sol_equation}, and fix
\(\widetilde p\) as in Proposition \ref{local_solution_sensor_equation}. Under
\(\widetilde{\mathbb P}\), the first level of
\(\boldchar{(}\widetilde{\mathbf B},\widetilde{\mathbf U}\boldchar{)}\) is given by
\begin{equation}\label{first_lvl_B_U_tilde}
\boldchar{(}\widetilde{\mathbf B},\widetilde{\mathbf U}\boldchar{)}^{(1)}
=
\left(
\widetilde B^1,\ldots,\widetilde B^{d_B},
\overline U^1+Kf_1(\widetilde Z),\ldots,
\overline U^{d_Y}+Kf_{d_Y}(\widetilde Z)
\right).
\end{equation}
By construction, this process has the required Volterra Gaussian law with kernel
\(K\).

Define \(\boldchar{(}\widetilde{\mathbf B},\overline{\mathbf U}\boldchar{)}\) to be the geometric
rough path constructed above the \(\mathbb P_K\)-Volterra Gaussian process
\[
\left(
\widetilde B^1,\ldots,\widetilde B^{d_B},
\overline U^1,\ldots,\overline U^{d_Y}
\right).
\]
By Lemma \ref{regularity_Kf}, we have
\(
    Kf(\widetilde Z)\in \mathcal C^q([0,T],\mathbb R^{d_Y})
\)
for some \(q\) such that
\(
    1 / q+1 /p_* >1.
\)
Thus the Young integrals involving \(Kf(\widetilde Z)\) and
\(\boldchar{(}\widetilde{\mathbf B},\overline{\mathbf U}\boldchar{)}\) are well-defined.

For notational convenience, set
\[
    A:=Kf(\widetilde Z).
\]
Then the higher levels of \(\boldchar{(}\widetilde{\mathbf B},\widetilde{\mathbf U}\boldchar{)}\) can
be written as
\begin{equation}\label{second_level_B_U}
\begin{aligned}
    \boldchar{(}\widetilde{\mathbf B},\widetilde{\mathbf U}\boldchar{)}^{(2)}_{s,t}
    &=
    \boldchar{(}\widetilde{\mathbf B},\overline{\mathbf U}\boldchar{)}^{(2)}_{s,t}
    +
    \int_s^t
    \boldchar{(}\widetilde{\mathbf B},\overline{\mathbf U}\boldchar{)}^{(1)}_{s,u}
    \otimes dA_u
    \\
    &\quad
    +
    \int_s^t
    A_{s,u}\otimes
    d\boldchar{(}\widetilde{\mathbf B},\overline{\mathbf U}\boldchar{)}^{(1)}_u
    +
    \int_s^t A_{s,u}\otimes dA_u
    \\
    &=
    \boldchar{(}\widetilde{\mathbf B},\overline{\mathbf U}\boldchar{)}^{(2)}_{s,t}
    +
    T^1_{s,t},
    \\
    \boldchar{(}\widetilde{\mathbf B},\widetilde{\mathbf U}\boldchar{)}^{(3)}_{s,t}
    &=
    \boldchar{(}\widetilde{\mathbf B},\overline{\mathbf U}\boldchar{)}^{(3)}_{s,t}
    +
    \int_s^t
    \boldchar{(}\widetilde{\mathbf B},\overline{\mathbf U}\boldchar{)}^{(2)}_{s,u}
    \otimes dA_u
    \\
    &\quad
    +
    \int_s^t
    T^1_{s,u}\otimes
    d\boldchar{(}\widetilde{\mathbf B},\overline{\mathbf U}\boldchar{)}^{(1)}_u
    +
    \int_s^t
    T^1_{s,u}\otimes dA_u
    \\
    &=
    \boldchar{(}\widetilde{\mathbf B},\overline{\mathbf U}\boldchar{)}^{(3)}_{s,t}
    +
    T^2_{s,t}.
\end{aligned}
\end{equation}
We now show that \(\widetilde Z\) is controlled by
\(\boldchar{(}\widetilde{\mathbf B},\overline{\mathbf U}\boldchar{)}\). Let \(\widetilde R\) denote
the controlled-path remainder of \(\widetilde Z\) as an element of
\(
    \mathscr D_{\boldchar{(}\widetilde{\mathbf B},\widetilde{\mathbf U}\boldchar{)},p_*}
    (\mathbb R^{d_Z}),
\)
and let \(\overline R\) denote the corresponding remainder with respect to
\(
    \boldchar{(}\widetilde{\mathbf B},\overline{\mathbf U}\boldchar{)}.
\)
Then
\begin{align*}
\overline R^{(0)}_{s,t}
&=
\widetilde Z^{(0)}_{s,t}
-
\widetilde Z^{(1)}_s
\boldchar{(}\widetilde{\mathbf B},\overline{\mathbf U}\boldchar{)}^{(1)}_{s,t}
-
\widetilde Z^{(2)}_s
\boldchar{(}\widetilde{\mathbf B},\overline{\mathbf U}\boldchar{)}^{(2)}_{s,t}
\\
&=
\widetilde Z^{(0)}_{s,t}
-
\widetilde Z^{(1)}_s
\boldchar{(}\widetilde{\mathbf B},\widetilde{\mathbf U}\boldchar{)}^{(1)}_{s,t}
-
\widetilde Z^{(2)}_s
\boldchar{(}\widetilde{\mathbf B},\widetilde{\mathbf U}\boldchar{)}^{(2)}_{s,t}
\\
&\quad
+
\widetilde Z^{(1)}_s
\begin{pmatrix}
0\\
A_{s,t}
\end{pmatrix}
+
\widetilde Z^{(2)}_s T^1_{s,t}
\\
&=
\widetilde R^{(0)}_{s,t}
+ \widetilde Z^{(1)}_s
\begin{pmatrix}
0\\
A_{s,t}
\end{pmatrix}
+
\widetilde Z^{(2)}_sT^1_{s,t}.
\end{align*}
The Young integral estimates for the terms defining \(T^1\), together with the regularity
of \(Kf(\widetilde Z)\), imply that the last two terms have the required
\(\mathcal C^{p_*/\kappa}\)-regularity. Hence
\(
    \overline R^{(0)}\in \mathcal C^{p_*/\kappa}.
\)
Similarly,
\begin{align*}
\overline R^{(1)}_{s,t}
&=
\widetilde Z^{(1)}_{s,t}
-
\widetilde Z^{(2)}_s
\boldchar{(}\widetilde{\mathbf B},\overline{\mathbf U}\boldchar{)}^{(1)}_{s,t}
\\
&=
\widetilde Z^{(1)}_{s,t}
-
\widetilde Z^{(2)}_s
\boldchar{(}\widetilde{\mathbf B},\widetilde{\mathbf U}\boldchar{)}^{(1)}_{s,t}
+
\widetilde Z^{(2)}_s
\begin{pmatrix}
0\\
A_{s,t}
\end{pmatrix}
\\
&=
\widetilde R^{(1)}_{s,t}
+
\tilde Z^{(2)}_s\begin{pmatrix}
0\\
A_{s,t}
\end{pmatrix}.
\end{align*}
Therefore
\[
    \overline R^{(1)}
    \in
    \mathcal C^{p_*/(\kappa-1)}.
\]
Finally,
\[
\overline R^{(2)}_{s,t}
=
\widetilde R^{(2)}_{s,t}
=
D\sigma(\widetilde Z_t)\sigma(\widetilde Z_t)
-
D\sigma(\widetilde Z_s)\sigma(\widetilde Z_s)
\in
\mathcal C^{p_*/(\kappa-2)}.
\]
Thus
\[
    \widetilde Z
    \in
    \mathscr D_{\boldchar{(}\widetilde{\mathbf B},\overline{\mathbf U}\boldchar{)},p_*}
    (\mathbb R^{d_Z}).
\]
It remains to identify the equation solved by \(\widetilde Z\) when the driver is
\(\boldchar{(}\widetilde{\mathbf B},\overline{\mathbf U}\boldchar{)}\). Since
\[
    \widetilde Z_t
    =
    z_0
    +
    \int_0^t
    g(\widetilde Z)_u\,
    d\boldchar{(}\widetilde{\mathbf B},\widetilde{\mathbf U}\boldchar{)}_u,
\]
the definition of the rough integral, together with
\eqref{first_lvl_B_U_tilde} and \eqref{second_level_B_U}, yields
\begin{align*}
\widetilde Z_t
&=
z_0
+
\lim_{|\mathcal P|\to0}
\sum_{i}
\sum_{j=0}^{\kappa-1}
g(\widetilde Z)^{(j)}_{u_i}
\boldchar{(}\widetilde{\mathbf B},\overline{\mathbf U}\boldchar{)}^{(j+1)}_{u_i,u_{i + 1}}
\\
&\quad
+
\lim_{|\mathcal P|\to0}
\sum_{i}
\begin{pmatrix}
0\\
A_{u_i,u_{i+1}}
\end{pmatrix}
\\
&\quad
+
\lim_{|\mathcal P|\to0}
\sum_{i}
\left[
g(Z)^{(1)}_{u_i}T^1_{u_i,u_{i +1}}
+
g(Z)^{(2)}_{u_i}T^2_{u_i,u_{i +1}}\right].
\end{align*}
The first limit is precisely
\[
    \int_0^t
    g(\widetilde Z)_u\,
    d\boldchar{(}\widetilde{\mathbf B},\overline{\mathbf U}\boldchar{)}_u.
\]
Moreover, the terms involving \(T^1\) and \(T^2\) vanish as the mesh of the
partition tends to zero. Indeed, by the Young--Loeve estimates and the
boundedness of the coefficients appearing in
\(g(\widetilde Z)^{(1)}\) and \(g(\widetilde Z)^{(2)}\),
\[
\sum_{i}
\left|
g(\widetilde Z)^{(i)}_{u_i}T^i_{u_i,u_{i +1}}
\right|
\to 0,
\qquad i=1,2,
\]
as \(|\mathcal P|\to0\). Hence
\[
    \widetilde Z_t
    =
    z_0
    +
    Kf(\widetilde Z)(t)
    +
    \int_0^t
    g(\widetilde Z)_u\,
    d\boldchar{(}\widetilde{\mathbf B},\overline{\mathbf U}\boldchar{)}_u.
\]
Thus
\[
\left(
\widetilde Z,
\boldchar{(}\widetilde{\mathbf B},\overline{\mathbf U}\boldchar{)}\right), \left(
\widetilde\Omega,
\widetilde{\mathcal F},
(\widetilde{\mathcal F}_t)_{t\in[0,T]},
\widetilde{\mathbb P}
\right)
\]
is a weak solution of \eqref{sensor_equation}.
By Proposition \ref{pathwise_uniqueness_proxy_equation}, pathwise uniqueness
holds. Therefore Proposition \ref{Yamada_Watanabe} yields strong existence.
In particular, there exists a measurable map
\[
\varsigma:
\mathbb R^{d_Z}
\times
\Omega G_p([0, T], \mathbb R ^I)
\to
C([0,T],\mathbb R^{d_Z})
\]
such that \(\varsigma\) is
\(
\mathscr B(\mathbb R^{d_Z})
\otimes
\mathscr B\left(
\Omega G_p([0, T], \mathbb{R}^{d_I})
\right)
/
\mathscr B\left(
C([0,T],\mathbb R^{d_Z})
\right)
\)
measurable, and, for every \(t\in[0,T]\), its restriction is
\(
\widehat{\mathscr B}_t
/
\mathscr B_t\left(
C([0,T],\mathbb R^{d_Z})
\right)
\text{-measurable}.
\)
Moreover,
\[
    \widetilde Z^{(0)}_\cdot
    =
    \varsigma \left(
        z_0,
        \boldchar{(}\widetilde{\mathbf B},\overline{\mathbf U}\boldchar{)}
    \right)_\cdot
    \qquad
    \widetilde{\mathbb P}\text{-a.s.}
\]
Let \(\Gamma\) be the set of pairs \((z,\mathbf i)\) such that, setting
\(
    H_\cdot := \varsigma(z,\mathbf i)_\cdot,
\)
and equipping \(H\) with the controlled-path structure induced by the equation,
namely
\[
    H^{(1)} := g(H),
    \qquad
    H^{(2)} := Dg(H)g(H)
\]
in the case \(p\in[3,4)\), the controlled path
\[
    \mathbf H := (H,H^{(1)},H^{(2)})
\]
belongs to \(\mathscr D_{\mathbf i ,p_*}(\mathbb R^{d_Z})\) and satisfies
the fixed-point identity
\[
    H_t
    =
    z
    +
    Kf(H)(t)
    +
    \int_0^t g(H)_u\,d\mathbf i_u,
    \qquad t\in[0,T].
\]
By the preceding argument,
\[
    \left(
        Z_0,
        \boldchar{(}\widetilde{\mathbf B},\overline{\mathbf U}\boldchar{)}
    \right)
    \in \Gamma,
    \qquad
    \widetilde{\mathbb P}\text{-a.s.}
\]
Moreover, by construction,
\[
    \mathcal L_{\widetilde{\mathbb P}}
    \left(
        z_0,
        \boldchar{(}\widetilde{\mathbf B},\overline{\mathbf U}\boldchar{)}
    \right)
    =
    \mathcal L_{\mathbb P}
    \left(
        z_0,
        \boldchar{(}\mathbf B,\mathbf U\boldchar{)}
    \right).
\]
Therefore, $\mathbb P\text{-a.s.}$, 
\(
    \left(
        Z_0,
        \boldchar{(}\mathbf B,\mathbf U\boldchar{)}
    \right)
    \in \Gamma.
\)
Returning to the original probability space
\((\Omega,\mathcal F,\mathbb P)\), define
\[
    Z^{(0)}_\cdot
    :=
    \varsigma\left(
        Z_0,
        \boldchar{(}\mathbf B,\mathbf U\boldchar{)}
    \right)_\cdot .
\]
Then equip \(Z\) with the controlled-path structure
\[
    Z^{(1)} := g(Z),
    \qquad
    Z^{(2)} := Dg(Z)g(Z).
\]
By the previous inclusion, the controlled path
\(
    Z := (Z^{(0)},Z^{(1)},Z^{(2)})
\)
belongs to
\(
    \mathscr D_{\mathbf I,p_*}
    (\mathbb R^{d_Z})
\)
and satisfies \eqref{sensor_equation}. Since \(\varsigma\) is non-anticipative, \(Z^{(0)}\),
and hence the controlled path \(Z\), is adapted to the filtration
generated by \(Z_0\) and \(\boldchar{(}\mathbf B,\mathbf U\boldchar{)}\). Therefore \(Z\) is a
strong solution. This proves the claim.
\end{proof}

\section{Construction of the enlarged Volterra Gaussian process}\label{app:enlarged_VGP}
In this section, we construct the enhancement to a geometric rough path of the $\mathbb{R}^{2d}$ valued Gaussian process $(I, V)$ introduced in Section \ref{section_preliminaries}. As above, we assume that $I$ satisfies Conditions \ref{condition_covariance_rho_var}, \ref{condition_filtration} and has sample paths of finite $p$-variation. 
We first recall that the Gaussian rough path construction of
\cite{friz2010differential} yields a canonical geometric rough path lift
of \(I\) and $V$. Therefore we only have to take care of the mixed terms in the enhancement.\newline
We recall that, in order to show that \(\hat{\mathbf I}\) defines a geometric rough path, it is enough to verify the following three properties (see Section 2 in  \cite{cass2022combinatorial}). Let $\kappa := \floor{p}$, for every \(i=1,\ldots,\kappa\) and every \(0\leq s\leq t\leq T\),
\[
    \hat{\mathbf I}^{(i)}\in
    \mathcal C^{p/i}\bigl(\Delta_{[0,T]},\mathbb R^{\otimes i}\bigr),
\]
and, for every \(0\leq s\leq u\leq t\leq T\),
\[
    \operatorname{proj}_{\leq \kappa}
    \bigl(
        \hat{\mathbf I}_{s,u}\otimes \hat{\mathbf I}_{u,t}
    \bigr)
    =
    \hat{\mathbf I}_{s,t}.
\]
Finally, for every \(0\leq s\leq t\leq T\), \(\hat{\mathbf I}\) satisfies the shuffle identity
\[
    \Delta_{\shuffle}\hat{\mathbf I}_{s,t}
    =
    \hat{\mathbf I}_{s,t}\otimes \hat{\mathbf I}_{s,t}.
\]
In this construction, we restrict ourselves to the case \(\kappa=3\). The remaining cases follow from the same argument, with simpler calculations.

Let \(i\) denote a component of \(I\), and let \(j\) denote the component \(\ell\) of \(V\). We define the mixed second-level coordinates by
\[
    \hat{\mathbf I}^{(2),ij}_{s,t}
    :=
    \int_s^t I^i_{s,u}\,dV^\ell_u,
    \qquad
    \hat{\mathbf I}^{(2),ji}_{s,t}
    :=
    I^i_{s,t}V^\ell_{s,t}
    -
    \int_s^t I^i_{s,u}\,dV^\ell_u .
\]
With this definition, the shuffle identity at level two is satisfied by construction. Moreover, the required regularity follows from the It\^o isometry, while the second property follows from the additivity of the It\^o integral.

At the third level, the mixed objects are constructed as Young integrals of the second-level objects against the first level \(\hat I\). Indeed, whenever at least one of the coordinates corresponds to a Brownian component, the relevant second-level elements \(\hat{\mathbf I}^{(2),ij}\) and the corresponding first-level elements \(\hat I^k\) have complementary Young regularity. The shuffle identity at level three then follows from the integration-by-parts formula for Young integrals, together with the shuffle identity already established at the lower levels. Finally, the second property follows from the linearity and additivity of Young integration, while the required regularity follows from the Young--Loeve estimates.

The next Proposition collects this consructions and proves the existence of a Fernique estimate for $\hat{\mathbf I}$.
\begin{proposition}[Existence of the geometric joint lift and Fernique estimate]\label{existence_geometric_lift}
Let \(\hat I\) be the Gaussian process defined above. Then, for every
\(p\in(2\rho,4)\), there exists a Borel measurable map
\(
    \hat{\mathbf I}: \Omega \to 
    \Omega G_p\bigl([0,T],\mathbb R^{2d}\bigr)
\)
that is the enhancement of \(\hat I\), in the sense that
\(
    \text{proj}_1(\hat{\mathbf I}_{s,t})
    =
    \hat I_t-\hat I_s,
    \) for \(
    0\leq s\leq t\leq T.
\)
Moreover, there exists \(\eta>0\) such that
\[
    \mathbb E
    \left[
        \exp
        \left(
            \eta
            \|\hat{\mathbf I}\|_{p;[0,T]}^2
        \right)
    \right]
    <
    \infty .
\]
\end{proposition}
\begin{proof}
Set
\[
    \varsigma(s,t)
    :=
    \|C_I\|_{\rho\text{-var};[s,t]^2}
    +
    \|C_V\|_{\rho\text{-var};[s,t]^2}.
\]
Since \(V\) is Brownian, there exist constants
\(0<c\leq C<\infty\) such that, for every
\(0\leq s\leq t\leq T\),
\[
    c|t-s|
    \leq
    \|C_V\|_{\rho\text{-var};[s,t]^2}
    \leq
    C|t-s|.
\]
In what follows, \(M\) denotes a finite constant that may change from
line to line.

Writing
\(
    \hat I=(I^1,\ldots,I^d,V^1,\ldots,V^d),
\)
we have, for every \(j\in\{1,\ldots,2d\}\),
\[
    \mathbb E\left[
        \left|\hat I^j_{s,t}\right|^2
    \right]
    \leq
    M\varsigma(s,t).
\]
By Proposition 25 in \cite{friz2010differential}, it is enough to prove
that, for \(n=1,2,3\),
\[
    \left\|
        \operatorname{proj}_n
        \bigl(\log\widehat{\mathbf I}_{s,t}\bigr)
    \right\|_{L^2}
    \leq
    M\varsigma(s,t)^{n/2}.
\]
The pure \(I\)- and pure \(V\)-coordinates satisfy these estimates by
the known covariance \(\rho\)-variation bounds. It therefore remains
to consider the mixed coordinates.

At level two, the new terms are the mixed iterated integrals. For
\(i,k\in\{1,\ldots,d\}\), Itô's isometry gives
\[
\begin{aligned}
    \mathbb E\left[
        \left|
            \int_s^t I^i_{s,u}\,dV^k_u
        \right|^2
    \right]
    &=
    \int_s^t
    \mathbb E\left[
        \left|I^i_{s,u}\right|^2
    \right]du                                                        \\
    &\leq
    M|t-s|\,
    \|C_I\|_{\rho\text{-var};[s,t]^2}                                \\
    &\leq
    M\varsigma(s,t)^2.
\end{aligned}
\]
The reversed mixed integral is handled by integration by parts:
\[
    \int_s^t V^k_{s,u}\,dI^i_u
    =
    V^k_{s,t}I^i_{s,t}
    -
    \int_s^t I^i_{s,u}\,dV^k_u.
\]
By Hölder's inequality and Gaussian hypercontractivity,
\[
\begin{aligned}
    \left\|V^k_{s,t}I^i_{s,t}\right\|_{L^2}
    &\leq
    \left\|V^k_{s,t}\right\|_{L^4}
    \left\|I^i_{s,t}\right\|_{L^4}                                   \\
    &\leq
    M
    \left\|V^k_{s,t}\right\|_{L^2}
    \left\|I^i_{s,t}\right\|_{L^2}                                   \\
    &\leq
    M\varsigma(s,t).
\end{aligned}
\]
Thus all mixed level-two coordinates are bounded in \(L^2\) by
\(M\varsigma(s,t)\).

At level three, the logarithmic expansion and integration by parts
show that every mixed coordinate is a finite linear combination of
terms of the form
\[
    \hat I^{i_1}_{s,t}
    \hat I^{i_2}_{s,t}
    \hat I^{i_3}_{s,t},
    \qquad
    \hat I^{i_1}_{s,t}\mathcal A_{s,t},
    \qquad
    \int_s^t \mathcal A_{s,u}\,dV^k_u,
\]
where \(i_1,i_2,i_3\in\{1,\ldots,2d\}\),
\(k\in\{1,\ldots,d\}\), and \(\mathcal A\) denotes a scalar pure or
mixed level-two coordinate already constructed above. Terms with
final integration against a component of \(I\) can be reduced to
this form by integration by parts. For example,
\[
    \int_s^t
    \left(
        \int_s^u I^i_{s,r}\,dV^k_r
    \right)dI^j_u
    =
    I^j_{s,t}
    \int_s^t I^i_{s,r}\,dV^k_r
    -
    \int_s^t I^i_{s,u}I^j_{s,u}\,dV^k_u.
\]

By Gaussian hypercontractivity,
\[
    \left\|
        \hat I^{i_1}_{s,t}
        \hat I^{i_2}_{s,t}
        \hat I^{i_3}_{s,t}
    \right\|_{L^2}
    \leq
    M\varsigma(s,t)^{3/2}.
\]
Moreover, Hölder's inequality, Gaussian hypercontractivity, and the
level-two estimate give
\[
\begin{aligned}
    \left\|
        \hat I^{i_1}_{s,t}\mathcal A_{s,t}
    \right\|_{L^2}
    &\leq
    \left\|\hat I^{i_1}_{s,t}\right\|_{L^4}
    \left\|\mathcal A_{s,t}\right\|_{L^4}                             \\
    &\leq
    M\varsigma(s,t)^{1/2}\varsigma(s,t)                              \\
    &=
    M\varsigma(s,t)^{3/2}.
\end{aligned}
\]
Finally, by Itô's isometry, Gaussian hypercontractivity, and the
level-two estimate,
\[
\begin{aligned}
    \mathbb E\left[
        \left|
            \int_s^t \mathcal A_{s,u}\,dV^k_u
        \right|^2
    \right]
    &=
    \int_s^t
    \mathbb E\left[
        |\mathcal A_{s,u}|^2
    \right]du                                                        \\
    &\leq
    M\int_s^t \varsigma(s,u)^2\,du                                  \\
    &\leq
    M|t-s|\varsigma(s,t)^2                                          \\
    &\leq
    M\varsigma(s,t)^3.
\end{aligned}
\]
Hence all mixed level-three coordinates are bounded in \(L^2\) by
\(M\varsigma(s,t)^{3/2}\).

Combining the pure and mixed estimates gives, for \(n=1,2,3\),
\[
    \left\|
        \operatorname{proj}_n
        \bigl(\log\widehat{\mathbf I}_{s,t}\bigr)
    \right\|_{L^2}
    \leq
    M\varsigma(s,t)^{n/2}.
\]
Thanks to Condition~\ref{condition_filtration}, the logarithmic
coordinates belong to the corresponding homogeneous Wiener chaoses (see Appendix C in \cite{friz2010differential}).
Proposition 25 in \cite{friz2010differential} therefore applies to
the iterated integrals constructed above, giving a Borel measurable
geometric rough path
\[
    \widehat{\mathbf I}
    : \Omega \to
    \Omega G_p\bigl([0,T],\mathbb R^{2d}\bigr)
\]
for every \(p\in(2\rho,4)\). Moreover, for every such \(p\), there
exists \(\eta>0\) such that
\[
    \mathbb E\left[
        \exp\left(
            \eta
            \|\widehat{\mathbf I}\|_{p;[0,T]}^2
        \right)
    \right]
    <\infty.
\]
\end{proof}

\section{Proof of Lemma \ref{lem:holder_roughness_from_Kstar}}\label{app:proof_holder_roughness}
We collect here the proof of quantitative H\"older roughness of the driving Gaussian
process. Following the argument of \cite{cass2015smoothness}, we first use the
lower bound on conditional variances to obtain a small-ball estimate for its
directional increments. This estimate then yields almost-sure
H\"older roughness, together with finite moments of every order for
the inverse roughness modulus.
\begin{proof}[Proof of Lemma \ref{lem:holder_roughness_from_Kstar}]
If $I$ satisfies Condition \ref{condition_4} (A), then the statement is proved in Corollary 5.10 in \cite{cass2015smoothness}. We focus only on Condition \ref{condition_4} (B). 
For $i \in [1:d_I]$, let \(V^i\) be the standard Brownian motion associated to $I^i$. Then, as isonormal Gaussian processes over \(L^2([0,T], \mathbb R^{d_I})\), $V$ satisfies
\[
    I_t^i
    =
    V^i\bigl(K^*\mathbbm 1_{[0,t]}\bigr),
    \qquad
    i=1,\ldots,d_I.
\]
We first establish a small-ball estimate on a fixed interval. Fix
\(m\in\mathbb N\), and choose linearly independent functions
\(
    \psi_1,\ldots,\psi_m\in C_c^\infty((0,1)).
\)
For an interval \(\mathcal I=[a,a+\delta]\subseteq[0,T]\), define
\[
    \varphi_j^{\mathcal I}(u)
    :=
    \psi_j\left(\frac{u-a}{\delta}\right),
    \qquad
    u\in {\mathcal I}.
\]
For \(c=(c_1,\ldots,c_m)\in\mathbb R^{d_m}\), set
\[
    \varphi_c^{\mathcal I}:=\sum_{j=1}^m c_j\varphi_j^{\mathcal I}.
\]
Since the functions \(\psi_1,\ldots,\psi_m\) are linearly independent,
finite-dimensional norm equivalence yields constants \(a_m,b_m>0\),
depending also on \(T\), such that
\[
    \|\varphi_c^{\mathcal I}\|_{\infty;{\mathcal I}}
    \geq a_m|c|
\quad \text{and} \quad   \|\varphi_c^{\mathcal I}\|_{C^\gamma({\mathcal I})}
    \leq b_m\delta^{-\gamma}|c|.
\]
Applying Condition \ref{condition_4} (B) to \(\varphi_c^{\mathcal I}\) gives
\[
    a_m|c|
    \leq
    C\delta^{-\eta}
    \bigl\|K^*(\varphi_c^{\mathcal I}\mathbbm 1_{\mathcal I})\bigr\|_{L^2}^\vartheta
    \bigl(b_m\delta^{-\gamma}|c|\bigr)^{1-\vartheta}.
\]
Consequently, there exists \(c_m>0\) such that, for $\zeta = \eta/\vartheta +\gamma(1-\vartheta)/\vartheta$
\begin{equation}\label{coercivity_estimate_1}
    \bigl\|K^*(\varphi_c^I\mathbbm 1_{\mathcal I})\bigr\|_{L^2}
    \geq
    c_m\delta^\zeta|c|.
\end{equation}
Fix \(v \in\mathbb R^{d_I}\) with \(|v|=1\), and define the centred
Gaussian random variables
\[
    z_j^{{\mathcal I},v}
    :=
    \sum_{i =1}^{d_I}
    v_i
    V^i\bigl(K^*(\varphi_j^{\mathcal I}\mathbbm 1_{\mathcal I})\bigr),
    \qquad j=1,\ldots,m.
\]
Let \(\Sigma_{\mathcal I}\) denote the covariance matrix of
\(
    z^{{\mathcal I},v}=(z_1^{{\mathcal I},v},\ldots,z_m^{{\mathcal I},v}).
\)
For every \(c\in\mathbb R^{d_m}\), independence and identical distribution of
the components give
\[
\begin{aligned}
    c^\top\Sigma_{\mathcal I}c
    &=
    \sum_{i =1}^{d_I} v_i^2
    \bigl\|K^*(\varphi_c^{\mathcal I}\mathbbm 1_{\mathcal I})\bigr\|_2^2  =
    \bigl\|K^*(\varphi_c^{\mathcal I}\mathbbm 1_{\mathcal I})\bigr\|_2^2.
\end{aligned}
\]
It follows from \eqref{coercivity_estimate_1} that
\(
    c^\top\Sigma_{\mathcal I}c
    \geq c_m^2\delta^{2\zeta}|c|^2,
\)
and hence
\begin{equation}\label{determinant_bound}
    \det\Sigma_I
    \geq c_m^{2m}\delta^{2m\zeta}.
\end{equation}
By the definition of the isonormal process associated with \(I\),
\[
    z_j^{{\mathcal I},v}
    =
    \int_{\mathcal I}\varphi_j^{\mathcal I}(u)\,d\langle v,I_u\rangle.
\]
Since \(\varphi_j^{\mathcal I}\) vanishes at both endpoints of \(\mathcal I\), integration by parts
yields
\[
    z_j^{I,v}
    =
    -\int_{\mathcal I}
    \langle v,I_u-I_a\rangle
    \frac{d}{du}\varphi_j^{\mathcal I}(u)\,du.
\]
Moreover,
\[
    \int_{\mathcal I } \left|\frac{d}{du}\varphi_j^I(u)\right|\,du
    =
    \int_0^1 \left|\frac{d}{dr}\psi_j(r)\right| \,dr.
\]
Therefore,
\begin{equation} \label{oscillation_bound}
    |z_j^{\mathcal I,v}|
    \leq
    C_m
    \operatorname{osc}_{\mathcal I}(\langle v ,I\rangle),
\end{equation}
where
\(
    \operatorname{osc}_{\mathcal I}(\langle v,I\rangle)
    :=
    \sup_{u,w\in \mathcal I}
    |\langle v,I_u-I_w\rangle|.
\)

The density of \(z^{\mathcal I, v}\) is bounded from above by
\[
    (2\pi)^{-m/2}(\det\Sigma_{\mathcal I})^{-1/2}.
\]
Combining this with \eqref{determinant_bound} and \eqref{oscillation_bound}, we obtain
\begin{equation}\label{probability_oscillation}
    \mathbb P\left(
        \operatorname{osc}_{\mathcal I}(\langle v,I\rangle)\leq r
    \right)
    \leq
    C_m
    \left(\frac{r}{\delta^\zeta}\right)^m.
\end{equation}
This estimate is uniform over all intervals \(\mathcal I\) and all
\(v\in S^{d_I-1}\), here $S^{d_I-1}$ is the unit sphere in $\mathbb R^{d_I}$.

We now pass to all intervals and all directions. Let
\(
    \delta_n:=T2^{-n},
\)
and let \(\mathcal P_n\) denote the dyadic partition of \([0,T]\) into
intervals of length \(\delta_n\). For \(\iota >0\), define the event
\[
    \mathcal B_n(\iota)
    :=
    \left\{
        \begin{array}{c}
        \text{there exist }\mathcal I\in\mathcal P_n
        \text{ and } v\in S^{d_I-1}\text{ such that}\\[2mm]
        \operatorname{osc}_{\mathcal I}(\langle v ,I\rangle)
        \leq \iota \delta_n^\theta
        \end{array}
    \right\},
\]
with $\theta \in (\zeta, 2/p \wedge 1)$.
For \(N\geq1\), set
\(
    A_N:=\{\|I\|_{\infty;[0,T]}\leq N\}.
\)

Fix \(n\) and let \(\mathcal N_{n,N,\iota}\) be a
\(
    \iota \delta_n^\theta/(4N)
\)-net of \(S^{d_I-1}\). Its cardinality satisfies
\begin{equation}\label{boundN_n_iota}
    \# \mathcal N_{n,N,\iota}
    \leq
    C_d
    \left(\frac{N}{\iota \delta_n^\theta}\right)^{d_I-1},
    \end{equation}
for all sufficiently small \(\iota\).

Suppose that \(\mathcal B_n(\iota)\cap A_N\) occurs, and let \(v\) be a direction
appearing in the definition of \(\mathcal B_n(\iota)\). Choose
\(\widehat v\in\mathcal N_{n,N,\iota}\) such that
\(
    | v-\widehat v|
    \leq {\iota \delta_n^\theta} / {4N}.
\)
On \(A_N\),
\[
\begin{aligned}
    \operatorname{osc}_{\mathcal I}(\langle\widehat v,I\rangle)
    &\leq
    \operatorname{osc}_{\mathcal I }(\langle v,I\rangle)
    +
    2N|v-\widehat v| \leq
    \frac{3}{2}\iota \delta_n^\theta.
\end{aligned}
\]
Using \eqref{probability_oscillation}, \eqref{boundN_n_iota}, and the fact that
\(
    |\mathcal P_n| =T/\delta_n,
\)
we obtain
\begin{equation}\label{probability_event_A}
\begin{aligned}
    \mathbb{P}(\mathcal B_n(\iota)\cap A_N)
    &\leq
    C_m
    \frac{T}{\delta_n}
    \left(\frac{N}{\iota \delta_n^\theta}\right)^{d_I-1}
    \iota^m\delta_n^{m(\theta-\zeta)} \\
    &\leq
    C_m
    N^{d_I-1}\iota ^{m-d_I+1}\delta_n^\beta,
\end{aligned}
\end{equation}
where
\(
    \beta
    :=
    m(\theta-\zeta)-(d_I-1)\theta-1.
\)
Since \(\theta> \zeta\), we may choose \(m\) sufficiently large that
\(\beta>0\). Hence
\[
    \sum_{n=0}^\infty\delta_n^\beta<\infty,
\]
and \eqref{probability_event_A} gives
\begin{equation}\label{probability_bound_B_A_2}
    \mathbb P\left(
        \bigcup_{n\geq0}\mathcal B_n(\iota)\cap A_N
    \right)
    \leq
    C_mN^{d_I-1}\iota^{m-d_I+1}.
\end{equation}

It remains to relate the events \(\mathcal B_n(\iota)\) to \(L_\theta(I)\). Suppose that
\(L_\theta(I)\leq\lambda\). By the definition of the infimum, there exist
\(s\in[0,T]\), \(0<\epsilon\leq T/2\), and
\(v\in S^{d_I-1}\) such that
\begin{equation}\label{sup_bound_v_X}
    \sup_{\epsilon/2<|t-s|<\epsilon}
    |\langle v,I_t-I_s\rangle|
    <2\lambda\epsilon^\theta.
\end{equation}
If \(s\leq T/2\), consider
\(
    \tilde{\mathcal I}=(s+ \epsilon/2,s+\epsilon),
\)
if \(s>T/2\), consider instead
\(
    \tilde{\mathcal I}=(s-\epsilon,s-\epsilon/2).
\)
In either case, \(\tilde{\mathcal I}\subseteq[0,T]\) has length
\(\epsilon/2\) and is contained in the annulus appearing in
\eqref{sup_bound_v_X}.

There exists a dyadic interval \(\mathcal I\subset \tilde{\mathcal{I}}\), of length \(\delta_n\), such
that
\begin{equation}\label{length_interval}
    \frac{ \epsilon}{16}
    <
    \delta_n
    \leq
    \frac{ \epsilon}{8}.
\end{equation}
For \(u,w\in \mathcal I\), \eqref{sup_bound_v_X}, \eqref{length_interval} imply
\[
\begin{aligned}
    |\langle v,I_u-I_w\rangle|
    &\leq
    |\langle v,I_u-I_s\rangle|
    +
    |\langle v,I_w-I_s\rangle| \\
    &<
    4\lambda\epsilon^\theta \\
    &\leq
    4\cdot16^\theta\lambda\delta_n^\theta.
\end{aligned}
\]
Consequently,
\begin{equation}\label{L_theta_inclusion}
    \{L_\theta(I)\leq\lambda\}
    \subseteq
    \bigcup_{n\geq0}\mathcal B_n(C_\theta\lambda),
    \qquad
    C_\theta:=4\cdot16^\theta.
\end{equation}

Combining \eqref{probability_bound_B_A_2} and \eqref{L_theta_inclusion}, we find that, for every sufficiently
large \(m\),
\begin{equation}\label{bound_roughness_exponent}
    \mathbb P\bigl(L_\theta(I)\leq\lambda,\ A_N\bigr)
    \leq
    C_mN^{d-1}\lambda^{m-d+1}.
\end{equation}
Take \(N=\lambda^{-1}\). Since \(I\) is a continuous Gaussian process,
Fernique's theorem gives constants \(c,C>0\) such that
\(
    \mathbb P(A_N^c)
    \leq
    C\exp(-cN^2).
\)
It follows from \eqref{bound_roughness_exponent} that
\[
    \mathbb P\bigl(L_\theta(I)\leq\lambda\bigr)
    \leq
    C_m\lambda^{m-2d_I+2}
    +
    C\exp(-c\lambda^{-2}).
\]
Since \(m\) can be chosen arbitrarily large, for every \(\iota >0\) there exists
\(C_\iota >0\) such that
\begin{equation}\label{estimate_probability_final}
    \mathbb P\bigl(L_\theta(I)\leq\lambda\bigr)
    \leq
    C_\iota\lambda^\iota,
    \qquad 0<\lambda\leq1.
\end{equation}
In particular,
\(
    \mathbb P(L_\theta(I)=0)=0,
\)
so \(I\) is almost surely \(\theta\)-Hölder rough.

Finally, for \(l>0\), choose \(\iota >l\). We can now write,
\[
\begin{aligned}
    \mathbb E[L_\theta(I)^{-l}]
    &=
    l\int_0^\infty
    r^{l-1}
    \mathbb P\bigl(L_\theta(I)^{-1}>r\bigr)\,dr \\
    &=
    l\int_0^\infty
    r^{l-1}
    \mathbb P\bigl(L_\theta(I)<r^{-1}\bigr)\,dr.
\end{aligned}
\]
The integral over \(r\in(0,1)\) is finite, while, by \eqref{estimate_probability_final},
\[
    \int_1^\infty
    r^{l-1}
    \mathbb P\bigl(L_\theta(I)<r^{-1}\bigr)\,dr
    \leq
    C_\iota \int_1^\infty r^{l-\iota-1}\,dr
    <\infty.
\]
Therefore
\(
    \mathbb E[L_\theta(I)^{-l}]<\infty
\)
for every \(l>0\).
\end{proof}
\bibliographystyle{amsplain}
\bibliography{biblio}

@article{decreusefond1999stochastic,
  author  = {Decreusefond, Laurent and {\"U}st{\"u}nel, Ali S.},
  title   = {Stochastic Analysis of the Fractional Brownian Motion},
  journal = {Potential Analysis},
  volume  = {10},
  number  = {2},
  pages   = {177--214},
  year    = {1999},
  doi     = {10.1023/A:1008634027843}
}

@article{friz2010differential,
  author  = {Friz, Peter K. and Victoir, Nicolas},
  title   = {Differential Equations Driven by Gaussian Signals},
  journal = {Annales de l'Institut Henri Poincar{\'e}, Probabilit{\'e}s et Statistiques},
  volume  = {46},
  number  = {2},
  pages   = {369--413},
  year    = {2010},
  doi     = {10.1214/09-AIHP202}
}

@book{bain2009fundamentals,
  author    = {Bain, Alan and Crisan, Dan},
  title     = {Fundamentals of Stochastic Filtering},
  series    = {Stochastic Modelling and Applied Probability},
  volume    = {60},
  edition   = {1},
  publisher = {Springer},
  address   = {New York, NY},
  year      = {2009},
  doi       = {10.1007/978-0-387-76896-0},
  isbn      = {978-0-387-76895-3}
}

@article{cass2013integrability,
  author  = {Cass, Thomas and Litterer, Christian and Lyons, Terry},
  title   = {Integrability and Tail Estimates for Gaussian Rough Differential Equations},
  journal = {The Annals of Probability},
  volume  = {41},
  number  = {4},
  pages   = {3026--3050},
  year    = {2013},
  doi     = {10.1214/12-AOP821}
}

@article{crisan2013robust,
  author  = {Crisan, Dan and Diehl, Joscha and Friz, Peter K. and Oberhauser, Harald},
  title   = {Robust Filtering: Correlated Noise and Multidimensional Observation},
  journal = {The Annals of Applied Probability},
  volume  = {23},
  number  = {5},
  pages   = {2139--2160},
  year    = {2013},
  doi     = {10.1214/12-AAP896}
}

@article{friz2016jain,
  author  = {Friz, Peter K. and Gess, Benjamin and Gulisashvili, Archil and Riedel, Sebastian},
  title   = {The Jain--Monrad Criterion for Rough Paths and Applications to Random Fourier Series and Non-Markovian H{\"o}rmander Theory},
  journal = {The Annals of Probability},
  volume  = {44},
  number  = {1},
  pages   = {684--738},
  year    = {2016},
  doi     = {10.1214/14-AOP986}
}

@book{kilbas1993fractional,
  author    = {Samko, Stefan G. and Kilbas, Anatoly A. and Marichev, Oleg I.},
  title     = {Fractional Integrals and Derivatives: Theory and Applications},
  publisher = {Gordon and Breach Science Publishers},
  address   = {Yverdon},
  year      = {1993},
  isbn      = {978-2-88124-864-1}
}

@article{cass2022combinatorial,
  author  = {Cass, Thomas and Driver, Bruce K. and Litterer, Christian and Rossi Ferrucci, Emilio},
  title   = {A Combinatorial Approach to Geometric Rough Paths and Their Controlled Paths},
  journal = {Journal of the London Mathematical Society},
  volume  = {106},
  number  = {2},
  pages   = {936--981},
  year    = {2022},
  doi     = {10.1112/jlms.12589}
}

@article{deya2009rough,
  title={Rough Volterra equations 1: the algebraic integration setting},
  author={Deya, Aur{\'e}lien and Tindel, Samy},
  journal={Stochastics and Dynamics},
  volume={9},
  number={03},
  pages={437--477},
  year={2009},
  publisher={World Scientific}
}

@article{cass2015smoothness,
 ISSN = {00911798},
 URL = {http://www.jstor.org/stable/24519169},
 author = {Thomas Cass and Martin Hairer and Christian Litterer and Samy Tindel},
 journal = {The Annals of Probability},
 number = {1},
 pages = {188--239},
 publisher = {Institute of Mathematical Statistics},
 title = {SMOOTHNESS OF THE DENSITY FOR SOLUTIONS TO GAUSSIAN ROUGH DIFFERENTIAL EQUATIONS},
 urldate = {2026-07-30},
 volume = {43},
 year = {2015}
}

@article{Jin2015,
 ISSN = {00255718, 10886842},
 URL = {http://www.jstor.org/stable/45106636},
 author = {BANGTI JIN and RAYTCHO LAZAROV and JOSEPH PASCIAK and WILLIAM RUNDELL},
 journal = {Mathematics of Computation},
 number = {296},
 pages = {2665--2700},
 publisher = {American Mathematical Society},
 title = {VARIATIONAL FORMULATION OF PROBLEMS INVOLVING FRACTIONAL ORDER DIFFERENTIAL OPERATORS},
 urldate = {2026-08-14},
 volume = {84},
 year = {2015}
}

@book{friz2010multidimensional,
  author    = {Friz, Peter K. and Victoir, Nicolas B.},
  title     = {Multidimensional Stochastic Processes as Rough Paths: Theory and Applications},
  series    = {Cambridge Studies in Advanced Mathematics},
  volume    = {120},
  publisher = {Cambridge University Press},
  address   = {Cambridge},
  year      = {2010},
  isbn      = {978-0-521-87607-0}
}

@book{nualart2009malliavin,
  author    = {Nualart, David},
  title     = {Malliavin Calculus and Its Applications},
  series    = {CBMS Regional Conference Series in Mathematics},
  volume    = {110},
  publisher = {American Mathematical Society},
  address   = {Providence, RI},
  year      = {2009},
  doi       = {10.1090/cbms/110}
}

@book{nourdin2012normal,
  author    = {Nourdin, Ivan and Peccati, Giovanni},
  title     = {Normal Approximations with Malliavin Calculus: From Stein's Method to Universality},
  series    = {Cambridge Tracts in Mathematics},
  volume    = {192},
  publisher = {Cambridge University Press},
  address   = {Cambridge},
  year      = {2012},
  isbn      = {978-1-107-01777-1}
}

@article{inahama2014malliavin,
  title={Malliavin differentiability of solutions of rough differential equations},
  author={Inahama, Yuzuru},
  journal={Journal of Functional Analysis},
  volume={267},
  number={5},
  pages={1566--1584},
  year={2014},
  publisher={Elsevier}
}

@incollection{nualart1989partial,
  author    = {Nualart, David and Zakai, Moshe},
  title     = {The Partial Malliavin Calculus},
  booktitle = {S{\'e}minaire de Probabilit{\'e}s XXIII},
  editor    = {Az{\'e}ma, Jacques and Yor, Marc and Meyer, Paul-Andr{\'e}},
  series    = {Lecture Notes in Mathematics},
  volume    = {1372},
  pages     = {362--381},
  publisher = {Springer},
  address   = {Berlin, Heidelberg},
  year      = {1989},
  doi       = {10.1007/BFb0083986},
  url       = {https://www.numdam.org/item/SPS_1989__23__362_0/}
}

@article{friz2018differential,
  title={Differential equations driven by rough paths with jumps},
  author={Friz, Peter K and Zhang, Huilin},
  journal={Journal of Differential Equations},
  volume={264},
  number={10},
  pages={6226--6301},
  year={2018},
  publisher={Elsevier}
}

@book{karatzas1991brownian,
  author    = {Karatzas, Ioannis and Shreve, Steven E.},
  title     = {Brownian Motion and Stochastic Calculus},
  series    = {Graduate Texts in Mathematics},
  volume    = {113},
  edition   = {2},
  publisher = {Springer},
  address   = {New York, NY},
  year      = {1991},
  doi       = {10.1007/978-1-4612-0949-2}
}

@article{friz2010rough,
  author  = {Friz, Peter K. and Oberhauser, Harald},
  title   = {Rough Path Stability of (Semi-)Linear {SPDEs}},
  journal = {Probability Theory and Related Fields},
  volume  = {158},
  number  = {1--2},
  pages   = {401--434},
  year    = {2014},
  doi     = {10.1007/s00440-013-0483-2}
}

@article{crandall1992users,
  author  = {Crandall, Michael G. and Ishii, Hitoshi and Lions, Pierre-Louis},
  title   = {User's Guide to Viscosity Solutions of Second Order Partial Differential Equations},
  journal = {Bulletin of the American Mathematical Society},
  volume  = {27},
  number  = {1},
  pages   = {1--67},
  year    = {1992},
  doi     = {10.1090/S0273-0979-1992-00266-5}
}

@article{crisan2020high,
  title={A high order time discretization of the solution of the non-linear filtering problem},
  author={Crisan, Dan and Ortiz-Latorre, Salvador},
  journal={Stochastics and Partial Differential Equations: Analysis and Computations},
  volume={8},
  number={4},
  pages={693--760},
  year={2020},
  publisher={Springer}
}

@article{alos2001stochastic,
  title={Stochastic calculus with respect to Gaussian processes},
  author={Alos, Elisa and Mazet, Olivier and Nualart, David},
  journal={The Annals of Probability},
  volume={29},
  number={2},
  pages={766--801},
  year={2001},
  publisher={Institute of Mathematical Statistics}
}

@article{bouleau1986proprietes,
  title={Propri{\'e}t{\'e}s d'absolue continuit{\'e} dans les espaces de Dirichlet et applications aux {\'e}quations diff{\'e}rentielles stochastiques},
  author={Bouleau, Nicolas and Hirsch, Francis},
  journal={S{\'e}minaire de probabilit{\'e}s de Strasbourg},
  volume={20},
  pages={131--161},
  year={1986}
}

@misc{decreusefond2003stochasticintegrationrespectvolterra,
      title={Stochastic Integration with respect to Volterra processes}, 
      author={L. Decreusefond},
      year={2003},
      eprint={math/0302047},
      archivePrefix={arXiv},
      primaryClass={math.PR},
      url={https://arxiv.org/abs/math/0302047}, 
}

@article{coutin1999abstract,
  title={Abstract nonlinear filtering theory in the presence of fractional Brownian motion},
  author={Coutin, Laure and Decreusefond, Laurent},
  journal={The Annals of Applied Probability},
  volume={9},
  number={4},
  pages={1058--1090},
  year={1999},
  publisher={Institute of Mathematical Statistics}
}

@article{kleptsyna1999linear,
  title={Linear filtering with fractional Brownian motion in the signal and observation processes},
  author={Kleptsyna, Marina L and Kloeden, Peter E and Anh, Vo Van},
  journal={International Journal of Stochastic Analysis},
  volume={12},
  number={1},
  pages={85--90},
  year={1999},
  publisher={Wiley Online Library}
}

@article{kleptsyna1998linear,
  title={Linear filtering with fractional Brownian motion},
  author={Kleptsyna, ML and Kloeden, PE and Anh, VV},
  journal={Stochastic analysis and applications},
  volume={16},
  number={5},
  pages={907--914},
  year={1998},
  publisher={Taylor \& Francis}
}

@article{ml2000general,
  title={General approach to filtering with fractional Brownian noises—application to linear systems},
  author={Kleptsyna, M L and Le Breton, A and Roubaud, MC},
  journal={Stochastics: An International Journal of Probability and Stochastic Processes},
  volume={71},
  number={1-2},
  pages={119--140},
  year={2000},
  publisher={Taylor \& Francis}
}

@article{amirdjanova2006new,
  title={New method for optimal nonlinear filtering of noisy observations by multiple stochastic fractional integral expansions},
  author={Amirdjanova, Anna and Chivoret, S},
  journal={Computers \& Mathematics with Applications},
  volume={52},
  number={1-2},
  pages={161--178},
  year={2006},
  publisher={Elsevier}
}

@article{ocone1983multiple,
  title={Multiple integral expansions for nonlinear filtering},
  author={Ocone, Daniel},
  journal={Stochastics: An International Journal of Probability and Stochastic Processes},
  volume={10},
  number={1},
  pages={1--30},
  year={1983},
  publisher={Taylor \& Francis}
}

@book{lyons2007differential,
  author    = {Lyons, Terry J. and Caruana, Michael and L{\'e}vy, Thierry},
  title     = {Differential Equations Driven by Rough Paths},
  series    = {Lecture Notes in Mathematics},
  volume    = {1908},
  publisher = {Springer},
  address   = {Berlin, Heidelberg},
  year      = {2007},
  doi       = {10.1007/978-3-540-71285-5}
}

@article{GUBINELLI2010693,
  author  = {Gubinelli, Massimiliano},
  title   = {Ramification of Rough Paths},
  journal = {Journal of Differential Equations},
  volume  = {248},
  number  = {4},
  pages   = {693--721},
  year    = {2010},
  doi     = {10.1016/j.jde.2009.11.015}
}

@article{barles2006new,
  title={A new stability result for viscosity solutions of nonlinear parabolic equations with weak convergence in time},
  author={Barles, Guy},
  journal={Comptes Rendus Mathematique},
  volume={343},
  number={3},
  pages={173--178},
  year={2006},
  publisher={Elsevier}
}

@article{hairer2013regularity,
  author  = {Hairer, Martin and Pillai, Natesh S.},
  title   = {Regularity of laws and ergodicity of hypoelliptic {SDEs} driven by rough paths},
  journal = {The Annals of Probability},
  volume  = {41},
  number  = {4},
  pages   = {2544--2598},
  year    = {2013},
  mrnumber = {3112925}
}

@incollection{clark1978design,
  author    = {Clark, John M. C.},
  title     = {The Design of Robust Approximations to the Stochastic Differential Equations of Nonlinear Filtering},
  editor    = {Skwirzynski, J. K.},
  booktitle = {Communication Systems and Random Process Theory},
  series    = {NATO Advanced Study Institute Series E: Applied Sciences},
  volume    = {25},
  pages     = {721--734},
  publisher = {Sijthoff \& Noordhoff},
  address   = {Alphen aan den Rijn},
  year      = {1978}
}

@article{clark2005robust,
  title={On a robust version of the integral representation formula of nonlinear filtering},
  author={Clark, John MC and Crisan, Dan},
  journal={Probability theory and related fields},
  volume={133},
  number={1},
  pages={43--56},
  year={2005},
  publisher={Springer}
}

@Inbook{Gawarecki2001,
author="Gawarecki, L.
and Mandrekar, V.",
title="On the Zakai Equation of Filtering with Gaussian Noise",
bookTitle="Stochastics in Finite and Infinite Dimensions: In Honor of Gopinath Kallianpur",
year="2001",
publisher="Birkh{\"a}user Boston",
address="Boston, MA",
pages="145--151",
isbn="978-1-4612-0167-0",
doi="10.1007/978-1-4612-0167-0_8",
url="https://doi.org/10.1007/978-1-4612-0167-0_8"
}

@article{Livieri02092018,
  author  = {Livieri, Giulia and Mouti, Saad and Pallavicini, Andrea and Rosenbaum, Mathieu},
  title   = {Rough Volatility: Evidence from Option Prices},
  journal = {IISE Transactions},
  volume  = {50},
  number  = {9},
  pages   = {767--776},
  year    = {2018},
  doi     = {10.1080/24725854.2018.1444297}
}

@misc{floch2022roughnessimpliedvolatility,
      title={Roughness of the Implied Volatility}, 
      author={Fabien Le Floc'h},
      year={2022},
      eprint={2207.04930},
      archivePrefix={arXiv},
      primaryClass={q-fin.MF},
      url={https://arxiv.org/abs/2207.04930}, 
}

\end{document}